\documentclass[12pt,leqno]{article}
\usepackage{amsthm,amsfonts,amssymb,amsmath,eufrak,oldgerm}
\usepackage{epsfig}
\usepackage{fontenc}
\renewcommand\d{\partial}
\renewcommand\a{\alpha}
\renewcommand\b{\beta}
\renewcommand\o{\omega}
\newcommand\R{\mathbb R}

\newcommand\D{{${\mathcal D}_\eps$}}
\def\g{\gamma}

\def\txi{{\tilde \xi}}
\def\ess{{\rm ess}}

\def\eps{\varepsilon }
\def\e{\varepsilon} 
\newcommand\kernel{\hbox{\rm Ker}}

\newcommand\errfn{\textrm{errfn}}
\newcommand\br{\begin{rem}}
\newcommand\er{\end{rem}}
\newcommand\bp{\begin{pmatrix}}
\newcommand\ep{\end{pmatrix}}
\newcommand\be{\begin{equation}}
\newcommand\ee{\end{equation}}
\newcommand\ba{\begin{equation}\begin{aligned}}
\newcommand\ea{\end{aligned}\end{equation}}

\newcommand{\CalB}{\mathcal{B}}
\newcommand{\CalE}{\mathcal{E}}
\newcommand{\CalF}{\mathcal{F}}
\newcommand{\CalC}{\mathcal{C}}

\newcommand{\CalJ}{\mathcal{J}}
\newcommand{\CalK}{\mathcal{K}}
\newcommand{\CalL}{\mathcal{L}}
\newcommand{\CalN}{\mathcal{N}}
\newcommand{\CalO}{\mathcal{O}}
\newcommand{\CalR}{\mathcal{R}}
\newcommand{\CalS}{\mathcal{S}}
\newcommand{\CalT}{\mathcal{T}}
\newcommand{\CalV}{\mathcal{V}}

\newcommand{\RR}{{\mathbb R}}

\newcommand{\ZZ}{{\mathbb Z}}
\newcommand{\TT}{{\mathbb T}}

\newcommand{\Id}{{\rm Id }}
\newcommand{\Range}{{\rm Range }}
\newcommand{\Span}{{\rm Span }}

\newcommand{\Res}{{\rm Residue}}

\newtheorem{theo}{Theorem}[section]
\newtheorem{prop}[theo]{Proposition}
\newtheorem{cor}[theo]{Corollary}
\newtheorem{lem}[theo]{Lemma}

\newtheorem{ass}[theo]{Assumption}

\newtheorem{rem}[theo]{Remark}
\newtheorem{rems}[theo]{Remarks}

\numberwithin{equation}{section}

\title{
Galloping instability of viscous shock waves
}

\author{\sc \small Benjamin Texier\thanks{
Universit\'e Paris 7/Denis Diderot, Institut de Math\'ematiques de Jussieu, UMR CNRS 7586;
texier@math.jussieu.fr:
Research of B.T.  was partially supported
under NSF grant number DMS-0505780.
} and
Kevin Zumbrun\thanks{Indiana University, Bloomington, IN 47405;
kzumbrun@indiana.edu:
Research of K.Z. was partially supported
under NSF grant number DMS-0300487.
}}

\begin{document}

\maketitle

\begin{abstract}
Motivated by physical and numerical observations of time oscillatory
``galloping'', ``spinning'', and
``cellular'' instabilities of detonation waves,
we study Poincar\'e--Hopf bifurcation of traveling-wave
solutions of viscous conservation laws.
The main difficulty is the absence of a spectral
gap between oscillatory modes and essential spectrum, preventing
standard reduction to a finite-dimensional center manifold.
We overcome this by direct Lyapunov--Schmidt reduction,
using detailed pointwise bounds on the linearized solution operator
to carry out a nonstandard implicit function construction in
the absence of a spectral gap.
The key computation is a space-time stability
estimate on the transverse linearized solution operator reminiscent of 
Duhamel estimates carried out on the full solution operator 
in the study of nonlinear stability of spectrally stable traveling waves.
\end{abstract}

\tableofcontents

\bigbreak
\section{Introduction}

Motivated by physical and numerical observations of 
time-oscillatory ``galloping''
or ``pulsating'' instabilities of detonation waves
\cite{MT, BMR, FW, MT, ALT, AT, F1, F2},
we study the Poincar\'e--Hopf bifurcation of viscous 
shock waves in one spatial dimension.
Our main result is to obtain
a rigorous (nonlinear) characterization in terms of 
spectral information.
The complementary problem of verifying this spectral information has been 
studied already in \cite{LyZ1, LyZ2}, where it was shown that transition 
from stability to instability of viscous detonation
waves generically involves a Poincar\'e--Hopf type bifurcation in the 
spectral configuration of the linearized operator about the wave.
By essentially the same analysis we obtain also a corresponding 
multi-dimensional result applying to planar viscous shock 
fronts traveling in a cylinder of finite cross-section,
with artificial Neumann or periodic boundary conditions. 
This gives a simplified mathematical model for time-oscillatory instabilities
observed in detonation waves moving within a duct, which, besides
the longitudinal galloping instabilities described above, 
include also transverse ``cellular'' or ``spinning'' instabilities.
The method of analysis appears to be of general application,
in particular, 
with suitable elaboration, 
to extend to the 
originally motivating case of viscous detonation waves of the 
reactive Navier--Stokes equations with physical viscosity.

{}From a mathematical standpoint, the main issue in our
analysis is that the linearized
equations about a standing shock wave have no spectral gap
between convective modes corresponding to essential spectra
of the linearized operator $L$ about the wave and oscillatory
modes corresponding to pure imaginary point spectra.
This prohibits the usual PDE analysis by center manifold reduction to
a finite-dimensional subspace; likewise, at the linearized level, decay
to the center subspace is at time-algebraic rather than
time-exponential rate, so that oscillatory and other modes are
strongly coupled.
We overcome this difficulty by the introduction of a nonstandard
Implicit Function Theorem framework suitable for infinite-dimensional
Lyapunov--Schmidt reduction in the absence of spectral gap
(Section \ref{LS}),
augmented by a finite-dimensional ``weak'' Implicit Function Theorem 
suitable for finite-dimensional bifurcation analysis in situations
of limited regularity (Section \ref{Brouwer}).
The latter, based on the Brouwer Fixed-point Theorem rather than the standard
Contraction-mapping construction, seems of interest in its own right.

The study of bifurcation from stability is a natural followup
to our previous work on stability of viscous shock and 
detonation waves; see, e.g., 
\cite{ZH, MaZ3, MaZ4, HRZ, LyZ1, LyZ2, LRTZ}.
Interestingly, the key estimate needed to apply our bifurcation
framework turns out to be a space-time stability
estimate on the transverse linearized solution operator (meaning the
part complementary to oscillatory modes) quite similar to estimates 
on the full solution operator arising in
the stability analysis of viscous traveling waves. 
To carry out this estimate requires rather detailed pointwise information on
the Green function of the linearized operator about the wave;
see Sections \ref{refined} and \ref{1dproof}.
Indeed, we use the full power of the pointwise semigroup techniques
developed in \cite{ZH, MaZ1, MaZ3}.

\subsection{Shocks, detonations, and galloping}\label{sdg}

The equations of compressible gas dynamics in one spatial dimension,
like many equations in continuum mechanics, 
take the form of hyperbolic conservation laws
\be\label{Euler}
u_t+f(u)_x=0
\ee
(Euler equations)
or hyperbolic--parabolic conservation laws
\be\label{NS}
u_t+f(u)_x= (B(u)u_x)_x
\ee
(Navier--Stokes equations), depending whether second-order
transport effects-- in this case, viscosity and heat-conduction--
are neglected or included.
Here, $x,t\in \RR^1$ are spatial location and time, 
$u\in \RR^n$ is a vector of densities of 
conserved quantities--
mass, momentum, and energy in the case of gas dynamics, 
$f\in \RR^n$ is a vector of corresponding fluxes,
and $B\in \RR^{n\times n}$ is a matrix of
transport coefficients.

Such equations are well known to support traveling-wave solutions of
form
\be\label{shock}
u(x,t)=\bar u(x-st):=\begin{cases} u_- & x-st\le 0,\\
u_+ & x-st> 0\\
\end{cases}
\ee
(discontinuous) and
\be\label{profile}
u(x,t)=\bar u(x-st), \qquad
\lim_{z\to \pm \infty}\bar u(z)=u_\pm
\ee
(smooth), respectively, known as ideal and viscous shock waves.
These waves may be observed physically; indeed, the corresponding
physical objects appear to be quite stable \cite{BE}.

Similarly, the equations of reacting, compressible gas dynamics
take the form
\ba\label{rEuler}
u_t+\tilde f(u,z)_x -qK\phi(u)z &=0\\
z_t
+(v(u,z)z)_x 
+ K\phi(u)z &= 0
\ea
(reactive Euler, or Zeldovich--von Neumann--Doering equations (ZND)) or
\ba\label{rNS}
u_t+\tilde f(u,z)_x -qK\phi(u)z &=(\tilde B(u,z)u_x)_x,\\
z_t
+(v(u,z)z)_x 
+ K\phi(u)z &= (\tilde D(u,z)z_x)_x
\ea
(reactive Navier--Stokes equations (rNS)), where $x,t\in \RR^1$, 
$u\in \RR^n$ is as before, $z\in \RR^r$ is a vector of
mass fractions of different reactant species, $v\in \RR^1$ is
fluid velocity, and $\tilde f$,
$\tilde B$ and $\tilde D$ are flux vectors and transport matrices
depending (through $z$) on chemical makeup of the gas, with
$\tilde f(u,0)=f(u)$, $\tilde B(u,0)=B(u)$.
The matrix $K$ models reaction dynamics, 
$q$ is a constant heat release coefficient
($q > 0$ corresponding to an exothermic reaction),
and $\phi\in \RR^1$ is an ignition function serving to 
``turn on'' the reaction: zero for temperatures
below a certain critical temperature and positive
for temperatures above.

Equations \eqref{rEuler} and \eqref{rNS} support traveling
wave solutions 
\be\label{detprofile}
u(x,t)=\bar u(x-st)
\ee
analogous to \eqref{shock} and \eqref{profile}, known
as ideal and viscous detonation waves, which may likewise
be observed physically.
However, in contrast to the shock wave case, detonations appear
to be rather unstable.
What is typically observed is not the planar, steadily progressing
solution \eqref{detprofile}, but rather a nearby solution varying 
time-periodically about \eqref{detprofile}: that is, an apparent
{\it bifurcation with exchange of stability}.
In the typical experimental setting of a detonation moving along
a duct, these may be longitudinal ``galloping'', or ``pulsating'',
instabilities for which the planar structure is maintained, but
the form of the profile changes time-periodically, or they may
be transverse instabilities
in which the front progresses steadily in the longitudinal (i.e., axial)
direction, but develops time-oscillatory
structure in the transverse (cross-sectional) directions of a form
depending on the geometry of the cross-section:
``cellular'' instabilities for a polygonal (e.g. rectangular) cross-section;
``spinning'' instabilities for a circular cross-section,
with one or more ``hot spots'', or ``combustion heads'',
moving spirally along the duct.

Stability of shocks and detonations may be studied within a unified 
mathematical framework;  see, e.g., 
\cite{Er1, Er2, Ko1, Ko2, D, BT, FD, LS, T, K, M1, M2, M3, FM, MeKochel, CJLW} 
in the inviscid case \eqref{Euler}, \eqref{rEuler}, and
\cite{Sat, Go1, Go2, KM, KMN, MN,
Liu86, Liu97, GX, SX, GZ, ZH, Br1, Br2, BrZ, BDG, KK, Z1, Z2, Z3, Z4,
MaZ2, MaZ3, MaZ4, GMWZ1, GMWZ2, GMWZ3, HZ, BL, LyZ1, LyZ2, JLW, HuZ1, HuZ2, 
BHRZ, PZ, FS, LRTZ}
in the viscous case \eqref{NS}, \eqref{rNS}.
In particular, stability in each case has been shown to reduce to
spectral considerations accessible by standard normal-modes analysis
(see Section \ref{spectral}).
However, they are studied with different motivations, and
yield different results.
Historically, it appears that instability motivated the physical
study, which has focused on the detonation case and spectral stability
criteria \cite{BE}.  
By contrast, the mathematical study has focused
on stability and the somewhat simpler shock wave case, with the main
(difficult!) issue being to establish full nonlinear stability 
assuming that spectral stability holds in the form of a suitable
``Lopatinski'' or ``Evans'' condition 
\cite{ M3, MeKochel, Z1, GMWZ2}.

Regarding bifurcation, there is strong evidence that the detonation 
instabilities described above correspond to a 
{\it Poincar\'e--Hopf bifurcation}
(indeed, our description above benefits much by hindsight).
In particular, Erpenbeck \cite{Er2}
and Bourlieux, Majda, and Roytburd \cite{BMR} have carried
out formal asymptotics in support of this viewpoint for the one-dimensional
(galloping, or pulsating) case, in the context of the ZND equations.
More recently, Kasimov and Stewart \cite{KS}
have carried out a definitive study in the multidimensional case by
numerical linearized normal modes analysis, also for the ZND equations, in
which they demonstrate that the onset of such instabilities indeed
corresponds to crossing of the imaginary axis by a conjugate pair of
eigenvalues toward the unstable (positive real part) side, with, moreover,
excellent correspondence between observed nonlinear oscillations 
and the associated normal modes.  
Depending whether the maximally unstable mode is in the
longitudinal or transverse direction, the oscillatory behavior is
seen to be of galloping or spinning (cellular) type.
However, up to now, there has been 
carried out no rigorous analysis of these phenomena.
To fill this gap is the object of the present study.

\subsection{Spectral considerations}\label{spectral}

The main difficulty in the analysis of shock or detonation waves
is the lack of a spectral gap between stationary
or oscillatory modes and essential spectrum 
of the linearized operator about the wave.
In the inviscid, hyperbolic case \eqref{Euler}, \eqref{rEuler},
the linearized problem has essential spectrum filling
the imaginary axis.
In the viscous case \eqref{NS}, \eqref{rNS}, the essential spectrum
is, rather, tangent at the origin to the imaginary axis.

For example, 
without loss of generality taking a standing-wave solution $u=\bar u(x)$, 
$s=0$ (i.e., working in coordinates moving with the wave), 
and linearizing \eqref{NS} about $\bar u$, we obtain,
taking $B\equiv I$ for simplicity,
$$
u_t =Lu:= u_{xx}- (Au)_x,  \qquad
A:= df(\bar u(x)),
$$
where the asymptotically constant-coefficient operator $L$ 
is what we have called the linearized operator about the wave.
A standard result of Henry
(\cite{He}, Theorem A.2, chapter 5) on spectrum of operators
 with asymptotically constant coefficients asserts that
the rightmost (i.e. largest real part) envelope of the essential
spectrum with respect to any $L^p$ is the envelope of the 
union of the rightmost envelopes
of the spectra of the limiting, constant-coefficient operators 
at $\pm \infty$: in this case,
\be\label{Lplusminus}
L_\pm := \partial_x^2 - A_\pm \partial_x,
\qquad
A_\pm:= df(u_\pm).
\ee
The $L^2$ spectra of $L_\pm$ may be computed by Fourier transform to
be the curves traced out by dispersion relations
\be\label{disp}
\lambda(\xi):= ia_j^\pm \xi - \xi^2,
\qquad
\xi\in \RR^1,
\ee
where $a_j^\pm$ are the eigenvalues of $A_\pm$.
These eigenvalues are {\it real} and {\it nonzero} 
under the physical assumptions
that \eqref{Euler} be hyperbolic in a neighborhood of $u_\pm$ 
and the shock be noncharacteristic.
Likewise \cite{Sat}, $L$ has always an $L^2$ 
zero-eigenvalue with eigenfunction
$\bar u'(x)$, associated with translation of the wave.

In the absence of a spectral gap, {\it standard stability and
bifurcation theorems do not apply}, and so we must carry
out a refined analysis.

\subsubsection{Spatio-temporal description}\label{xt}

The spectral configuration just described translates in
the $x$-$t$ plane, in the stable case that there exist no
other eigenvalues $\lambda$ of $L$ in the nonnegative complex
half-plane $\Re \lambda\ge 0$, to the
following description developed in \cite{Liu86, Liu97, ZH, Z4, MaZ3}
of the Green function $G(x,t;y)$ associated with $\partial_t-L$.
A point source, or delta-function initial datum, originating at $y>0$
will propagate initially as an approximate superposition of Gaussians
with constant mass (total integral), 
centered along hyperbolic characteristics $dx/dt=a_j^+$
determined by the asymptotic system at $+\infty$.
Those propagating in the positive direction will continue out to $+\infty$;
those propagating in the negative direction will continue until they
strike the shock layer at approximately $x=0$, whereupon they will be
transmitted and reflected along outgoing characteristic directions
$dx/dt=a_k^-$, $a_k^-<0$ and $dx/dt=a_k^+$, $a_k^+>0$.
In addition, there will be deposited at the shock layer a certain
amount of mass in the stationary eigenmode $\bar u'(x)$
(integral $\int \bar u(x)dx= u_+-u_-\ne 0$), corresponding
to translation of the background wave.

We refer the reader to \cite{ZH, MaZ3} for further discussion and details.
For the present purpose, it suffices to note that the Green function
may be modeled qualitatively as the sum of terms
\be\label{Kterm}
K(x,t;y):= t^{-1/2}e^{-(x-y-at)^2/4t}
\ee
and
\be\label{Eterm}
J(x,t;y):= \bar u'(x) \errfn ((-y-at)/2t^{1/2})
\ee
propagating with noncharacteristic speed $a < 0$,
where
\begin{equation} \label{errfunction} 
\errfn (z) := \frac{1}{2\pi} \int_{-\infty}^z e^{-\xi^2} d\xi,
\end{equation} 
the first modeling moving heat kernels (Gaussian signals), and
the second the excitation of the zero-eigenfunction by incoming
signals from $y\ge 0$.

{}From this description, we see explicitly that decaying modes decay
time-algebraically and not exponentially, in agreement with the
lack of spectral gap.  Moreover, the low-frequency/large-time--space
behavior is approximately hyperbolic, propagating along characteristics
associated with the endstates $u_\pm$.

\subsubsection{Evans/Lopatinski determinants and transition to instability
}\label{evans}

Further information may be obtained 
at the spectral level 
by the comparison of Evans and Lopatinski determinants $D$ and $\Delta$
for the viscous and inviscid problems \eqref{NS} and \eqref{Euler}.
The Evans function $D(\lambda)$, defined as a Wronskian of functions spanning
the decaying manifolds of solutions of the eigenvalue equation
$$
(L-\lambda)u=0
$$
associated with $L$ at $x\to +\infty$ and $x\to -\infty$ is an
analytic function with domain containing $\{\Re \lambda \ge 0\}$,
whose zeroes away from the essential spectrum
correspond in location and multiplicity with eigenvalues
of $L$ with respect to any $L^p$.
Its behavior is also closely linked with that of the resolvent kernel
of $L$, i.e., the Laplace transform with respect to time of the
Green function $G$;
 see \cite{AGJ, GZ, ZH, ZS, Z1, Z2} for history and further details.
The corresponding object for the inviscid linearized problem is
the Lopatinski determinant $\Delta(\lambda)$ defined in
\cite{K, M1, M2, MeKochel}, a homogeneous function of degree one:
in the one-dimensional case, just linear.

An important relation between these two objects, 
established in \cite{GZ, ZS, MeZ2}, is the low-frequency expansion
\be\label{ZS}
D(\lambda)=\gamma \Delta(\lambda) + o(|\lambda|)
\ee
quantifying the above
observation that low-frequency behavior is essentially hyperbolic,
where $\gamma$ is a constant measuring transversality of the
profile $\bar u$ as a connecting orbit of equilibria $u_\pm$
in the associated traveling-wave ODE, nonvanishing for transversal
connections.

A corresponding expansion 
\be\label{detZS}
D_{rNS}(\lambda)=\gamma \Delta_{CJ}(\lambda) + o(|\lambda|)
\ee
holds for the Evans function $D_{rNS}$ associated with detonation
wave solutions of \eqref{rNS}, where $\Delta_{CJ}$ denotes the
Lopatinski condition, not for \eqref{rEuler}, but for the still simpler
inviscid--instantaneous-reaction-rate
{\it Chapman--Jouguet} model in which detonations are modeled
by piecewise constant solutions \eqref{shock} across which combustion
proceeds instantaneously; see \cite{Z1, LyZ1, LyZ2, JLW, CJLW} for further
details.
That is, the low-frequency behavior is not only hyperbolic but 
``instantaneous'', with small-scale details of the profile structure
lost.
Both \eqref{ZS} and \eqref{detZS} extend to the corresponding
multidimensional problem of a planar shock moving in $\RR^d$
(different from the finite-cross-sectional case considered here). 

As calculated in \cite{M1} and \cite{LyZ1}, respectively,
neither $\Delta$ nor $\Delta_{CJ}$ vanishes for $\lambda \ne 0$
for an ideal-gas equation of state.\footnote{Indeed, this holds for 
multidimensions as well; see \cite{M1, JLW}.}
Likewise, $\gamma$ does not vanish for any choice of parameters
for ideal gas dynamics by a well-known result of Gilbarg \cite{Gi}, 
or, by results of \cite{GS} for reactive gas dynamics with 
ideal gas equation of state in the ZND-limit $|B|\to 0$.
{}From these observations, combined with \eqref{ZS}--\eqref{detZS}, 
we may deduce that $D$ and $D_{rNS}$ 
have for all choices of physical parameters precisely
one zero at $\lambda=0$, corresponding to translation-invariance
of the background equations.

In particular, as physical parameters are varied, starting from
a stable viscous traveling-wave, transition to instability,
signalled by passage from the stable complex half-plane $\{\Re \lambda<0\}$ 
to the unstable complex half-plane $\{\Re \lambda>0\}$ of one
or more eigenvalues of the linearized operator about the wave,
{\it cannot occur through passage of a real eigenvalue through
the origin, but rather must occur through the passage of 
one or more nonzero complex conjugate pairs through the imaginary axis:}
that is, a Poincar\'e--Hopf-type configuration.
This observation, made in \cite{LyZ1, LyZ2}, gives rigorous
corroboration at the spectral level of the numerical observations
of Kasimov and Stewart \cite{KS}.
What remains is to convert this spectral information into a
rigorous nonlinear existence result.

\subsubsection{Formulation of the problem}\label{formule}
{}From physical/numerical observations,
one expects for an ideal gas equation of state
that such transition seldom or never occurs for shocks
but frequently occurs for detonations.
However, from a mathematical point of view, the situation
for shocks and detonations appears to be entirely parallel.
We may thus phrase a common mathematical problem:
\medbreak

{\bf (P)} {\it Let $\bar u^\eps$ denote a family of traveling-wave (either
shock or detonation) solutions indexed by bifurcation parameter $\eps\in \RR^1$,
with associated linearized operators $L(\e)$ transitioning at $\eps=0$ 
from stability to instability
as follows. 

Assuming the ``generic'' spectral situation that 
\begin{itemize}\item[{\rm (i)}] the $L^2$ 
essential spectrum of each $L(\e)$ is contained in 
$\{\Re \lambda <0\} \cup \{0\}$,
\item[{\rm (ii)}] the translational zero-eigenvalue of each $L(\e)$ is simple in the
sense that the associated Evans function $D_\eps$
vanishes at $\lambda=0$ with multiplicity one, and 
\item[{\rm (iii)}] a single complex conjugate pair of $L^2$ eigenvalues 
of $L(\e)$ (zeroes of $D_\eps$)
$$
\lambda_\pm(\eps)= \gamma(\eps)\pm i\tau(\eps),
\qquad \gamma(0)=0, \quad \tau(0)\ne 0, \quad (d\gamma/d\eps)(0) > 0
$$
crosses the imaginary axis with positive
speed, and the rest of the 
$L^2$ point spectrum of $L(\e)$ is located in $\{\Re \lambda <0\},$
\end{itemize} show that there occurs a
Poincar\'e--Hopf bifurcation from the family of traveling wave solutions $\bar u^\e$ to nearby time-oscillatory solutions.}
\medbreak

 We examine this problem in a series of different contexts.

\subsection{Model analysis I: the scalar case}\label{model1}

Let us first recall the situation considered in \cite{TZ}, of a 
one-parameter family of standing-wave solutions $\bar u^\eps(x)$
of a smoothly-varying family of equations
\begin{equation}
\label{epseqn}
u_t =\CalF(\e, u):= u_{xx}-F(\e, u, u_x)
\end{equation}
(possibly shifts 
$F(\eps, u, u_x):= f(u,u_x)-s(\eps) u_x$
of a single equation written in coordinates
$x\to x-s(\e)t$ moving with traveling-wave solutions
of varying speeds $s(\e)$),
with linearized operators
$L(\e) :=\partial \CalF/\partial u|_{u=\bar u^\e}$
for which a spectral gap may be recovered in an appropriate
exponentially-weighted norm 
$$
\|f\|_{H^{2}_\eta}^2:=
\sum_{j=0}^2
\|(d/dx)^j f(x)\|_{L^{2}_\eta}^2,  \quad \|f\|_{L^{2}_\eta}:=
\|e^{\eta (1+|x|^2)^{1/2}}f(x)\|_{L^{2}},
\quad \eta>0.
$$
We call this the {\it weighted norm condition}.

This approach, introduced by Sattinger \cite{Sat} applies
to the case (see Section \ref{xt}) that signals in the far 
field are convected under the linearized evolution equation
$u_t=L(\e) u$ {\it inward} toward the background profile,
hence time-exponentially decaying in the weighted norm
$\|\cdot\|_{H^2_\eta}$ penalizing distance from the origin.
The method encompasses both \eqref{NS} and \eqref{rNS}
in the scalar case $u$, $z\in \RR^1$, with artificial 
viscosity $B=D=I$. It has also interesting applications
to certain reaction-diffusion systems and the related
Poincar\'e--Hopf phenomenon of ``breathers'' \cite{NM, IN, IIM}.
However, it does not apply to either shock or detonation
waves in the system case $u\in \RR^n$, $n\ge 2$.

In this case, we have the following rather complete result,
including stability along with bifurcation description.

\begin{prop} [\cite{TZ}] \label{oldPH}
Let $\bar u^\eps$, \eqref{epseqn}  
be a family of traveling-waves and systems 
satisfying the weighted norm condition 
and assumptions {\rm (P),} with $F\in C^4$.
Then, for $a\ge 0$ sufficiently small and $C>0$ sufficiently large,
there are $C^1$ functions 
$\e(a)$, $\e(0)=0$, and $T^*(a)$, $T^*(0)=2\pi/\tau(0)$,
and a $C^1$ 
family of solutions  
\be\label{persoln}
u^{a}(x,t)= {\bf u}^a(x-\sigma^a t,t)
\ee
of \eqref{epseqn} with $\e=\e(a)$,
where ${\bf u}^a(\cdot, t)$ is time-periodic with period $T^*(a)$ 
and $\sigma^a$ is a constant drift, such that
\be \label{desc}
C^{-1}a\le \|{\bf u}^a(\cdot, t)-\bar u^{\eps(a)}\|_{H^2_\eta}\le Ca
\ee
for all $t\ge 0$.
Up to fixed translations in $x$, $t$, for $\eps$ sufficiently small,
these are the only
nearby solutions of this form, as measured in $H^2_\eta$.
Moreover, solutions ${\bf u}^a$ are time-exponentially phase-asymptotically
orbitally stable with respect to $H^2_\eta$ if $d\e /da>0$,
in the sense that perturbed
solutions converge time-exponentially to a specific shift in $x$ and $t$
of the original solution, and unstable
if $d\e/da<0$.\footnote{In \eqref{desc}, 
we have repaired an obvious error in \cite{TZ}, where $u^0$ appears in
place of $u^\eps$, and in \eqref{persoln} eliminated a redundant
time-periodic translation $\theta^a(t)$.}
\end{prop}

Proposition \ref{oldPH} was established in \cite{TZ} by 
center-manifold reduction, with the main issue being to accomodate
the underlying group invariance of translation.
The basic idea is to coordinatize $u$ as $(v,\a),$ where $\a$ parameterizes the group invariance, then work on the quotient space $v,$ reducing the problem
from relative to standard Poincar\'e--Hopf bifurcation, afterward
recovering the location $\alpha$ (a function of $a$ and $t$) driven by the solution
on the quotient space, by quadrature.
The integral of the periodic driving term yields a periodic part that may be subsumed in the profile and a drift $\sigma^a t$.
See \cite{TZ} for further discussion and details.

\br\label{decay}
\textup{
A consequence of \eqref{desc} (by Sobolev embedding) is 
\be\label{dec}
|{\bf u}^a(x, t)-\bar u^\eps| \le Ce^{-\eta|x|}.
\ee
That is, the existence result in weighted norm space includes
quite strong information on the structure of the wave.
On the other hand, the stability result is somewhat weakened by the
appearance of a spatial weight, being restricted to exponentially
decaying perturbations.
}
\er

\subsection{Model analysis II: systems of conservation laws}\label{model2}

The purpose of the present paper is to extend the analysis initiated
in \cite{TZ} for scalar models to the more physically realistic 
system case.
For simplicity, we restrict to the somewhat simpler
case of viscous shock solutions of systems
with artificial viscosity $B\equiv I$; however,
the method of analysis in principle applies also in the general case;
see discussion, Section \ref{discussion}.

Specifically, consider a one-parameter 
family of standing viscous shock solutions 
\be\label{profintro}
u(x,t)=\bar u^\eps(x),
\qquad \lim_{z\to \pm \infty} \bar u^\eps(z)=u_\pm^\eps
\quad \hbox{\rm (constant for fixed $\eps$)},
\ee 
of a smoothly-varying family of conservation laws 
\begin{equation}
\label{sysepseqn}
u_t =\CalF(\e, u):= u_{xx}- F(\e, u)_x,
\qquad u\in \RR^n,
\end{equation}
with associated linearized operators
\be\label{Ldefintro}
L(\e) :=\frac{\partial \CalF}{\partial u}|_{u=\bar u^\e}
= -\partial_x A^\eps(x) + \partial_x^2,
\ee
denoting
\begin{equation} \label{A} A^\eps(x):= F_u(\bar u^\eps(x), \eps), \qquad A^\eps_\pm:=\lim_{z\to \pm \infty} A^\eps(z)=F_u(u^\eps_\pm, \eps).
\end{equation}

We take
$\bar u^\eps$ to be of standard {\it Lax type}, 
meaning that the hyperbolic
convection matrices $A^\eps_+$ and $A^\eps_-$ 
at plus and minus spatial infinity have, respectively, $p-1$ negative and
$n-p$ positive real eigenvalues for $1\le p \le n$, where $p$ is the 
characteristic family associated with the shock: in other words, there
are precisely $n-1$ outgoing hyperbolic characteristics in the far field.

This is the only type occurring for gas dynamics with standard (e.g., ideal
gas) equation of state; for reacting gas dynamics, the corresponding
object is a {\it strong detonation}, which is the only (nondegenerate)
type occurring in the ZND limit $B, D\to 0$ \cite{GS}.
The special features of Lax-type shocks (resp. strong detonations), 
as compared to more general undercompressive shocks 
(resp. weak detonations) that can occur in other settings,
turn out to be important for the analysis 
(see Remark \ref{linuc}),
in sharp contrast with the generality of \cite{TZ}.

Note, for the system case $n\ge 2$, 
that there is at least one outgoing characteristic,
so that the weighted norm methods of the previous section do not apply.
In particular, we see no way to construct a center manifold for
this problem, and suspect that one may not exist; the slow
(time-algebraic) decay rate of outgoing modes evident in our description
of the Green function in Section \ref{xt} suggests an essential obstacle
to such construction.

On the other hand, the conservative form of system \eqref{sysepseqn}
implies that relative mass $\int u(x) dx$ is conserved
for all time for perturbations  $u=  \tilde u - \bar u^\eps$, 
$\tilde u$ satisfying \eqref{sysepseqn}.
Thus, at least at a formal level,
there is a convenient invariant subspace of \eqref{sysepseqn}
consisting of perturbations with zero excess mass, on which 
which the zero eigenvalue associated with translational invariance
is removed.
For, recall that $\int \bar u'(x) dx= u_+-u_- \ne 0$, so that the
associated zero-eigenfunction is not in the subspace.
On the other hand, nonzero eigenfunctions always have zero mass \cite{ZH},
since $\lambda \phi= L(\e) \phi$ implies $\lambda \int \phi(x)dx=0$
by divergence form of $L(\e)$, hence the crossing nonzero imaginary
eigenvalues $\lambda_\pm (\eps)$ in (P) persist.
Thus, the spectral scenario (P) translates in the zero-mass subspace to
a standard Poincar\'e--Hopf scenario, with no additional zero-eigenvalue,
for which the translational group-invariance need not be taken into
account.  Recall, that this was the main issue in the analysis of
the scalar case.

In short, the two problems (scalar vs. system) {\it have essentially
complementary mathematical difficulties}, hence little technical contact.
Accordingly, the analysis has a quite different flavor in the system
case, depending on both the full, pointwise Green function bounds
of \cite{MaZ3} and the special, conservative structure of the equations.

Our result in this case, and the main result of the paper, is as follows.

\begin{theo} \label{newPH}
Let $\bar u^\eps$, \eqref{sysepseqn}
be a family of traveling-waves and systems 
satisfying assumptions {\rm (P),} with $F\in C^2$.
Assume further that $\bar u^\e$ is of Lax type,
with $\sigma(A_\pm^\e)$ 
real, nonzero, and simple, $A_\pm^\eps$ as
defined in {\rm \eqref{A}}.
Then, for $a\ge 0$ sufficiently small and $C > 0$ sufficiently
large, there are $C^1$ functions
$\eps(a)$, $\e(0)=0$, and $T^*(a)$, $T^*(0)=2\pi/\tau(0)$, and a $C^1$ family of 
solutions $u^{a}(x,t)$ of \eqref{sysepseqn} 
with $\e=\e(a)$, time-periodic of period $T^*(a)$, such that
   \begin{equation} \label{xbd}
   C^{-1} a \, \leq \,
   \sup_{x \in \R} \, (1 + |x|) \big| u^{a}(x,t) - \bar u^{\e(a)}(x) \big| \, \leq Ca, \quad
 \mbox{for all $t \geq 0.$}
   \end{equation}
Up to fixed translations in $x$, $t$, for $\eps$ sufficiently
small, these are the only nearby solutions as measured in 
norm $\|f\|_{X_1}:=\|(1+|x|)f(x)\|_{L^\infty(x)}$ that are
time-periodic with period
$T\in [T_0, T_1]$, for any fixed $0<T_0<T_1<+\infty$.
If $u_+^\eps\ne u_-^\eps$, they are the only 
nearby solutions of the more general form \eqref{persoln}.
\end{theo}

Note that the statement of Theorem \ref{newPH}
is considerably weaker than that of Proposition \ref{oldPH},
asserting no stability information, and only the relatively weak
structural information of algebraic decay \eqref{xbd} at $\pm \infty$,
as compared to the exponential bound \eqref{dec}.
As suggested by our formal ``zero-mass'' discussion, the periodic
solutions constructed have zero mean drift $\sigma$, in contrast 
to \eqref{persoln}.
However, somewhat surprisingly, they do {not} appear to
have zero excess mass. 
See Remarks 
\ref{zeromass} and \ref{restricted} for further discussion and explanation.

\subsection{Model analysis III: flow in an infinite cylinder}\label{model3}
One-dimensional traveling-wave solutions \eqref{profile} 
of \eqref{NS} with $B\equiv I$
may alternatively be viewed as the restriction to one dimension of
a planar viscous shock solution
\be\label{planar}
u(x,t)=\bar u(x_1-st)
\ee
of a multidimensional system of viscous conservation laws
\be\label{multiNS}
u_t + \sum f^j(u)_{x_j}= \Delta_x u, 
\qquad u\in \RR^n, \, x\in \RR^d, \, t\in \RR^+
\ee
on the whole space. 
Likewise, traveling-wave solutions \eqref{planar} may be viewed
as planar traveling-wave solutions of \eqref{multiNS} on an infinite cylinder
$$
\CalC:= \{x:\, (x_1, \tilde x)\in \RR^1\times \Omega\},
\qquad
\tilde x=(x_2, \dots, x_d)
$$
with bounded cross-section $\Omega\in \RR^{d-1}$,
under artificial {\it Neumann boundary conditions}
$$
\partial u/\partial_{\tilde x} \cdot \nu_\Omega=0
\quad \hbox{\rm for} \quad \tilde x\in \partial \Omega,
$$
or, in the case that $\Omega$ is rectangular, {\it periodic boundary
conditions}, $\Omega=T^{d-1}$, $T^{d-1}$ the rectangular torus.

We take this as a simplified mathematical model for flow in
a duct, in which we have neglected boundary-layer phenomena along the
wall $\partial \Omega$ in order to isolate the oscillatory phenomena
of our main interest.
Consider a one-parameter family of standing planar viscous shock 
solutions $\bar u^\eps(x_1)$ of a smoothly-varying family of conservation laws 
\begin{equation}
\label{multisysepseqn}
u_t =\CalF(\e, u):= \Delta_x u- \sum_{j=1}^d F^j(\e, u)_{x_j},
\qquad u\in \RR^n
\end{equation}
in a fixed cylinder $\CalC$, with Neumann (resp. periodic) 
boundary conditions--
typically, shifts 
$\sum F^j(\eps, u)_{x_j}:= \sum f^j(u)_{x_j} -s(\eps) u_{x_1}$
of a single equation \eqref{multiNS}
written in coordinates
$x_1\to x_1-s(\e)t$ moving with traveling-wave solutions
of varying speeds $s(\e)$--
with linearized operators
$L(\e) :=\partial \CalF/\partial u|_{u=\bar u^\e}$.
As in the previous subsection, we take
$\bar u^\eps$ to be of standard Lax type, considered as
a shock wave in one dimension.
For simplicity, we take $\Omega=\TT^{d-1}$.

Then, we have the following result generalizing Theorem \ref{newPH}.

\begin{theo}\label{multinewPH}
Let $\bar u^\eps$, \eqref{multisysepseqn}
be a family of traveling-waves and systems 
satisfying assumptions {\rm (P),} with $F\in C^2$, $\Omega=\TT^{d-1}$.
Assume further that $\bar u^\e$ is of Lax type,
with $\sigma(A_\pm^\e)$ 
real, nonzero, and simple, $A_\pm^\eps$ as
defined in {\rm \eqref{A}}.
Then, for $a\ge 0$ sufficiently small and $C > 0$ sufficiently
large, there are $C^1$ functions
$\eps(a)$, $\e(0)=0$, and $T^*(a)$, $T^*(0)=2\pi/\tau(0)$, and a $C^1$ 
family of 
solutions $u^{a}(x_1,t)$ of \eqref{sysepseqn} 
with $\e=\e(a)$, time-periodic of period $T^*(a)$, such that
 $$
   C^{-1} a \leq  \,
   \sup_{x_1 \in \R} \, (1 + |x_1|) \big| u^{a}(x_1,t) - \bar u^{\e(a)}(x_1) \big| \, \leq Ca, \quad \mbox{for all $t \geq 0.$}
 $$
Up to fixed translations in $x$, $t$, for $\eps$ sufficiently
small, these are the only nearby solutions as measured in 
norm $\|f\|_{X_1}:=\|(1+|x_1|)f(x)\|_{L^\infty(x)}$ that are
time-periodic with period
$T\in [T_0, T_1]$, for any fixed $0<T_0<T_1<+\infty$.
If $u_+^\eps\ne u_-^\eps$, they are the only 
nearby solutions of the more general form \eqref{persoln}.
\end{theo}

\br\label{multireg}
\textup{
Extension to a circular cross-section
$\Omega:=\{|\tilde x|\le R\}$ follows by 
Bessel function expansion as in \cite{KS}.
The extension to general cross-sectional geometries is discussed in
Section \ref{gencross}.
}
\er

\subsection{Idea of the proof}\label{idea}

We now briefly discuss the ideas behind the proof of 
Theorems \ref{newPH} and \ref{multinewPH}.
Under the spectral assumptions (P) (stated in Section \ref{formule}), there exist smooth
(in $\eps$) $L(\eps)$-invariant projections onto the two-dimensional eigenspace $\Sigma^\e$ 
of $L(\eps)$ associated with
the pair of crossing eigenvalues 
$\lambda_\pm(\eps)= \gamma(\eps) \pm i\tau(\eps)$
and its complement $\tilde \Sigma^\e$.
Projecting onto these subspaces, and rewriting in polar
coordinates the flow on $\Sigma^\e$,
we may thus express \eqref{sysepseqn} in standard fashion as
\ba\label{PHeq1}
\dot r&= \gamma(\eps)r + N_r(\e, r, \theta, v),\\
\dot \theta&= \tau(\eps) + N_\theta(\e, r, \theta, v),\\
\dot v&= \tilde L(\eps)v + N_v(\e, r, \theta, v),\\
\ea
where the ``transverse linearized operator'' $\tilde L(\e)$ is 
the restriction of $L(\eps)$ to $\tilde \Sigma^\e$,
and $N_j$ are higher-order terms coming from the nonlinear part
of \eqref{sysepseqn}: $N_r$ and $N_v$ quadratic order
in $r,v$ and their derivatives, and $N_\theta$ linear order.

We are precisely interested in the case that $\tilde L$ has no spectral
gap, i.e., $\Re \sigma(L)\not \le -\theta<0$
for any $\theta>0$.
One may think for example to the finite-dimensional case that $\tilde L(0)$
has additional pure-imaginary spectra besides 
$\lambda_\pm(0)=\pm i\tau(0)$.
Thus, the center manifold may be of higher dimension, and
one cannot follow the usual course of reducing to a center manifold
involving only $r$, $\theta$, $\eps$ and applying the standard, 
two-dimensional Poincar\'e--Hopf Theorem.
Instead, we proceed by a direct analysis, combining
the two-dimensional Poincar\'e return map construction
with Lyapunov--Schmidt reduction.

\subsubsection{Return map construction}\label{return}
Specifically, truncating $|v|\le Cr$, $C \gg 1$, in the arguments of 
$N_r$, $N_\theta$, we obtain 
$$
|N_\theta|\le C_2 r,
$$
and therefore $ \dot \theta\ge \tau(0)/2>0 $ 
for $\eps$, $r$ sufficiently small.
Thus, in seeking periodic solutions 
\be\label{per}
(r,\theta, v)(T)= (r,\theta, v)(0),
\ee
we may eliminate $\theta$, solving
for $T(\e,a,b)$ as a function of initial data 
$$
(a,b):=(r,v)(0)
$$
and the bifurcation parameter $\eps$, with 
$ T(0,0,0)=2\pi/\tau(0)$, and seek solutions \eqref{per}
as fixed points
\be\label{it}
(a,b)=(r, v)(T(\e,a,b))
\ee
of the {Poincar\'e return map} $(r,v)(T(\e,\cdot, \cdot))$.

Using Duhamel's formula/variation of constants, we may express
\eqref{it} in a standard way (see, e.g., \cite{HK}) as
\ba\label{Pieqns}
0&=f(\e,a,b)=
\big(e^{\gamma(\eps)T(\e,a,b)} -1\big)a + N_1,\\
0&=g(\e,a,b)=
\big(e^{\tilde L (\eps)T(\e,a,b)} -\Id\big)b + N_2,\\
\ea
where
\ba\label{Duhamel}
N_1&:= \int_0^{T(\e,a,b)} e^{\gamma(\eps)(T(\e,a,b)-s)}
N_r(\e,r, \theta, v)(s)ds,\\
N_2&:= \int_0^{T(\e,a,b)} e^{\tilde L(\eps)(T(\e,a,b)-s)}
N_v(\e,r, \theta, v)(s)ds\\
\ea
are, formally,  quadratic order terms in $a$, $b$. 
Note that $(r,\theta, v)$ are functions of $(\e,a,b)$
as well as $s$, through the flow of \eqref{PHeq1}, with
$$
\|(r, v)\|\le C\|(a,b)\|
$$
for any ``reasonable'' norm in the sense that
\eqref{PHeq1} are locally well-posed
(in practice, no restriction).
By quadratic dependence of $N_j$, we have, evidently,
that $(\e,a,b)=(\e,0,0)$ is a solution of \eqref{Pieqns} for all $\eps$.

\subsubsection{Standard reduction}\label{standard}
Continuing in this standard fashion, we should next like
to perform a Lyapunov-Schmidt reduction, using the Implicit
Function Theorem to solve the $g$ equation for
$b$ in terms of $(\e,a)$, i.e., to find a function
$$
b=B(\e,a), 
\qquad
B(0,0)=0,
$$
satisfying
$g(\e,a,B(\e,a))\equiv 0$: equivalently, to
reduce to the nullcline of $g.$

If this were possible, from $|N_2|\le C|(a,b)|^2$ in \eqref{Pieqns}(ii),
we would find, further, that 
$ |B(\e,a)|\le C_3 |a|^2$,
from which straightforward Duhamel/Gronwall estimates on
\eqref{PHeq1}(i),(ii) would yield
$$
|v|\le C_2r^2
$$
for $ s\in [0,T]$, justifying a posteriori the truncation 
$|v|\le Cr$, of $N_\theta$ performed
in the first step, provided that $(\e,a,b)$ are
taken sufficiently small: 
specifically, small enough that $r$ remains less than or
equal to  $C/C_2$ for $s\in [0,T]$.

Substituting into the $f$ equation, we would then obtain
a reduced, scalar bifurcation problem
$$
0=f^*(\e,a):=f(\e,a,B(\e,a)),
\qquad
f^*(\e,0)\equiv 0
$$
that would be solvable in the usual way (i.e., by dividing out $a$ and applying
the Implicit Function Theorem a second time, using
$$
\partial_\eps \big(a^{-1}\pi\big)(0,0)=\partial_\eps (\gamma/\tau)(0)
= \partial_\eps \gamma(0)/\tau(0)\ne 0
$$
to solve for $\eps$ as a function of $a$; 
see, e.g., \cite{HK}, or Section \ref{bifapp}). The Poincar\'e return map system \eqref{Pieqns} would then be solved, and  a solution of \eqref{Pieqns} would generate a periodic solution of \eqref{PHeq1}.  

\subsubsection{Regime of validity}\label{validity}
Let us ask ourselves now under what circumstances this
basic
reduction procedure may actually be carried out.
Formal differentiation of \eqref{Pieqns} yields
$$
\partial_b g(0,0,0) = 
e^{2\pi(\tilde L/\tau)(0) } -\Id.
$$
In the finite-dimensional ODE case, the usual condition of application of the (standard)
Implicit Function Theorem is thus that 
$(e^{2\pi(\tilde L/\tau)(0)} -\Id)$ be invertible, i.e.,
that $\tilde L(0)$ have no ``resonant'' oscillatory modes,
i.e., pure imaginary
eigenvalues that are integer multiples of the crossing
modes $\lambda_\pm(0)=\pm i\tau(0)$.
Note that this includes interesting cases for which $\tilde L$
has no spectral gap, and thus that do not yield to
the standard center-manifold reduction.
%
The analogous criterion in the infinite-dimensional (e.g., PDE) setting
is that all possible resonant modes $ni\tau(0)$, $n\in \ZZ$, 
{\it lie in the resolvent set of} $\tilde L$.

This is encouraging, and shows that the requirements of Lyapunov-Schmidt
reduction are much less than those for center-manifold reduction:
in particular, only spectral separation (from both $\pm i\tau(0)$ and
their aliases $ni\tau(0)$) rather than spectral gap
is needed.
Unfortunately, in the setting \eqref{sysepseqn}, \eqref{multisysepseqn}
of our interest
we have not even spectral separation, since essential spectra of 
$\tilde L(\eps)$ accumulates at $\lambda=0$ for every $\eps$.
In this case, or others for which separation fails,
we must follow a different approach.

\subsubsection{Refined analysis}\label{refined}

Alternatively, we may rewrite \eqref{Pieqns}(ii) formally as
a fixed-point equation
\be\label{ref}
b= \big(\Id - e^{\tilde L(\eps) T(\e,a,b)}\big)^{-1} N_2(\e,a,b),
\ee
then apply the Contraction-mapping Principle to
carry out an Implicit Function construction ``by hand''
(reminiscent of, but not exactly the standard proof of
the Implicit Function Theorem).
In the standard case that $|e^{\tilde Lt}|$ is exponentially
decaying ($\sim$ spectral gap), 
we may rewrite \eqref{ref} using Neumann expansion as
\be\label{neuversion}
b= \sum_{j=0}^{\infty} e^{j\tilde L(\eps) T(\e,a,b)} N_2(\e,a,b).
\ee
The basis for our analysis is the
simple observation that, more generally, provided the series on the
righthand side converges (conditionally) for $\|(a,b)\|$ bounded,
then (i) (by scaling argument) the righthand side is contractive in $b$ for 
$\|(\e,a,b)\|$ sufficiently small, and (ii) (by standard telescoping sum
argument, applying 
$\big(\Id - e^{\tilde L(\eps) T(\e,a,b)}\big)$ to \eqref{ref})
the resulting solution 
$$ 
b=B(\e,a)
$$
guaranteed by the Contraction-mapping Principle
is in fact a solution of the original equation \eqref{Pieqns}(ii).
See Section \ref{framework} for further details.

In the context of \eqref{sysepseqn}, by divergence form of the
equation, $N_2=(n_2)_x$, with $n_2= O(|(r,v)|^2)$.  For simplicity
of discussion, model $n_2$ as a bilinear form in $(r,v)$, 
so that (by quadratic scaling)
$n_2\in L^1$ for $v\in L^2$, from which we obtain 
$N_2= n_x$,
$n\in L^1$ for $b\in L^2$.
Thus, replacing the righthand side of \eqref{neuversion} by its continuous
approximant 
$$
N_2 \,  +  \, T^{-1}\int_{T}^\infty e^{\tilde L(\eps) t} N_2 \, dt \quad = \quad 
N_2 \, + \, T^{-1}\int_{T}^\infty \big( e^{\tilde L(\eps) t} 
\partial_x\big) n \, dt,
\qquad
n\in L^1,
$$
we see that convergence (as a map from $L^2$ to $L^2$) reduces, roughly, to
a space-time stability estimate
\be\label{toyest}
\Big\|\int_{-\infty}^{+\infty} \Big( \int_T^\infty
\tilde G_y(x,t;y) \, dt \Big) \, n(y) \, dy \Big\|_{L^2(x)}\le C\|n\|_{L^1}
\ee
where $\tilde G$ is the Green kernel associated with transverse solution
operator $e^{\tilde Lt}.$ Estimate \eqref{toyest} is quite similar to those arising for the full solution operator
(again, through Duhamel's principle;
see, for example, \cite{Z4, MaZ4})
in the study of nonlinear stability of spectrally stable waves.
That \eqref{toyest} holds for all $n \in L^1,$ is equivalent to 
\be\label{Gbd1}
\sup_y \Big\| \int_T^\infty   
 \tilde G_y(\cdot,t;y)dt \Big\|_{L^2} \le C.
\ee
(One implication is a consequence of the triangle inequality; for the other take a sequence converging to a Dirac function in $L^1$.)
This would hold, for example, if
\be\label{Gbd2}
\sup_y \int_T^\infty   
\| \tilde G_y(\cdot,t;y)dt \|_{L^2} \le C.
\ee

Note that the neglected first term $N_2,$ corresponding to $j = 0$ in \eqref{neuversion}, is $C^2$ from $L^2$ to $L^2$
by the smoothing action of \eqref{Duhamel}(ii); see Section
\ref{returnsection}.

In the arguments above, we have used strongly the specific structure of
the nonlinearity, both quadratic dependence and divergence form, since the $L^2$-operator norm
$\big|e^{\tilde Lt}\big|_{L^2\to L^2}$ does not decay in $t$.
However, this is still not enough.
For, recall the one-dimensional discussion of Section \ref{xt}, 
modeling $\tilde G$ qualitatively as the sum of a convected heat kernel $K$ as in \eqref{Kterm} 
and an error-function term $J$ as in \eqref{Eterm}.
{}From the standard bounds $\|K_y(\cdot, t;y)\|_{L^2}= C t^{-3/4}$, $ \| J_y(\cdot,t;y) \|_{L^2} \leq C t^{-1/2},$
 hence 
(the support of $K$ and $J$ being essentially disjoint)
$\int_T^\infty \tilde G_y(x,t;y)dt $ is not absolutely convergent
in $L^2(x),$
and the bound \eqref{Gbd2} does not hold.  

Instead, we must show {\it conditional convergence}, using
detailed knowledge of the propagator $\tilde G$ to
identify cancellation in 
$\int_T^\infty \tilde G_y(x,t;y)dt  $. 
For example, using
$ K_y= a^{-1}(K_t-K_{yy})$, we find that
\ba\label{cancel}
\int_T^t K_y(x,s;y) \, ds &=
a^{-1}\big( \int_T^t K_t(x,s;y)\, ds -
 \int_T^t K_{yy}(x,s;y)\, ds \big)\\
&= a^{-1}\big(
K(x,s;y)|_T^t - \int_T^t K_{yy}(x,t;y)\, ds\big)\\
\ea
is the sum of $-a^{-1}K(x,T;y) \in L^2(x)$ with terms
$a^{-1}K(x,t;y)\sim t^{-1/4}$ and 
$a^{-1}\int_T^t K_{yy}(x,t;y)\, ds\sim \int_T^t s^{-5/4}ds$
that are respectively decaying and absolutely convergent in $L^2(x)$,
hence converges uniformly with respect to $y$ in $L^2(x)$.

A similar cancellation argument yields convergence of the error-function 
term, along with uniform boundedness with respect to $y$.
Convergence of the error-function term is not uniform with respect to $y$, a
detail which necessitates further technicalities:  in particular,
the weighted estimate 
$$
|v(x)|\le C(1+|x|)^{-1}
$$
leading eventually to \eqref{xbd}.
Likewise, there are further issues associated with $\eps$-regularity
of the solution, needed for bifurcation analysis of the reduced equation.
However, the main idea is contained in calculation \eqref{cancel}.

The multi-dimensional case goes similarly, with the 
computation of the critical neutral, zero transverse wave number
reducing to the one-dimensional case, and all others to the case
of a spectral gap.

\br\label{cancelrmk}
\textup{
It is readily checked that the above arguments go through
in the general case that $n_2$ is quadratic order for $|v|_{L^\infty}\le C$ 
and not bilinear in $(v,r)$, 
substituting $v\in L^2\cap L^\infty$ 
($\Rightarrow n_2\in L^1$)
for $v\in L^2$, working in the $L^2\cap L^\infty$ norm in place of $L^2$,
and substituting for \eqref{Gbd1} the estimate
$$
\sup_y \Big\|\int_T^\infty 
\tilde G_y(x,t;y)dt \Big\|_{L^2\cap L^\infty(x)}\le C.
$$
Similar cancellation estimates 
are important in the study of asymptotic behavior of stable viscous
shock waves \cite{Liu91, SX, Liu97, Raoofi, HRZ}.
}
\er

\br\label{tangent}
\textup{
As in the result of Proposition \ref{oldPH} 
obtained by center manifold reduction, 
the family of periodic solutions obtained, up to translation, 
lies tangent to the 
$(r,\theta)$-plane corresponding to linearized oscillatory behavior.
}
\er

\br\label{zeromass}
\textup{
As noted earlier, operator $S=e^{\tilde L T}$
by divergence form of the equations
preserves the property of zero mass, $\int u dx=0$,
hence each finite approximant $\sum_{j=0}^N S^j$
to $(I-S)^{-1}$ also preserves zero mass.
However, $\int u\, dx=0$ is not closed with 
respect to $\|\cdot\|_{L^2}$, as may be seen
by the example
$$
K(x,s;y)|_T^t  \to -K(x,T;y)
$$
in $L^2$ as $t\to \infty$; that is, mass may ``escape at infinity''
in the limiting process.
Comparing with \eqref{cancel}, we see that the principal term
in $\lim_{N\to \infty}\sum_{j=0}^N S^j \partial_x$ 
is exactly of this form, hence the $L^2$-limit $(I-S)^{-1}$ 
does {not} preserve zero mass, even though each finite approximant does.
Accordingly, there is no reason that the solution 
$b=b(a,\eps)=(I-S)^{-1}{\CalN_2}(b,a,\eps)$ obtained by
reduction should have zero mass:
in the model case $n=u^2$, $G=K$, it can be seen that it does not.
Thus, the idea of carrying out a simplified bifurcation analysis in the 
invariant subspace of zero-mass perturbations, though appealing at a 
formal level, 
in general cannot succeed.
See also Remark \ref{restricted}.
}
\er

\subsection{Discussion and open problems}\label{discussion}

Theorems \ref{newPH} and \ref{multinewPH} together
with the spectral observations of \cite{LyZ1, LyZ2}
give rigorous validation in a simplified context
of the formal and numerical observations of \cite{BMR, KS}.
An interesting problem for future work is to extend these
results to the originally-motivating case of detonation 
waves of the full, reacting compressible Navier--Stokes equations.

We expect that our one-dimensional analysis will extend in 
straightforward fashion, combining tools developed in 
\cite{MaZ3, MaZ4, Z1} and \cite{LRTZ} to treat, respectively,
nonreacting gas dynamics with physical, partial viscosity 
and reacting gas dynamics with artificial viscosity.
Likewise, we expect that we can readily treat flow in a cylinder 
for the physical equations with artificial Neumann or
periodic boundary conditions.
However, the treatment of 
physical, no-slip boundary conditions
(presumably associated with characteristic viscous boundary layers) 
involves technical and philosophical difficulties beyond the
scope of the present analysis, as yet unresolved even
for nonreacting, incompressible flows: for example, even
the construction of a background traveling profile becomes problematic
in this setting.
We point out that viscous boundary effects (since also viscosity) 
are neglected also in the ZND setting of \cite{BMR, KS}.

A second natural direction for future investigation
is the question of stability
of the periodic waves whose existence we have established
here.  
In the absence of a spectral gap,
our method of analysis does not directly yield stability
as in the case of center manifold reduction, but at
best partial information on the location of point spectrum
associated with oscillatory modes, with stability presumably
corresponding to the standard condition $d\eps/da>0$.
The hope is that we could combine such information
with an analysis like that carried out for stationary
waves in \cite{ZH, MaZ3},
adapted from the autonomous to the time-periodic setting: 
that is, a generalized Floquet analysis in the PDE setting
and in the absence of a spectral gap.
We consider this a quite exciting direction for further development
of the theory.

Bifurcation in the absence of a spectral gap has been considered
by a number of authors in different settings; in particular,
it has been studied systematically by
Ioss et al \cite{I1, I2, IA, IK, IM, SS1, SS2} in various contexts
using an alternative ``spatial dynamics'' approach.
It would be very interesting to investigate what results could
be obtained by this technique in the context of viscous shock waves:
more generally, to relate it at a technical level to the one used here.
Considered in this larger context,
the interest of the present approach is that it gives a simple and explicit
connection between stability and bifurcation,
in the spirit of center-manifold reduction and formal asymptotics,
but adapted to the boundary case of zero spectral gap and time-algebraic decay.
On the other hand, the price of this ``direct'' approach is
that the needed estimates may in practice (as here) 
be rather delicate to obtain.

We mention also a recent work of Kunze and Schneider \cite{KuS}
in which they analyze pitchfork bifurcation in the absence of a
spectral gap using Sattinger's weighted-norm method,
as described in Section \ref{model1}, but in a situation where
convection is outward, away from the profile layer.
This entails the use of ``wrong-way'' exponentially decaying
weights at the linearized level, introducing a spatially-exponentially
growing multiplier in quadratic-order source terms, 
in combination with separate,
compensating estimates at the nonlinear level.
For similar arguments in the context of nonlinear stability,
see, e.g., \cite{PW, Do}.
This interesting approach has been used successfully in the shock-wave
context to treat the scalar undercompressive case \cite{Do}; however,
it does not appear to generalize to the system case.

\begin{rems}\label{final}
\textup{
1.  Birtea et al in \cite{BPRT} 
successfully carry out a rather complicated bifurcation analysis
in the ODE setting without explicit knowledge of the background
solution.  Similar techniques might perhaps be useful
in treating flow in a duct with physcial, nonslip
boundary conditions, for which description of
the background flow is itself problematic.
}

\textup{
2. In carrying out our nonstandard Implicit Function
construction for shock waves, we faced the problem of non-uniform
convergence of series \eqref{neuversion} due to lack
of spatial localization. 
We remedied this problem by additional weighted-norm estimates.
However, the example of another famous nonstandard Implicit
Function construction, namely, Nash--Moser iteration,
suggests the alternative approach of introducing a ``localizing''
step dual to the smoothing step in the Nash--Moser scheme.
It would be interesting to see if this approach could also be
carried out, thus avoiding the need for additional analysis
associated with pointwise bounds.
We note that
Nash--Moser iteration combining 
temporal
 with the usual frequency
cutoffs has been carried out by Klainerman \cite{Kl}
in the context of a nonlinear wave equation.
}
\end{rems}

{\bf Plan of the paper.} 
In Section \ref{framework}, we formalize the reduction
procedure set out in Section \ref{refined} as an abstract 
bifurcation framework suitable for application to general discrete
dynamical systems in the absence of spectral gap.
In Section \ref{linest}, we recall the pointwise
estimates furnished by the methods of \cite{ZH, MaZ3}
on the Green function $\tilde G$ associated with the 
transverse linearized solution operator $e^{\tilde Lt}$,
then, in Sections \ref{returnsection} and \ref{1dproof}
use these to verify that the shock bifurcation problem
indeed fits the hypotheses of our abstract framework.
This establishes a Lipschitz version of Main Theorem \ref{newPH};
we improve this to $C^1$ by a bootstrap argument in Section \ref{C1}.
In Section \ref{multidcase}, we describe the extension to 
multidimensions, verifying Theorem \ref{multinewPH}.

\medskip

{\bf Acknowledgement.}
Thanks to Claude Bardos for pointing out the
references \cite{I1, I2}, and to Bj\"orn Sandstede
and Arnd Scheel for general discussions of the
spatial dynamics method.
Thanks to Walter Craig for pointing out the reference 
\cite{Kl}, to Thierry Gallay for pointing out \cite{KuS},
and to Tudor Ratiu for pointing out \cite{BPRT}.
K.Z. thanks the Ecole Polytechnique F\'ed\'erale de Lausanne for
their hospitality in the course of a two-week visit within
their 2006 special semester on fluid dynamics,
during which much of the writing was done.  B.T. thanks Indiana University
for their hospitality during several collaborative visits in
which the main analysis was carried out.

\medskip
{\bf Notes.} 
Since the completion of this work, there have been several
further developments.
In \cite{SS3}, Sandstede and Scheel 
recover and somewhat sharpen our results
using spatial dynamics techniques, 
answering the question posed in
Section \ref{discussion} of what results may be obtained
by these methods,
obtaining the additional information of exponential localization
of solutions and exchange of spectral stability.
In \cite{TZ3, TZ4}, we extend our results
to shock and detonation waves of systems with physical viscosity,
at the same time greatly sharpening and simplifiying the basic
cancellation estimate.
The latter yields also exponential localization
and appears to shed light on the more general question of
the technical relation between spatial dynamics methods and
the ``temporal dynamics'' method used here; see Remark 1.10, \cite{TZ3}.

\section{Abstract bifurcation framework}\label{framework}

We begin by formalizing 
the approach described
in Section \ref{refined}, providing an abstract 
framework for the analysis of bifurcations in
the absence of a spectral gap for discrete dynamical systems
in the general form \eqref{Pieqns}: that is, the portion 
of our analysis occurring {\it after} the 
reduction via Poincar\'e return-map construction to a fixed-point problem.
The specifics of the return-map construction are discussed separately
in Section \ref{returnsection}.
We carry out the analysis at the level of Lipschitz regularity;
this can be improved to $C^r$ regularity by substituting everywhere
$C^r$ for Lipschitz regularity in the various assumptions.

\subsection{Generalized Lyapunov--Schmidt reduction}\label{LS}

 Given two Banach spaces $X$ and $Y,$ we denote by 
${\cal L}(X,Y)$ the space of linear continuous functions from $X$ to $Y,$
and by ${\cal L}(X)$ the space of linear continuous functions from $X$ to $X.$
 Let ${\cal B}_1,$ ${\cal B}_2,$ $X_1,$ $X_2$ be Banach spaces,
$$ X_2 \hookrightarrow X_1 \hookrightarrow \CalB_1, \qquad X_2 \hookrightarrow \CalB_2 \hookrightarrow {\cal B}_1,$$ 
with norms 
 $$\|\cdot\|_{{\cal B}_1}\le \|\cdot\|_{X_1} \le \| \cdot \|_{X_2}, \qquad \| \cdot \|_{{\cal B}_1} \le \|\cdot\|_{{\cal B}_2}\le \|\cdot\|_{X_2},$$ and for which 
the unit ball in $X_1$ is closed in $\CalB_1$.

 Consider a system of difference equations
\ba
\label{diff}
\Delta a&= f(\e, a, b)= (R(\e,a,b)- \Id) a+ N_1(\e,a,b),\\
\Delta b&= g(\e, a, b)= (S(\e,a,b)- \Id)b+ N_2(\e,a,b),\\
\ea
associated with a discrete dynamical system
\ba
\label{dyn}
\hat a&= a +f(\e, a, b)= R(\e, a,b)a+ N_1(\e, a,b),\\
\hat  b&= b+ g(\e, a, b)= S(\e, a,b)b+ N_2(\e, a,b),\\
\ea
with
$$ (f,g)(\e,0,0) \equiv 0,$$
where $\eps\in \RR^m$ is a bifurcation parameter and $R(\e, a,b)$ and $S(\e, a,b)$ are primary and transverse ``linearized''
solution operators for one time-step of \eqref{dyn},
\begin{eqnarray} \nonumber R : & (\e,a,b) \in \R^m \times \R^n \times {X}_1 & \to R(\e,a,b) \in {\cal L}(\RR^n), \\
 S : & (\e,a,b) \in \R^m \times \R^n \times {X}_1 & \to S(\e,a,b) \in {\cal L}({\cal B}_1) \cap {\cal L}(X_1), \nonumber 
 \end{eqnarray}
 and $N_1$ and $N_2$ are
 nonlinear maps
 \begin{eqnarray}\nonumber N_1: & (\e, a, b) \in \R^m \times \R^n \times {X}_1 & \to N_1(\e,a,b) \in \R^n, \\ \nonumber N_2: & (\e,a,b) \in \R^m \times \R^n \times X_1 & \to N_2(\e,a,b) \in X_2. \end{eqnarray}
 We assume that $R,$ $N_1$ and $N_2$ are continuous in $(\e,a,b),$ as maps from $\R^m \times \R^n \times {\cal B}_1$ to ${\cal L}(\R^n),$ $\R^n$ and $X_2$ respectively.  

 We assume moreover that $N_1$ and $N_2$ are of quadratic order, and 
Lipschitz, in $(a,b),$ with
  \begin{equation} \label{ass-N1}
 \begin{aligned}
 |N_1(\e,a,b)| & \le C(|a|+\|b\|_{{X}_1})^2,
 \\
  |\partial_a N_1(\e,a,b)|_{{\cal L}(\R^n, \R^n)} + |\partial_b N_1(\e,a,b)|_{{\cal L}(\CalB_1,\RR^n)}
& \le C(|a|+\|b\|_{{X}_1}), \end{aligned}
\end{equation}
and
\begin{equation} \label{ass-N2}
 \begin{aligned}\|N_2(\e,a,b)\|_{X_2} & \le C(|a|+\|b\|_{X_1})^2,  \\
  |\partial_a N_2(\e, a,b)|_{{\cal L}(\RR^n,\CalB_2)} + |\partial_b N_2(\e,a,b)|_{{\cal L}(\CalB_1,\CalB_2)}
& \le C(|a|+\|b\|_{X_1}). \end{aligned} \end{equation} 

We also assume that $R$ satisfies the 
 bound
 \begin{equation} \label{ass-R}
  | R(\e,a,b) - R(\e,0,0) |  \leq C( | a | + \| b \|_{{X}_1}),
  \end{equation}
and that $R$ is Lipschitz with respect to $a$ and $b.$

 In Sections \ref{basicred} and \ref{bifu}, 
we assume in addition the $\eps$-Lipschitz bounds 
\begin{equation} \label{lip-e-N1}
 |\d_\e N_1(\e,a,b)|_{{\cal L}(\R^m, \R^n)} \leq C( | a| + \| b \|_{{X}_1})^2,
 \end{equation}
 \begin{equation} \label{lip-e-N2}
 |\d_\e N_2(\e,a,b)|_{{\cal L}(\R^m, {\cal B}_2)} \leq C(|a| + \| b \|_{X_1})^2,
 \end{equation}
 and
 \begin{equation} \label{lip-e-R}
  | \d_\e( R(\e,a,b) - R(\e,0,0)) |  \leq C( |\e| + |a| + \| b \|_{{X}_1}).
  \end{equation}
  
We assume that the constants in 
all the above estimates
are uniform in $\e,a,b,$ for $|\e| + |a| + \| b \|_{X_1}$ small enough.

The domains of definition of $R, S, N_1$ and $N_2,$ together with the domain of validity of the above estimates, can be restricted to a neighborhood of the origin in $\R^m \times \R^n \times X_1.$

We have in mind a general discrete dynamical system whose
linearized solution operator has a finite-dimensional
subspace $a$  that is spectrally separated from its complement
$b$ at $\eps=0$, in which case form \eqref{dyn} may always be achieved,
with $R=R(\eps)$ and $S=S(\eps)$.
We allow dependence of $R$, $S$ on $(a,b)$ as well as $\eps$
in order to admit the slightly more general class of systems
\eqref{Pieqns} derived by the return map construction of 
Section \ref{return}.

The case of our particular interest is 
{\it bifurcation from a simple eigenvalue}, 
of equilibria of \eqref{dyn}, or, equivalently, zeroes of \eqref{diff},
which corresponds to
$m=n=1$, $R(\e,0,0)$ differentiable at $\e =0,$ with
 \begin{equation} \label{ass-R-0}
 R(0,0,0)= 1,  \qquad \partial_\eps R(0,0,0)\ne 0.
 \end{equation}

\br\label{strongweak}
\textup{
The reason for introducing the different norms $X$ and ${\cal B}$
is that we shall ultimately carry out a fixed-point iteration
using the standard extension of the Contraction-mapping Principle
to the case of a map $\CalT$ that is bounded in norm $\|\cdot\|_{X_1}$ and
contractive in the weaker norm $\|\cdot\|_{{\cal B}_1}$ for bounded
$\|\cdot\|_{X_1}$, with the unit ball in $X_1$ closed in ${\cal B}_1$. 
This gives additional flexibility that is useful in applications.
}
\er

\br\label{notheta}
\textup{
In the return map construction of Section \ref{return},
it is more standard \cite{HK, TZ} to 
eliminate $t$ instead by 
using
$dt/d\theta= \tau(\eps)^{-1} + O(r)$
to rewrite \eqref{PHeq1} with ``$'$'' denoting $d/d\theta$ as
\ba\label{rthetaeqns2}
r'&= (\gamma/\tau)(\eps)r + N_1(r, \theta, v,\eps),\\
v'&= (\tilde L/\tau)(\eps)v + N_2(r, \theta, v,\eps),\\
\ea
leading to a return map \eqref{Pieqns} of form \eqref{dyn}
with $R$, $S$ depending only on $\eps$. However,
the nonlinearity
$$ 
N_2= \big((\tau+ N_\theta)^{-1}- \tau\big)\tilde L v +
(\tau+ N_\theta)^{-1}N_v
=O(r)\tilde Lv + O(N_v)
$$
in \eqref{rthetaeqns2} contains terms $\partial_x (rv)$ that do not
take $v\in L^2$ to $\partial_x L^1$, to which the
analysis of Section \ref{refined} does not apply.
It is for this technical reason that we allow $(a,b)$-dependence
of $R$, $S$, 
gaining needed flexibility
at the expense of slight notational and expositional inconvenience.
The quasilinear form of $S$ may also be useful in more general circumstances, 
allowing more general iteration schemes in situations of limited regularity.
}
\er

%

\br\label{Rcond}
\textup{
In the Poincar\'e--Hopf context (see Section \ref{returnsection}), 
$\e, a \in \R$ and 
$R(\e,a,b)= e^{\gamma(\eps)T(\e,a,b)}$
with $\gamma\in C^k$, $T$ Lipschitz, $T(0) \neq 0,$ $\gamma(0)=0$, $\gamma'(0) \neq 0,$ so that \eqref{ass-R}, \eqref{ass-R-0} and \eqref{lip-e-R} hold. 
}
\er

\subsubsection{The equilibrium problem}\label{equilibrium}
We seek to solve the equilibrium problem $(f,g)=(0,0)$ by
Lyapunov--Schmidt reduction, that is, to determine a function
$b=B(\e, a)$ satisfying
$$
g(\e, a, B(\e,a))\equiv 0.
$$
Substituting into \eqref{diff}, we would then
obtain a {\it reduced}, finite-dimensional equilibrium problem
(on nullcline, rather than center manifold)
$$
0=f^*(\e, a):=f(\e, a, B(\e, a)),
$$
presumably amenable to analysis by standard finite-dimensional techniques.

\subsubsection{Assumptions}\label{assumptions}

 We are interested in the (non-standard) case that $\mbox{Id} - S(\e,a,b)$ does not possess a bounded inverse in ${\cal L}({\cal B}_1).$ Our key assumption is as follows.
  
  \begin{ass} \label{B'''} For every $(\e,a,b) \in {\cal V},$ a neighborhood of the origin in $\R^m_\e \times \R^n_a \times X_1,$ the linear operator $\Id - S(\e,a,b) \in {\cal L}(X_1)$ possesses a right inverse $$(\Id - S(\e,a,b))^{-1}:  X_2 \to {\cal B}_1,$$
  that is bounded in ${\cal L}({\cal B}_2, {\cal B}_1)$ and ${\cal L}({X}_2, {X}_1)$ norms.
  \end{ass}
 
 Note that the range of 
$N_2$
 is contained in $X_2,$ the subspace of ${\cal B}_2$ where the right inverse is assumed to be defined.  
 
 \begin{prop} 
 \label{B-to-A} For Assumption {\rm \ref{B'''}} to hold, it is sufficient that, for all $(\e,a,b) \in {\cal V},$ the sequence $\sum_{j=0}^{N}S(\e,a,b)^j$ be conditionally convergent in ${\cal L}(X_2,{\cal B}_1)$ norm, and the limit operator $\sum_{j=0}^{\infty}S(\e,a,b)^j$ be bounded in ${\cal L}({\cal B}_2, {\cal B}_1)$ and ${\cal L}({X}_2,{X}_1)$ norms,    
in which case $(\Id -S(\e,a,b))^{-1}=\sum_{j=0}^\infty S(\e,a,b)^j$.
\end{prop}
 
 \begin{proof} Given $x \in X_2,$ let 
 \begin{equation} \label{inv-def} (I-S(\e,a,b))^{-1} x : = \lim_{N\to \infty} \sum_{j=0}^N S^j(\e,a,b) x,\end{equation}
where the limit is taken in ${\cal B}_1$. By a standard telescoping sums argument, for all $N,$ 
\begin{equation} \label{telescop}
 (\Id-S(\e,a,b)) \sum_{j=0}^N S^j(\e,a,b) x = x - S^{N+1}(\e,a,b) x. \end{equation}
Because $S(\e,a,b) \in {\cal L}({\cal B}_1),$ 
 $$ (\Id-S(\e,a,b))  \lim_{N \to \infty} \sum_{j=0}^N S^j(\e,a,b) x = \lim_{N \to \infty} (\Id-S(\e,a,b)) \sum_{j=0}^N S^j(\e,a,b) x . $$
{}From the above equality and \eqref{telescop},
 $$ (\Id-S(\e,a,b))(\Id-S(\e,a,b))^{-1}x = \lim_{N\to \infty} x - S^{N+1} (\e,a,b) x = x,$$
as (conditional)
convergence of $\sum_{j=0}^\infty S^j(\e,a,b) x$ implies that $S^j(\e,a,b) x$ converges to 0. Thus, the operator defined by \eqref{inv-def} is a right inverse of $\Id - S.$ It is bounded in ${\cal L}({\cal B}_2, {\cal B}_1)$ and ${\cal L}({X}_2, {X}_1)$ norms by assumption. 
\end{proof}

We will also need the following regularity assumption.

\begin{ass} \label{B'''bis} 
 For $\e,a,b$ in ${\cal V},$ the right inverse $(\Id - S)^{-1}$ 
given by Assumption {\rm \ref{B'''}} is Lipschitz continuous 
with respect to $b$ as a function from ${\cal B}_1$ to
${\cal L}({\cal B}_2, {\cal B}_1)$
and with respect to $(\eps,a)$ as a function from $\R^{m+n}$ to
${\cal L}(X_2, {\cal B}_1)$, in the sense that
\ba\label{lipdisplay}
   \sup_{\e,a,b \in {\cal V}}
\|\partial_{b} (\Id - S)^{-1}\|_{
{\cal L}({\cal B}_1, 
{\cal L}({\cal B}_2, {\cal B}_1))}
&< \infty,\\
   \sup_{\e,a,b \in {\cal V}}
\|\partial_{\eps, a} (\Id - S)^{-1}\|_{ {\cal L}(\R^m\times \R^n, 
{\cal L}(X_2, {\cal B}_1))}
&< \infty,\\
\ea
satisfying also
  \ba  \label{cont-0}
   \sup_{\e,a,b \in {\cal V}} \| (\Id - S)^{-1}(\e,a,b) \|_{{\cal L}(X_2,X_1)} 
&< \infty,\\
   \sup_{\e,a,b \in {\cal V}} \| (\Id - S)^{-1}(\e,a,b) 
\|_{{\cal L}(\CalB_2,\CalB_1)} &< \infty.
   \ea
 \end{ass}
  
\br\label{expdecay}
\textup{
Convergence in Proposition \ref{B-to-A} is implied by exponential decay of $\|S^n\|_{{\cal L}({\cal B}_2, {\cal B}_1)}$ as
$n\to \infty$, which yields absolute convergence by the ratio test.
This in turn is implied by stability plus spectral gap of
the transverse operator $S$, i.e., $|\sigma (S)|<1$, the
standard case.
In our case of interest, there is no spectral gap 
and $\|S^n\|_{{\cal L}(\CalB_2,\CalB_1)}$, or equivalently
$\|S^n N_2\|_{{\cal L}(\CalB_1)}$, 
is not exponentially but algebraically decaying, while
$\|S^n\|_{{\cal L}(\CalB_1)}$ is not decaying at all.  Thus, 
convergence will in general be conditional if it occurs, 
and also depends strongly on the specific structure of $N_2$. 
In the application to viscous shock waves, moreover,
Lipschitz continuity, \eqref{lipdisplay},
holds for the limit, but not the approximating finite sums;
see Section \ref{revisit}.
}
\er

\subsubsection{Reduction}\label{basicred}

 \begin{lem} \label{g+} Under Assumption {\rm \ref{B'''}}, the equation 
 \begin{equation} \label{b} g(\e,a,b) = 0, \qquad (\e,a,b) \in \R^m \times \R^n \times X_1,\end{equation}
 is equivalent to
 \begin{equation} \label{equiv'}
 b - (\Id -S(\e, a,b))^{-1}N_2(\e, a,b) \in \kernel(\Id -S(\e,a,b)) \cap X_1.
 \end{equation} 
 \end{lem}
 
 \begin{proof} 
 Applying to the left of \eqref{B'''} the right inverse $(\Id-S(\e,a,b))^{-1}$ given by Assumption \ref{B'''}, we obtain
$$
 \tilde b := (\Id -S(\e,a,b))^{-1}(\Id -S(\e,a,b)) b= (\Id -S(\e,a,b))^{-1}N_2(\e,a,b).
$$
Observing that $\tilde b - b$ belongs to $\kernel(\Id -S(\e,a,b)) \cap X_1$,
we obtain \eqref{equiv'}. Conversely, \eqref{equiv'} implies $g(\e,a,b) = 0,$ by definition of 
$(\Id - S(\e,a,b))^{-1}.$
\end{proof}

 Consider the map ${\cal T}: \R^m \times \R^n \times X_1 \times X_1 \to X_1,$ defined by,
 $$ {\cal T}(\e,a,b,\o) = \o + (\Id -S(\e, a,b))^{-1}N_2(\e, a,b).$$
 
 \begin{prop}\label{red'}
 Suppose that $N_2$ satisfies {\rm \eqref{ass-N2}} and {\rm \eqref{lip-e-N2},} and that Assumptions {\rm \ref{B'''}} and {\rm \ref{B'''bis}} are satisfied. Then for
  $|\e| + |a| + \|\omega\|_{X_1}$ sufficiently small, ${\cal T}$ has a unique fixed point $b=B(\e,a,\omega) \in X_1$, $B(\e,0,0) = 0,$
such that $B$ is Lipschitz in 
$(\e, a, \omega)$ with respect to norm $\|\cdot\|_{\CalB_1}$, with
\ba\label{Bybd}
\|B(\e,a,\omega)\|_{{\cal B}_1} &\le C(\|\omega\|_{{\cal B}_1} + |a|^2), \\
\|B(\e,a,\omega)\|_{X_1} &\le C(\|\omega\|_{X_1} + |a|^2), 
\ea
and,
\begin{equation} \label{Bybd2}
\| \d_\e B(\e,a,\omega)\|_{{\cal L}(\R^n,{\cal B}_1)} \le C(\|\omega\|_{X_1}^2
+ |a|^2).
\end{equation}
\end{prop}

\begin{proof} 
By \eqref{cont-0}(i) together with \eqref{ass-N2}(i), 
$$
\|\CalT(\e, a, b, \o)\|_{X_1}\le C(\|\omega\|_{X_1}+\|b\|_{X_1}^2 + |a|^2)
$$
for $|\e| + |a| + |\o|_{X_1}$ sufficiently small.
Thus, for $\e,a,\o$ fixed and $|\e| + |a| + |\o|_{X_1}$ small enough, 
${\cal T}$ maps a small ball $B(0, 2C(|\o\|_{X_1}+|a|))$ 
of $X_1$ centered at the origin into itself.

  To prove that ${\cal T}$ is a contraction, we evaluate $\|{\cal T}(b) - {\cal T}(b')\|_{{\cal B}_1},$ for small $b, b' \in X_1.$ This term can be bounded by the sum of 
\begin{equation} \label{m1} 
 \| (\Id - S)^{-1}(b) - (\Id - S)^{-1}(b') \|_{{\cal L}({\cal B}_2, {\cal B}_1)} \| N_2(b') \|_{{\cal B}_2},
 \end{equation}
 and
 \begin{equation} \label{m2}
 \| (\Id - S)^{-1}(b') \|_{{\cal L}({\cal B}_2, {\cal B}_1)} \| N_2(b) - N_2(b') \|_{{\cal B}_2}.
 \end{equation}
 By Lipschitz continuity (Assumption \ref{B'''}) and \eqref{ass-N2}(i), 
for $(\e,a,b) \in {\cal V}$ and $(\e,a,b') \in {\cal V},$ 
the term in \eqref{m1} is bounded by 
$C (|a| + \| b'\|_{X_1})^2 \|b - b'\|_{{\cal B}_1},$ 
where $C$ does not depend on $\e,a,b,b'.$ 
By \eqref{cont-0}(ii) and \eqref{ass-N2}(ii), 
for $|\e| + |a| + \|b\|_{X_1} + \|b'\|_{X_1}$ small enough, 
the term in \eqref{m2} is bounded by 
$C(|a| + \|b \|_{X_1} + \| b'\|_{X_1}) \| b - b' \|_{{\cal B}_1},$ 
where $C$ does not depend on $\e,a,b,b'.$

Combining, we find that ${\cal T}$ is contractive with respect 
to $b$ in the $\CalB_1$ norm, for $(\e,a,b)$ in a neighborhood of the 
origin in $\R^m \times \R^n \times X_1,$ possibly smaller than ${\cal V},$ 
and for $\| \o \|_{X_1}$ small enough. 
By the extension of the Contraction-mapping Principle 
described in Remark \ref{strongweak},
we may thus conclude the existence of a unique fixed point 
$b=B(\e,a,\o)\in X_1$ of ${\cal T}$ in a small ball of $X_1$ centered at the origin; 
specifically, considering the ball in $X_1$ as a $\CalT$-invariant closed
set in ${\cal B}_1$, we may apply the standard Contraction-mapping 
Principle with respect to $\|\cdot\|_{{\cal B}_1}$.
(Recall our assumption that the ball in $X_1$ be closed in ${\cal B}_1$.)
The uniqueness of that fixed point, together with the identity ${\cal T}(\e,0,0,0) = 0,$ imply that $B(\e,0,0) = 0,$ for $|\e|$ small.

By \eqref{ass-N2}(i)--(ii), we have
 $$
\begin{aligned}
 \| B \|_{\CalB_1}&\le \|\omega\|_{\CalB_1}+ C(|a|^2 + \|B\|_{\CalB_1}^2),\\
 \| B \|_{X_1}&\le \|\omega\|_{X_1}+ C(|a|^2 + \|B\|_{X_1}^2),\\
\end{aligned}
$$ 
from which, recalling by our fixed point construction that $B$ is
confined to a small ball in $X_1$, hence has small ${\cal B}_1$ and $X_1$ norms,
we obtain \eqref{Bybd}(i)--(ii). 
Likewise, smallness of $B$ together with \eqref{cont-0}(ii) gives
  \begin{equation} \label{Bcont-0}
   \sup_{\e,a,b \in {\cal V}} \| (\Id - S)^{-1}(\e,a,B(\e,a)) 
\|_{{\cal L}({\cal B}_2,{\cal B}_1)} < \infty.
   \end{equation}

 To establish Lipschitz regularity in $\e$, 
we evaluate 
 \begin{equation} \label{1-2.1} \begin{aligned}
  \| \d_\e B \|_{{\cal L}(\R^n,{\cal B}_1)} & \leq \| (\d_\e {\cal T})(\e, B(\e)) \|_{{\cal L}(\R^n,{\cal B}_1)} \\ &  + \| (\d_b {\cal T})(\e, B(\e)) \|_{{\cal L}({\cal B}_1)} \| \d_\e B \|_{{\cal L}(\R^n,{\cal B}_1)}. \end{aligned} 
  \end{equation} 
 The second term in the right-hand side of \eqref{1-2.1} is controlled by contractivity of ${\cal T}$ 
with respect to $b,$ 
hence bounded by $\theta  \| \d_\e B \|_{{\cal L}(\R^n,{\cal B}_1)}$ for $0<\theta<1$. The first term in the right-hand side of \eqref{1-2.1} is controlled by the sum of
  \begin{equation} \label{m4}  \| (\d_\e (\Id - S)^{-1})(\e, B(\e)) 
\|_{{\cal L}(\R^n,{\cal L}(X_2, {\cal B}_1))} \| 
N_2(\e, B(\e)) \|_{X_2}, \end{equation}
    and
    \begin{equation} \label{m5}
   \| (\Id - S)^{-1}(\e,B(\e)) \|_{{\cal L}({\cal B}_2, {\cal B}_1)} \| 
(\d_\e N_2)(\e, B(\e)) \|_{{\cal L}(\R^n,{\cal B}_2)}.
   \end{equation}
By Lipschitz regularity (with respect to $\e$) of $(\Id - S)^{-1}$ in 
${\cal L}(X_2, {\cal B}_1)$ norm, 
\eqref{Bybd}(ii),  
and bound \eqref{ass-N2}(i) on $N_2,$ 
the term \eqref{m4} is bounded by 
$ C(|a|^2 + \|\o\|_{X_1}^2).  $
Finally, 
by \eqref{Bcont-0} together with the Lipschitz bound \eqref{lip-e-N2} 
on $N_2,$ 
term \eqref{m5} is bounded by
$ C
(|a|^2 + \|\o\|_{X_1}^2)$ as well.
Combining these bounds 
 proves Lipschitz regularity of $B$ with respect to $\e,$ 
with Lipschitz constant 
$$
2 C (|a|^2 + \|\o\|_{X_1}^2)/(1-\theta).
$$
Summing all contributions gives \eqref{Bybd2}.
Using the facts that $(\Id - S)^{-1}$ and $N_2$ are Lipschitz 
with respect to $a, b$, 
and ${\cal T}$ is contractive with respect to $b,$ we can prove similarly that $B$ is Lipschitz with respect to $a$ and $\o$ as well.
\end{proof}

Noting that $0\in \kernel (\Id -S(\e,a,b))$ for any $(\e,a,b)$,
we may choose $\omega \equiv 0.$ Then $$g(\e,a,B(\e,a,0)) = 0,$$
  for $\e,a$ small. In general, this is the only reasonable choice, since
the rest of the kernel may depend on the solution $b$. It follows from Lemma \ref{g+} that if $\o(\e,a) \in \cap_{b} \kernel (\Id - S(\e,a,b)),$ for all $\e,a,$ and such that $\o$ is small enough for $\e,a$ small, then $g(\e,a,B(\e,a,\o(\e,a))) = 0,$ for all $\e,a$ small.   

\br\label{restricted}
\textup{
If $\Id-S$ has a non-trivial kernel, as holds always in the traveling-wave context
by translational group invariance, then, in order for $\Id - S$ to have a right inverse, 
$\Range(N_2)$ must lie in a complementary subspace to $\kernel (\Id - S),$ which
is preserved by $S$.
If also Assumption \ref{B'''} holds, then by the
choices $(I-S)^{-1}=\sum S^j$, $\omega=0$ we 
automatically project out the kernel and work on 
the (nonlinearly) invariant complementary subspace,
at least in the case (as for finite-dimensional ODE, or, more generally,
when zero is an isolated eigenvalue)
that the complementary subspace is closed with respect to $\|\cdot\|_{\CalB_1}$.
%
Note, however, that we do not a priori (i.e., by force)
restrict to invariant subspaces, but find this in a natural
way through the analysis.
And, in the case that the complementary subspace is {\it not} closed 
with respect to $\|\cdot\|_{\CalB_1}$, we may still obtain a result,
even though this result may lie outside the complementary subspace.
}

\textup{
For example, in the shock wave context, $\Range (N_2)$ is zero-mass, 
preserved by $S$, hence (recall, zero-eigenfunction $\bar u'$ 
has nonzero mass $u_+-u_-$) the set of zero-mass perturbations is a 
complementary $S$-invariant subspace to $\kernel (I-S)$.
However, it is not closed with respect to $\CalB_1=\|\cdot\|_{L^2}$,
and our reduction procedure (apparently) takes us out of this subspace;
see Remark \ref{zeromass}.
}
\er

\br\label{stabcase}
\textup{
In the case that $S$ is not stable, it may be split into
stable/unstable parts and Lyapunov--Schmidt reduction carried out
by a combined forward/backward scheme.  
In the parabolic PDE context,  there can exist only
finitely many unstable modes, hence the backward unstable flow
is well-posed and this algorithm may indeed be carried out.
}
\er

\subsection{Bifurcation} \label{bifu}

 Extending the discussion at the end of the previous section, 
 consider a Lipschitz
(with respect to $X_1$-norm) 
curve $\bar \o,$ defined on a neighborhood of the origin 
in $\R^m_\e \times \R^n_a,$ such that $\bar \o(\e,0) \equiv 0$ and
  \begin{equation} \label{omega} (\e,a) \mapsto \bar \o(\e,a) \in \bigcap_{\| b \|_{X_1} \leq b_0} \kernel (\Id - S(\e,a,b)) \cap X_1,\end{equation} 
for some $b_0 > 0.$ The trivial choice $\bar \omega\equiv 0$ shows that at least one such
curve always exists.
  
 \begin{prop}\label{lipenough}
Under the assumptions of Proposition {\rm \ref{red'}}, suppose in addition that $N_1$ satisfies {\rm \eqref{ass-N1}} and {\rm \eqref{lip-e-N1}}, that $R$ satisfies {\rm \eqref{ass-R}, \eqref{ass-R-0}} and {\rm \eqref{lip-e-R}}, 
and $a, \e\in \RR.$ 
Then, for each $\bar \omega$ as in \eqref{omega},
there exists a Lipschitz map $a \mapsto \e_{\bar \o}(a),$ defined on a ball 
of radius $r$ centered at the origin in $\R_a,$ 
such that $\e_{\bar \o}(0) = 0,$ and
 $$ (f,g)(\e_{\bar \o}(a), a, B(\e_{\bar \o}(a), a, \bar \o(\e_{\bar \o}(a), a))) \equiv 0,$$
where $B$ is the Lipschitz map given by Proposition {\rm \ref{red'}}.
Moreover, for $|a| \leq r,$
$\eps=\eps_{\bar \o}(a)$ is the unique solution of 
$
(f,g)(\cdot,a,B(\cdot, a, \bar \o(\cdot, a))=(0,0). 
$
Here, the size of $r$ depends in particular
 on the size of the 
Lipschitz constant of $\bar \omega$.
\end{prop}  

\begin{proof} There exists $\kappa > 0,$ such that, for $|\e| + |a|$ small, $\| \bar \o \|_{X_1} \leq \kappa |a|.$ Consider the map $(\e,a) \mapsto \bar B := B(\e,a,\bar \o(\e,a)),$ where $B$ is given by Proposition \ref{red'}. The map $\bar B$ is defined for $|\e| + |a|$ small, and $g(\e, a, \bar B) \equiv 0.$ Besides, from \eqref{Bybd}, the bound on $\bar \o,$ and \eqref{ass-N1}(i), there holds the bound $| N_1(\e,a,\bar B) | \leq C(\kappa) |a|^2,$ and 
 $$ f(\e,a,\bar B) = (R(\e,0,0) - 1) a + (R(\e,a,\bar B) - R(\e,0,0)) a + N_1(\e,a,\bar B),$$
 (recall that $R(0,0,0) = 1$). Let
 $$ \bar N_1 := R(\e,a,\bar B) - R(\e,0,0) + N_1(\e,a,\bar B)/a, \qquad \bar N_1(a = 0) := 0.$$
{} From \eqref{ass-R} and the above bound on $N_1(\bar B),$ we see that $|\bar N_1| \leq C |a|,$ for $\e,a$ small. The bounds \eqref{ass-N1}, \eqref{lip-e-N1} and \eqref{Bybd}-\eqref{Bybd2} imply that $N_1(\bar B)/a$ is Lipschitz in $\e,$ with Lipschitz semi-norm controlled by $|a|.$ The bounds \eqref{lip-e-R}, \eqref{Bybd} and \eqref{Bybd2}
 imply that $R(\e,a,\bar B) - R(\e,0,0)$ is also Lipschitz in $\e,$ with Lipschitz semi-norm controlled by $|\e| + |a|.$ Thus, $\bar N_1$ is Lipschitz in $\e,$ with Lipschitz semi-norm controlled by $|\e| + |a|.$ Similarly, the bounds \eqref{ass-N1}(ii), \eqref{Bybd}-\eqref{Bybd2} and \eqref{ass-R} imply that $\bar N_1$ is Lipschitz in $a$ as well. Let now  
 $$ \bar f := \gamma^{-1} f(\bar B)/a = \e + \gamma^{-1} (R(\e,0,0) - 1 - \gamma \e + \bar N_1),$$
 where $\gamma := \d_\e R(0,0,0) \neq 0$ by \eqref{ass-R-0}. Because $R(\e,0,0)$ is assumed to be differentiable at the origin, 
 \begin{equation} \label{small-o} R(\e,0,0) - 1 - \gamma \e = o(|\e|).\end{equation}
 Thus, for $|a|$ small enough, the map $\bar N := \gamma^{-1} (R(\e,0,0) - 1 - \gamma \e + \bar N_1)$ maps a small ball in $\R^m_\e$ to itself. The above description of $\bar N_1$ and \eqref{small-o} imply that $\bar N$ is contractive in $\e,$ for $\e,a$ small enough.  The Contraction-mapping argument used in Proposition \ref{red'} then implies the existence and uniqueness of a fixed point of $- \bar N,$ in a small ball centered at the origin in $\R^m_\e,$ for $|a|$ small. This fixed point $a \mapsto \e_{\bar \o}(a)$ is the desired solution to $f = 0.$ Indeed, it is Lipschitz in $a,$ because $\bar N$ is, and $\e_{\bar \o}(0) = 0$ by uniqueness. 
 \end{proof}

Proposition \ref{lipenough} concerns the general case.
Let us now consider a more specific situation for which
 $\kernel(\Id -S(\e,a,b))$ has a uniform structure.
 
 \begin{ass}\label{ker} The Kernel $\kernel(\Id -S(\e,a,b))$ is independent 
of $b$, with a locally Lipschitz basis $v_1(\eps, a), \dots, v_\ell(\eps,a)$ 
of fixed dimension $\ell$.
 \end{ass}

Under Assumption \ref{ker}, the maps
$$
(\sigma,\e,a) \mapsto \bar \omega^\sigma(\eps, a):= |a|\sum_{j=1}^\ell \sigma_j v_j(\eps, a),
$$
are uniformly Lipschitz in $(\eps,a)$ for $\sigma$ in 
$B(0,R)\subset \RR^\ell$ and $|\e| + |a|$ small, with combined range covering the cone 
$\|\omega\|_{X_1} \le C|a|$, any $C>0$, for $R>0$ sufficiently large.

\begin{cor}\label{cone} Under the assumptions of Proposition {\rm \ref{lipenough}}, suppose further that Assumption {\rm \ref{ker}} is satisfied.
Then, for each $|a|$ small enough and each $\sigma\in B(0,R) \subset \R^{\ell}$, in a neighborhood of the origin in 
$\R_\eps\times X_1,$ 
there exist a unique $(\e,b)$ such that 
 $$ \left\{ \begin{aligned} (f,g)(\e,  a, b) & = 0, \\ b - (\Id - S(\e,a,b))^{-1} N_2(\e,a,b) & = \bar \o^\sigma(\e,a).\end{aligned}\right.$$
 This equilibrium is given by $\e = \e_{\bar \o^\sigma}(a),$ $b = B(\e_{\bar \o^\sigma},a,\bar \o^\sigma),$ with the notations of Propositions {\rm \ref{red'}} and {\rm \ref{lipenough}}. Moreover, all solutions of $(f,g)=0$ with $\|b\|_{X_1}\le C|a|$
and $(\eps,a,b)$ sufficiently small are of this type,
apart from the trivial solution $(a,b)\equiv(0,0)$.
\end{cor}

\begin{proof}
Only the final assertion requires discussion.
Evidently, solutions with $\|\omega\|_{X_1}\le C|a|$,
$\omega:=b-(\Id -S(\eps, a, b))^{-1}N_2(\eps,a,b)$,
are of the described type.
Thus, we have only to observe that
$\|\omega\|_{X_1}=\|b\|_{X_1}+ O(\|b\|_{X_1}+|a|)^2$, by the bounds on $N_2$ (\eqref{ass-N2}(i)),
so that $\|b\|_{X_1} \le C|a|$ implies $\|\omega\|_{X_1}\le C_1|a|.$
\end{proof}

 That is, there occurs a ``restricted'' $\ell$-fold Lipschitz 
bifurcation of equilibria  
at $\eps=0$ from the trivial solution $(a,b)\equiv (0,0)$,
on a cone of expanding radius proportional to $|a|$ around
the base solution corresponding to $\bar \omega\equiv 0$.

\subsubsection{Group invariance}

The main situation that we have in mind is the case that
 $\kernel(\Id-S(\e,a,b))$ corresponds to group-invariance(s) 
of the underlying equations, so that, in particular, the kernel depends only on $\eps$. 

Let ${\cal E}^\e$ be the set of equilibria $(a,b) \in \R^n_a \cap X_1$ 
associated with the value $\e$ of the bifurcation parameter. 
Given $(a,b) \in {\cal E}^\e,$ let $\| (a,b) \| := |a| + \| b \|_{X_1}.$ 
 
\begin{ass} \label{group2} There exists an $\ell$-dimensional 
Lie Group $G,$ with Lie algebra ${\cal G},$ and $\e_0 > 0,$ such that, for $|\e| \leq \e_0,$ there exists a family of group actions
 $$ \Phi^\e: \quad G \times {\cal E}^\e \to {\cal E}^\e,$$
one-to-one in $g\in G$,
satisfying the following properties, where the identity element in $G$ is denoted by $e_G$:
 \begin{itemize}
 \item[{\rm (i)}] The map $\Phi^\e(\cdot,0,0)$ is differentiable at $e_G,$ uniformly in $\e,$ in the 
sense that, for all $\g \in C^1: \Omega \to G,$ where $\Omega$ is a neighborhood of the origin in $\R^\ell,$
  $$ \sup_{|\e| \leq \e_0} \frac{1}{|x|} \Big\| \Phi^\e(\g(x),0,0) - \d_1 \Phi^\e(e_G,0,0) \cdot (\g'(0) \cdot x) \Big\| \longrightarrow 0, \quad \mbox{as $x \to 0.$}$$ 
  \item[{\rm (ii)}] The derivative $\d_1 \Phi^\e(e_G,0,0)$ is an isomorphism
 $${\cal G}  \quad \tilde \longrightarrow \quad \{ 0_{\R^n_a} \} \times \big( \kernel (\Id - S(\e,a,b)) \cap X_1 \big),$$ and the norms of  $(\d_1 \Phi^\e(e_G,0,0))^{-1}$ are bounded uniformly in $\e,$ for $|\e| \leq \e_0.$
  \item[{\rm (iii)}] For $g$ in a neighborhood of $e_G,$ the map $\Phi^\e(g,\cdot,\cdot)$ is Lipschitz in $(a,b),$ for small $\| (a,b) \|,$ with Lipschitz semi-norm bounded uniformly in $\e,$ for $|\e| \leq \e_0.$  
\item[{\rm (iv)}] 
For $\|(a,b)\|$ sufficiently small,
the map $\Phi^\e(\cdot, a, b)$ is continuous in $g$ for $g$ in a neighborhood
of $e_G$.
  \end{itemize}
   \end{ass}

 \begin{prop} \label{group-inv} Suppose that $N_2$ satisfies {\rm \eqref{ass-N2}(i)}, that Assumptions {\rm \ref{B'''}} and {\rm \ref{group2}} are satisfied, and moreover that the right inverse satisfies bound {\rm \eqref{cont-0}(i).} Then, for any $C > 0,$ there exists $\rho > 0,$ such that, for each equilibrium $(a,b) \in {\cal E}^\e$ such that $\| (a,b) \| < \rho$ and $|\e| \leq \e_0,$ there exists $g \in G,$ such that $(\hat a, \hat b) := \Phi^\e(g)(a,b)$ satisfies $\| \hat b \|_{X_1} \leq C \| \hat a \|.$
 \end{prop}

 \begin{proof} 
We first observe, given 
$(a_0,b_0) \in {\cal E}^\e$ with $\|(a_0,b_0)\|$ sufficiently small, and $|\e| \leq \e_0,$
that there exists $g_0$ such that 
$ \| \Phi^\e(g_0,a,b) \|$ minimizes $ \| \Phi^\e(g,a,b) \|$
for $g \in G $ lying in a fixed neighborhood $\CalN$ of $e_G$.
For, by (i)--(ii), $\|\Phi^\e(g, 0,0)\|$ is bounded from
below by some $c_0>0$ for $g$ outside a smaller neighborhood $\CalN'$ of $e_G$,
whereas, by (iii),
$\| \Phi^\e(g,a,b)-\Phi^\e(g,0,0) \| \le C\|(a,b)\|$, without loss of generality
$C>1$.
Thus, for $\|(a,b)\|\le c_0/2C$, $\|\Phi^\e(g,a,b)\|\ge c_0/2> \|(a,b)\|$
for $g$ outside $\CalN'$, and so, by continuity, (iv), there exists a minimal
$g_0\in \CalN'$.

Let $(a_0,b_0) \in {\cal E}^\e,$ and let $(a,b)$ be 
such 
a minimizing equibrium,
with $(a,b) = \Phi^\e(g_0)(a_0, b_0)$ for some $g_0 \in G.$  Assume by way of contradiction that $|a| =  c_1 \|b\|_{X_1}$
for $c_1>0$ sufficiently small.
{}From $g(\eps, a, b)=0$, by Lemma \ref{g+}, we have 
$$
b-(\Id -S(\eps,a,b))^{-1}N_2(\eps, a,b)=\omega \in \kernel (\Id -S(\eps,a,b)) \cap X_1, $$
where, by \eqref{ass-N2}(i) and \eqref{cont-0}(i),
\be\label{smalloh}
\|(\Id -S(\eps,a,b))^{-1}N_2(\eps, a,b)\|_{X_1} \le C(|a|+\|b\|_{X_1})^2,
\ee
uniformly in $\e.$ 
From \eqref{smalloh}, we deduce that, if $\| b \|_{X_1}$ is small enough,
 $$(a,b) - (0,\o)  = o(\| b \|_{X_1}).$$
 By (ii), associated with $\o$ is $\delta^\e \in {\cal G},$ such that $ \d_1 \Phi^\e(e_G,0,0) \cdot \delta^\e = (0,\omega).$ By (i)-(ii), there exists $g_\e \in G,$ such that 
  $$ \Phi^\e(g_\e,0,0) - (0,\omega) = o(\| \omega \|_{X_1}),$$
 hence 
 $$ \Phi^\e(g_\e,0,0) - (a,b) =o(\| b \|_{X_1}),
$$
uniformly in $\e.$ 
But, then, Lipschitz continuity, (iii), together with the group property,
gives
$$ \begin{aligned} 
\Phi^\eps(g_\e^{-1}, a, b)&=
\Phi^\eps(g_\e^{-1}, \Phi^\eps(g_\e, 0,0) + o(\| b \|_{X_1}))\\
&=\Phi^\eps(g_\e^{-1}, \Phi^\eps(g_\e, 0,0)) + o( \| b \|_{X_1})\\
&= o( \| b \|_{X_1}),
\end{aligned} $$
uniformly in $\e,$ thus contradicting minimality of $(a,b).$
 \end{proof}

\begin{cor}\label{group-bif}
Under the assumptions of Proposition {\rm \ref{lipenough}},
suppose further that Assumption {\rm \ref{group2}} is satisfied, and that $\d_1 \Phi^\e(e_G,0,0)$ is Lipschitz in $\e.$
Then, for $(\eps, a, b)$ sufficiently small,
the solution $(\e_0(a), a,  B(\e_0(a), a, 0))$ 
determined by $\o \equiv 0$ in Proposition {\rm \ref{lipenough}} is the
unique solution of $(f,g)=0$ up to group invariance,
apart from the trivial solution $(a,b)\equiv(0,0)$.
\end{cor}

\begin{proof} Assumption \ref{group2} implies in particular that $\kernel (\Id - 
S(\e,a,b))$ is independent of $b.$ If $v_1, \dots, v_\ell$ is a basis of ${\cal G},$ then 
$$\d_1 \Phi^\e(e_G,0,0) \cdot v_1, \dots, \d_1 \Phi^\e(e_G,0,0) \cdot v_\ell,$$ is a 
Lipschitz basis of $\kernel (\Id - S(\e,a,b)).$ Thus, Assumption \ref{ker} is satisfied, 
and Corollary \ref{cone} implies that the solution $\e_0, b_0,$ where $b_0(x) = 
B(\e_0(a), a, 0),$ and $\e_0$ is given by $\o \equiv 0$ in Proposition \ref{lipenough}, 
is unique up to group
invariance, within the cone $\|b\|_{X_1}\le C|a|,$ for any $C > 0.$
For, by equivalence of dimension (Invariance of Domain Theorem), for fixed $a$, the 
manifold of equilibria constructed
in Corollary \ref{cone}, indexed by $\sigma\in \R^\ell$,  must agree locally with the 
continuous
manifold $\{\Phi^{\eps_0(a)}(g,  a,B(\e_0(a), 0))\}$, indexed by $g\in \R^\ell$, since 
uniqueness implies that the second is contained in the first within the cone in question.
Finally, Proposition \ref{group-inv} states that each equilibrium has a representative in 
the cone $\|b\|_{X_1}\le C|a|,$ giving uniqueness up to group invariance with no 
restriction.
\end{proof}

 That is, there occurs a complete $\ell$-fold Lipschitz 
bifurcation of equilibria  
at $\eps=0$ from the trivial solution $(a,b)\equiv (0,0)$
in a neighborhood of $(\eps, a, b)=(0,0,0)$.

\begin{rem} 
\textup{
Note that we carry out a complete relative bifurcation
analysis without factoring out the underlying group invariance
as in \cite{TZ}.
}
\end{rem}
 \bigskip
 
\subsubsection{Manifold of equilibria}\label{manifold}

The ideas of the previous subsection extend readily to
the more general case that $\{0\}\times \kernel (\Id-S)$ is 
the tangent space at zero of an $\ell$-dimensional
manifold of equilibria, not necessarily generated by group invariance.

Let $\| (a,b) \| := |a| + \| b \|_{X_1}$ as before. 

\begin{ass} \label{man} 
There exists an $\ell$-dimensional family indexed
by $\delta \in \RR^\ell$ of locally invertible Lipschitz
coordinate changes
$\Phi^{\eps,\delta}: \RR^n_a\times X_1 \to \RR^n_a\times X_1$
taking \eqref{diff} to a system of the same form 
\ba
\label{transdiff}
\Delta a&= (R(\e, \delta, a,b)- \Id) a+ N_1(\e, \delta ,a,b),\\
\Delta b&=  (S(\e, \delta ,a,b)- \Id)b+ N_2(\e, \delta, a,b),\\
\ea
$N_j(\e,\delta,0,0) \equiv 0$, such that:
 \begin{itemize}
 \item[{\rm (i)}] 
The maps $\Psi^\eps(\cdot)=(\Psi_a, \Psi_b)(\cdot):= (\Phi^{\eps,\cdot})^{-1}(0,0)$ 
and $\hat \Psi^\eps(\cdot):= (\Phi^{\eps,\cdot})(0,0)$ 
are each differentiable at $\delta=0$, uniformly in $\e$ in the 
sense that
$$
 \sup_{|\e| \leq \e_0} \frac{1}{|x|} \| \Psi^\e(x) - \partial_\delta
 \Psi^\e(0) \cdot x) \| \longrightarrow 0, \quad \mbox{as $x \to 0$.}
$$
 \item[{\rm (ii)}] 
 $\kernel(\Id-S(\e,a,b))=\Span\{\partial_\delta \Psi_b(0)\}$,
while $\partial \hat \Psi(0)$ is full rank.
  \item[{\rm (iii)}] For $\delta$ in a neighborhood of $0$, 
the map $\Phi^{\e,\delta}(\cdot,\cdot)$ is Lipschitz in $(a,b),$ for small $\| (a,b) \|.$
  \item[{\rm (iv)}] For $|\eps|$,  $\|(a,b\|$ sufficiently small,
the map $\Phi^{\e, \cdot}(a,b)$ is continuous in $\delta$ for $|\delta|$ 
sufficiently small.
  \end{itemize}
   \end{ass}

%

\br\label{manrmk}
\textup{
Under Assumption \ref{man}(i)--(ii), $\{0\}\times \kernel(\Id -S)$ 
is the tangent space at $\delta=0$
of the manifold of equilibria $\{\Psi(\delta): \, \delta\in \RR^\ell\}$.
}
\er

 \begin{prop} \label{eq-man} For system \eqref{transdiff},
suppose that $N_2$ satisfies {\rm \eqref{ass-N2}(i)}, and that Assumptions {\rm \ref{B'''}} and {\rm \ref{man}} are satisfied, uniformly in $\delta$. 
Then, for $|\e| \leq \e_0$ and any $C>0$,  
for each sufficiently small
equilibrium $(a,b)$ of \eqref{diff}, there exists $\delta \in B(0,r)$ 
such that $(\hat a, \hat b) := \Phi^{\e,\delta}(a,b)$ satisfies 
$\| \hat b \|_{X_1} \leq C \| \hat a \|.$ 
 \end{prop}
 
 \begin{proof} 
We first observe, given $\|(a_0,b_0)\|$ sufficiently small,
that there exists $\delta_0$ such that 
$ \| \Phi^{\e,\delta_0}(a,b) \|$ minimizes $ \| \Phi^{\e,\delta}(a,b) \|$
for $\delta $ lying in a neighborhood $\CalN$ about zero.
For, by (i)--(ii), $\|\hat \Phi^{\e,\delta}(0,0)\|$ is bounded from
below by some $c_0>0$ for $\delta$ outside a smaller neighborhood $\CalN'$ 
of zero, whereas, by (iii),
$\Phi^{\e,\delta}(a,b)-\Phi^{\e,\delta}(0,0)\le C\|(a,b)\|$, without loss of generality
$C>1$.
Thus, for $\|(a,b)\|\le c_0/2C$, $\|\Phi^{\e,\delta}(a,b)\|\ge c_0/2> \|(a,b)\|$
for $\delta$ outside $\CalN'$, 
and so, by continuity, (iv), there exists a minimal
$\delta_0\in \CalN'$.

Let $(a_0,b_0)$ be a sufficiently small equilibrium,
and let $(a,b)$ be such a minimizing equibrium, with
$(a,b)=\Phi^{\e,\delta}(a_0,b_0).$ 
Assume by way of contradiction that $|a|\le c_0 \|b\|_{X_1}$
for $c_0>0$ sufficiently small, or $a=o(b)$.
{}From $g(\eps, a, b)=0$, by Lemma \ref{g+}, we have 
$$
b-(\Id -S(\eps,a,b))^{-1}N_2(\eps, a,b)=\omega \in \kernel (\Id -S(\eps,a,b)) \cap X_1, $$
where, by \eqref{ass-N2}(i), 
$$
\|(\Id -S(\eps,a,b))^{-1}N_2(\eps, a,b)\|_{X_1} \le C(|a|+\|b\|_{X_1})^2
=o(b),
$$
hence $(a,b)=(0,\o)+ o(\o)$.

By (i)--(ii), associated with $\o$ is $\delta \in B(0,r)$ such that 
$\partial_\delta \Psi(0) \cdot \delta = (0,\omega),$ hence
$$
\Psi(\delta)= (0,\o)+o(\o)= (a,b)+o(\o).
$$ 
Applying $\Phi^{\e,\delta}$ to both sides, and recalling the definition
of $\Psi$ in (ii), we obtain
$$ \begin{aligned} 
(0,0)=\Phi^{\eps,\delta}\Psi(\delta)&= \Phi^{\eps,\delta}\big((a,b)+o(\o)\big)
= \Phi^{\eps,\delta}(a,b)+o(\o),\\
\end{aligned} $$ 
or
$$
\Phi^{\eps,\delta}(a,b)=o(\o),
$$
contradicting the assumed minimality of $(a,b).$ 
 \end{proof}

\begin{cor}\label{man-bif}
Under the assumptions of Proposition {\rm \ref{lipenough}},
suppose further that Assumption {\rm \ref{man}} is satisfied. 
Then, for $(\eps, a, b)$ sufficiently small,
the $\ell$-dimensional 
family of solutions $(\e_0(a,\delta), a,  B(\e_0(a,\delta), a, 0))$ 
determined by $\o \equiv 0$ in Proposition {\rm \ref{lipenough}} applied
to systems \eqref{transdiff} are the unique solutions 
of $(f,g)=0$, apart from the 
trivial solutions $(a,b)\in \{\Psi(\delta)\}$.
\end{cor}

\begin{proof} Identical with that of Corollary \ref{group-inv}.
\end{proof}

\begin{rem}\label{easyfamily}
\textup{
In the context of Hopf bifurcation, 
existence of such a family $\Phi^{\eps,\delta}$ 
amounts to the assumption that $\kernel(\Id-S)=\kernel (\tilde L)$ 
corresponds to (i.e., is the tangent manifold of) a Lipschitz manifold of
equilibria of the original differential system, with $\Phi^{\eps, \delta}$
affine maps consisting of translation in $\tilde u:=\bar u+ a\phi + b$
by points $\bar u^{\eps, \delta}$ on the manifold,
composed with linear coordinate changes to and from eigenvariables $(a,b)$.
}
\end{rem}



\subsection{Brouwer-based Implicit Function Theorem }\label{Brouwer}

 We show in this section that under weaker assumptions, we can still obtain a bifurcation result. For Assumption \ref{B'''bis}, we substitute the following, weaker assumption. 
 
  \begin{ass} \label{B'} The right inverse $(\Id - S)^{-1}$ given by Assumption {\rm \ref{B'''}} is continuous with respect to $\e,a,b \in {\cal V},$ in ${\cal L}(X_2, {\cal B}_1)$ norm, 
and satisfies estimates \eqref{cont-0} from Assumption {\rm \ref{B'''bis}}. 
It is moreover Lipschitz with respect to $b,$ as a function from ${\cal B}_1$ to ${\cal L}({\cal B}_2, {\cal B}_1),$ in the sense that
\be\label{blipdisplay}
   \sup_{\e,a,b \in {\cal V}}
\|\partial_{b} (\Id - S)^{-1}\|_{
{\cal L}({\cal B}_1, 
{\cal L}({\cal B}_2, {\cal B}_1))}
< \infty.
\ee
 \end{ass} 
  
 \begin{lem} 
If $S$ is continuous in $(\e,a)$ for $(\e,a,b) \in {\cal V}$ with
respect to the ${\cal L}(X_2, {\cal B}_1)$ norm
and uniformly bounded  in $\CalL(X_2)$, and 
 the sequence $\sum_{j= 0}^N S(\e,a,b)^j$ converges uniformly with respect to $\e,a,b 
 \in {\cal V}$
 in ${\cal L}(X_2,{\cal B}_1)$ norm, then $(\Id - S)^{-1}$ is continuous with respect to $\e,a,b \in {\cal V},$ in ${\cal L}(X_2, {\cal B}_1)$ norm.  
 \end{lem}
 
\begin{proof}
Since we have assumed
Lipschitz continuity in $b$ in the stronger norm $\CalL({\cal B}_1, \CalL({\cal B}_2, {\cal B}_1)$,
it is sufficient to show continuity in $(\e, a)$
for $(\e,a,b)\in \CalV$ with respect to $\CalL(X_2, {\cal B}_1)$,
which reduces, by continuity of uniform limits of continuous functions,
to establishing continuity of $S^j$.
Suppressing $b$, denote $S'=S(\e',a')$, $S=S(\e',a')$.
Then,
$$ \begin{aligned}
\|(S')^j-S^j\|_{\CalL(X_2,B_1)}&\le
\|\big((S')^{j-1}-S^{j-1}\big)S'\|_{\CalL(X_2,{\cal B}_1)}\\
&\qquad
+
\|S^{j-1}(S'-S)\|_{\CalL(X_2,{\cal B}_1)}\\
&\le
\|(S')^{j-1}-S^{j-1}\|_{\CalL(X_2,{\cal B}_1)}
\|S'\|_{\CalL(X_2)}\\
&\qquad
+
\|S^{j-1}\|_{\CalL({\cal B}_1)}
\|S'-S\|_{\CalL(X_2,{\cal B}_1)},\\
\end{aligned} $$ 
from which we obtain the result by induction on $j$.
\end{proof}
 
 The main difference between Assumption \ref{B'''bis} and Assumption \ref{B'} 
is that in the latter, we do \emph{not} assume Lipschitz regularity 
with respect to $\e, a,$
an advantage since $(\Id -S)^{-1}$ typically 
depends weakly or not at all on $b$ (see Remark \ref{notheta} and
Section \ref{revisit}).
In place of the Contraction-mapping principle of Proposition \ref{lipenough}, 
which used Lipschitz regularity in $\e, a,$ we will use the following 
Implicit Function Theorem.

 \begin{rem} \textup{
 In the following fixed point argument (Proposition \ref{prop-B}), it would be sufficient to assume that for fixed $x \in X_2,$ $(\Id - S)^{-1} x$
is continuous with respect to $\e,a,b,$ in ${\cal B}_1$ norm,
which would be implied by uniform convergence of
$\sum_{j= 0}^N S(\e,a,b)^jx$ in $\CalB_1$ norm for each fixed $x\in X_2$.
 }
 \end{rem}

 \begin{lem}[Brouwer-based IFT]\label{BIFT}
Let $F: \R^m \times \R^n \to \R^m$ be continuous in a neighborhood of the origin, 
$F(0,0)=0$, with $F(\delta,0)$ differentiable at $\delta=0,$ and
$\partial_\delta F(0,0)$ invertible.  Then, for $|a|$ sufficiently small,
there exists a solution $\delta=\delta(a)$ of $F(\delta(a),a)=0$, $\delta(0)=0$,
with $\delta(\cdot)$ continuous at $a=0.$
\end{lem}

\begin{proof}
Expand
$$
\d_\delta F(0,0)^{-1} F(\delta,a) = \delta + n_1(\delta) + n_2(\delta,a),$$
where 
 $$ n_1(\delta) := \d_\delta F(0,0)^{-1} \big(F(\delta,0)- F(0,0) - \delta \d_\delta F(0,0)\big).$$
 and
 $$ n_2(\delta,a) := \d_\delta F(0,0)^{-1} \big( F(\delta,a) - F(\delta,0) \big).$$ 
Both $n_1$ and $n_2$ are continuous in a neighborhood of the origin. By differentiability, $n_1=o(|\delta|)$, while by uniform continuity in a compact neighborhood of the origin,
$n_2\to 0$ as $|a|\to 0$. Thus, for $|a|$ sufficiently small,
$-n_1 - n_2$ takes a small ball centered at the origin in $\R^m$ to itself,
of radius 
$$
2\sup_{|\delta|\le c_0}|n_2(\delta, a)|=o(|a|),
$$
whence there exists a solution 
$\delta= (-n_1 + n_2)(\delta, a)$ by the Brouwer fixed-point Theorem. 
By construction, $\delta(a)=o(a)$, yielding 
$\delta(0) = 0$ and continuity of $\delta(\cdot)$ at $a=0$.
\end{proof}

Notice that differentiability in $\delta$ for $F$ is assumed at the origin only, and that no regularity is asserted for $\delta(\cdot)$ except at the origin. Assuming further regularity in $a,$ for instance H\"older continuity, 
we would obtain H\"older continuity at $a=0$ of $\delta(\cdot)$
by the same proof, 
incorporating the new bound $|n_2(\delta, a)|\le C|a|^\alpha$
coming from H\"older continuity to work on ball $B(0, 2C|a|^\alpha)$.

\subsubsection{Application to bifurcation}\label{bifapp}

 Lemma \ref{BIFT} is precisely what is needed to
apply the Lyapunov--Schmidt method in finite-parameter bifurcation analyses. Consider a continuous curve $\bar \o,$ as in \eqref{omega}, taking for example
$\bar \omega(\eps, a)\equiv 0$. 

 \begin{prop}\label{prop-B} Suppose that $N_1, N_2$ satisfy {\rm \eqref{ass-N1}, \eqref{ass-N2},} that $R$ satisfies {\rm \eqref{ass-R}, \eqref{ass-R-0}}, and that Assumptions {\rm \ref{B'''}} and {\rm \ref{B'}} hold. Suppose 
further that $a\in \RR.$ 
Then, there exists a map $a \mapsto \e_{\bar \o}(a),$ $\e_{\bar \o}(0) = 0,$ $\e_{\bar \o}$ continuous at $0,$ defined on a neighborhood of $0,$ and such that 
 $$ (f,g)(\e_{\bar \o}(a), a, B(\e_{\bar \o}(a), a, \bar \o(a))) \equiv 0,$$
 in a neighborhood of the origin in $\R_a,$ where $B$ is the continuous map given by Proposition {\rm \ref{red'}}. 
\end{prop}  

\begin{proof} The proof of Proposition \ref{red'}, with the exception of its last paragraph, still holds. 
Substituting continuity for Lipschitz continuity in the argument of
that final paragraph, we obtain a \emph{continuous} (in ${\cal B}_1$-norm) 
map $(\e,a) \mapsto \bar B:=B(\e,a,\bar \o(a))$ 
satisfying \eqref{Bybd}, and such that $g(\e,a,\bar B) \equiv 0,$ 
for $|\e| + |a|$ small. 
Then, following the notations of the proof of Proposition \ref{lipenough}, consider the continuous map $(\e,a) \mapsto \bar f(\e,a),$ $\bar f(0,0) = 0.$ By assumption on $R,$ $\bar f(\e,0)$ is differentiable at $\e = 0,$ and $\d_\e \bar f(0,0) = 1.$ We can apply Lemma \ref{BIFT}, and a solution to $\bar f= 0$ yields a solution to $f = 0.$
 \end{proof} 



\section{Linearized bounds}\label{linest}

In the remainder of the paper,
we apply the framework of Section \ref{framework} to
the problem of bifurcation of viscous shock waves, first in the
one-dimensional and then in the multi-dimensional case.
We begin by assembling the relevant linearized estimates on
the one-dimensional flow near a viscous shock.

\subsection{Hypotheses and notation}\label{notation}
As described in Section \ref{model2}, consider a one-parameter 
family of standing viscous shock solutions 
\be\label{prof}
u(x,t)=\bar u^\eps(x),
\qquad \lim_{z\to \pm \infty} \bar u^\eps(z)=u_\pm^\eps
\quad \hbox{\rm (constant for fixed $\eps$)},
\ee 
of a smoothly-varying family of conservation laws 
\begin{equation}
\label{sys}
u_t =\CalF(\e, u):= u_{xx}- F(\e, u)_x,
\qquad u\in \RR^n,
\end{equation}
with associated linearized operators
\be\label{Ldef}
L(\e) :=\partial \CalF/\partial u|_{u=\bar u^\e}
= -\partial_x A^\eps(x) + \partial_x^2,
\ee
$A^\eps(x):= F_u(\e,\bar u^\eps(x))$, denoting
$A^\eps_\pm:=\lim_{z\to \pm \infty} A^\eps(z)=F_u(\e,u^\eps_\pm)$.
Profiles $\bar u^\eps$ satisfy the standing-wave ODE
\be\label{ode}
u'=F(\eps, u)- F(\eps, u^\eps_-).
\ee

We recall the standard assumptions of \cite{ZH, Z4, Z1}:
\medbreak

\quad (H0) \quad  $F\in C^{k}$, $k\ge 2$.
\medbreak
\quad (H1)  \quad $\sigma (A^\eps_\pm)$ real, distinct, and  nonzero.
\medbreak
\quad (H2)  \quad  Considered as connecting orbits of \eqref{ode}, 
$\bar u^\eps$ are transverse and unique up to translation,
with dimensions of the stable subpace $S(A^\eps_+)$ 
and the unstable subspace $U(A^\eps_-)$ 
summing for each $\eps$ to $n+1$.
\medbreak

\br\label{laxrmk}
\textup{
Condition (H1) implies that $u^\eps_\pm$ are nonhyperbolic rest points of
ODE \eqref{ode}, $\Re \sigma(F_u(\eps, u^\eps_\pm))\ne 0$, 
whence, by standard ODE theory,
\be\label{profdecay}
|\partial_x^\ell (\bar u^\eps-u^\eps_\pm)(x)|\le Ce^{-\eta|x|},
\qquad
0\le \ell\le k+1,
\ee
for $x\gtrless 0$, some $\eta$, $C>0$; 
in particular, $|(\bar u^\eps)'(x)| \le Ce^{-\eta|x|}$.
Condition (H2) implies in part that
$\bar u^\eps$ is of standard {\it Lax type}, 
i.e., the hyperbolic convection matrices 
$A^\eps_\pm$ 
at $\pm\infty$ have, respectively, $p-1$ negative and
$n-p$ positive eigenvalues for $1\le p \le n$, where $p$ is the characteristic
family associated with the shock. 
For further discussion, see \cite{ZH, Z1, HZ}.
}
\er

To (H0)--(H2) we adjoin the generalized spectral (i.e., Evans function) 
condition:

\medbreak

\quad (\D) \quad
On a neighborhood of $\{\Re \lambda \ge 0\}\setminus \{0\}$,
the only zeroes of $D$ are (i) a zero of multiplicity one at
$\lambda=0$, and (ii) a crossing conjugate pair of zeroes
$\lambda_\pm(\eps)=\gamma(\eps)+i\tau(\eps)$ with
$\gamma(0)=0$, $\partial_\eps \gamma(0)>0$, and $\tau(0)\ne 0$.
\medbreak

\begin{lem}\label{HP}
Conditions {\rm (H0)--(H2)} and {\rm (\D)} are equivalent to conditions 
{\rm (P)(i)--(iii)} of the introduction (stated in Section {\rm \ref{formule}}) together with 
$F\in C^k$, $k\ge 2,$
simplicity and nonvanishing 
of $\sigma(A^\eps_\pm)$, and the Lax condition
\be\label{lax}
\dim S(A^\eps_+)+ \dim U(A^\eps_-)=n+1,
\ee
with $\bar u^\eps$ (linearly and nonlinearly) stable
for $\eps <0$ and unstable for $\eps\ge 0$.
\end{lem}

\begin{proof}
By (H0), (H2) and standard ODE theory, solution
$\bar u^\eps$ of \eqref{ode} is $C^3$ in $x$,
$C^2$ in $(x,\eps)$, and decays at exponential rate
in first two derivatives to endstates $u^\eps_\pm$
as $x\to \pm \infty$.
Thus, $L(\e)$ is asymptotically constant-coefficient,
and its essential spectrum may be computed as in 
\eqref{Lplusminus}--\eqref{disp}, Section \ref{spectral}, 
to lie entirely within a parabolic region:
$$
\sigma_{\ess}(L(\e))\subset
\{\Re \lambda \le -\theta |\Im \lambda|^2\},
\qquad \theta>0,
$$
verifying (P)(i).  On the other hand, computation \eqref{disp}
shows that (P)(i) implies (H2) in the generic situation that
$\sigma(A^\eps_\pm)$ is simple and nonvanishing.
Likewise, recalling \cite{AGJ, GZ} that to the right of the essential
spectrum boundary, zeroes of $D_\eps$ correspond to eigenvalues
of $L(\e)$, we find that
(\D)(i) and (ii) are exactly (P)(ii) and (iii).
Finally, recalling the {\it stability criterion} 
(see, e.g., \cite{ZH, Z4, MaZ3, MaZ4}) that, under (H0)--(H2),
linear and nonlinear stability of $\bar u^\eps$ is equivalent to
\medbreak
\quad (${\mathcal D}$) \quad
On $\{\Re \lambda \ge 0\}$, $D_\eps$ has a single
zero at $\lambda=0$,
\medbreak
\noindent
we find that $\bar u^\eps$ is stable precisely for $\eps<0$.
\end{proof}

\begin{rem}\label{ucrmk}
\textup{
The undercompressive shocks studied in \cite{HZ}
give an example of shocks satisfying (P) but not \eqref{lax}.
Their weaker decay properties are insufficient for
the convergence argument of Section \ref{refined};
see Remark \ref{linuc}.
}
\end{rem}

\br\label{simple}
\textup{
Under (H0)--(H2),
condition \D(i) is equivalent to 
$\langle \ell^\eps, (\bar u^\eps)'\rangle=\ell^\eps\cdot (u^\eps_+- u^\eps_-)
\ne 0$ for
$\ell^\eps$ (constant) orthogonal to $S(A^\eps_-)\cup U(A^\eps_+)$ \cite{GZ}.
Under the normalization $\langle \ell^\eps, (\bar u^\eps)'\rangle=1$,
operator $\Pi_0^\eps f:= (\bar u^\eps)' \langle \ell^\eps, f\rangle$
plays the role of a ``generalized spectral projection'' onto 
$\kernel L(\eps)=\Span \{(\bar u^\eps)'\}$, and $\ell^\eps$ the
role of a generalized left eigenfunction \cite{ZH}.
Note that $\ell^\eps$ lies outside the domain
of $\Pi_0^\eps$, a consequence of the absence of a spectral gap.
}
\er

Finally, we introduce the Banach spaces $\CalB_1=L^2$,
$\CalB_2=\partial_x L^1\cap L^2$, 
$$
X_1=\{f:\, |f(x)|\le C(1+|x|)^{-1}\},
$$
and
$$
X_2=\partial_x \{f:\, |f(x)|\le C(1+|x|)^{-2}\} \cap X_1,
$$
equipped with norms
$\|f\|_{\CalB_1}=\|f\|_{L^2}$,
$\|\partial_x f\|_{\CalB_2}=\|f\|_{L^1}+\|\partial_x f\|_{L^2}$,
$$\|f\|_{X_1}=\|(1+|x|)f\|_{L^\infty},
\quad \hbox{\rm and } \quad
\|\partial_x f\|_{X_2}=\|(1+|x|)^2f\|_{L^\infty}+\|\partial_x f\|_{X_1},
$$
where $\partial_x$ is taken in the sense of distributions.
By inspection, we have that 
$\CalB_2\subset \CalB_1$,
$X_2\subset X_1$,
$X_1\subset \CalB_1$,
$X_2\subset \CalB_2$,
and the closed unit ball in $X_1$ is closed as
a subset of $\CalB_1$.

\subsection{Projector bounds}\label{projector}
    
\begin{lem}\label{efns}
Under assumptions {\rm (H0)--(H2)}, {\rm (\D)},
%
associated with eigenvalues $\lambda_\pm(\eps)$ of $L(\e)$
are right and left eigenfunctions $\phi^\eps_\pm$ 
and $\tilde \phi^\eps_\pm \in C^k(x,\eps)$, $k\ge 2$ as in Assumption {\rm (H0)} from Section {\rm \ref{notation}},
exponentially decaying in up to $k$ derivatives as
$x\to \pm \infty$,  and $L(\e)$-invariant projection
\be\label{proj}
\Pi^\e f:= \sum_{j=\pm} \phi^\eps_j(x)\langle \tilde \phi^\eps_j, f\rangle
\quad
\ee
onto the total (oscillatory) eigenspace $\Sigma^\eps:=\Span \{\phi^\eps_\pm\}$,
bounded from $L^q$ or $\CalB_2$ to $W^{2,p}\cap X_2$ 
for any $1\le q,p\le \infty$.  Moreover,
\be\label{zeromassefn}
\phi^\eps_\pm= \partial_x \Phi^\eps_\pm,
\ee
with $\Phi^\eps\in C^{k+1}$ exponentially
decaying in up to $k+1$ derivatives as $x\to \pm \infty$.
\end{lem}

\begin{proof}
{}From simplicity of $\lambda_\pm$, we obtain either by standard 
spectral perturbation theory \cite{Kat} or by direct Evans-function
calculations \cite{GJ1, GJ2, ZH} that there exist $\lambda_\pm(\cdot) $,
$\phi^\eps_\pm(\cdot)\in L^2$ with the same smoothness $C^2(\eps)$ 
assumed on $F$.  The exponential decay properties in $x$ then follow
by standard asymptotic ODE theory; see, e.g., \cite{GZ, Z1}.
Finally, recall the observation of \cite{ZH} that, by divergence form
of $L(\e)$, we may integrate  $L(\e) \phi= \lambda \phi$ from
$x=-\infty$ to $x=+\infty$ to obtain 
$\lambda \int_{-\infty}^{+\infty}\phi(x)dx=0$, and thereby
(since $\lambda_\pm \ne 0$ by assumption)
\be\label{zmass}
\int_{-\infty}^{+\infty}\phi_\pm(x)dx=0,
\ee
from which we obtain by integration \eqref{zeromassefn} with the stated
properties of $\Phi_\pm$.
{}From \eqref{zeromassefn} and representation \eqref{proj}, we 
obtain by H\"older's inequality the stated bounds on projection $\Pi^\e$.
\end{proof}

Defining $\tilde \Pi^\eps:=\Id -\Pi^\eps$, $\tilde \Sigma^\eps:= \Range \tilde \Pi^\eps$,
and $\tilde L(\e):= L(\e)\tilde \Pi^\eps$, denote by
\be\label{kernel}
G(x,t;y):= e^{L(\e) t}\delta_y(x)
\ee
the Green kernel associated with the linearized
solution operator $e^{Lt}$ of the linearized evolution equations
$u_t=L(\e) u$, and
\be\label{transkernel}
\tilde G(x,t;y):= e^{\tilde L(\e) t}\tilde \Pi^\e \delta_y(x)
\ee
the Green kernel associated with the transverse linearized
solution operator $e^{\tilde L(\e) t}\tilde \Pi^\e$.
By direct computation,
$G = \CalO +  \tilde G$, where
\be \label{O}
\CalO(x,t;y):= e^{(\gamma(\eps)+i\tau(\eps)) t}\phi_+(x) \tilde \phi_+^t(y)
+
e^{(\gamma(\eps)-i\tau(\eps)) t}\phi_-(x) \tilde \phi_-^t(y).
\ee

\subsection{Short time estimates}\label{shorttime}

\begin{lem}\label{masspres}
Under assumptions {\rm (H0)--(H2)}, {\rm (\D)},
%
for 
$0 \leq t\le T$,
 any fixed $T>0$, $|\e| \leq \e_0,$ there exists $C >0,$ depending only on $\e_0$ and $T,$ such that
\begin{eqnarray}
 \| e^{\tilde L(\e) t} \tilde \Pi^\e f \|_{{\cal B}_1} & \leq & C \| f \|_{{\cal B}_1}, \label{BDuhamel1} \\
 \| e^{\tilde L(\e) t}\tilde \Pi^\e \partial_x f \|_{\CalB_2} & \leq & C( \| f \|_{L^1} + t^{-1/2} \|f\|_{{\cal B}_1}), \label{BDuhamel2}
 \end{eqnarray}
 and
 \begin{eqnarray}
\| e^{\tilde L(\e) t} \tilde \Pi^\e f \|_{X_1} & \leq & C \| f \|_{X_1}, \label{XDuhamel1} \\
 \| e^{\tilde L(\e) t} \tilde \Pi^\e \d_x f \|_{X_2} & \leq & C t^{-1/2} \sup_x (1 + |x|^2) |f (x)|. \label{XDuhamel2} 
 \end{eqnarray}
 Besides,
   \begin{eqnarray}
  \| \d_\e (e^{\tilde L(\e) t} \tilde \Pi^\e) f \|_{{\cal B}_1} & \leq &  C \| f \|_{{\cal B}_1}, \label{epsDuhamel1} \\
  \| \d_\e (e^{\tilde L(\e) t} \tilde \Pi^\e) \d_x f \|_{{\cal B}_2} & \leq &  C ( \| f \|_{L^1} + \| f \|_{{\cal B}_1}). \label{epsDuhamel2}
\end{eqnarray}
\end{lem}

\begin{proof} 
{}From standard parabolic semigroup bounds 
$$ 
|e^{L(\e) t}|_{L^2\to L^2} \le Ct,\qquad
|e^{L(\e) t}\partial_x|_{L^2\to L^2} \le Ct^{-1/2}, 
$$
and properties 
\begin{equation} \label{1-12.0} e^{L(\e) t}=e^{L(\e) t}\Pi^\e + e^{\tilde L(\e) t}\tilde \Pi^\e,
\end{equation}
and $\|\Pi^\e \partial_x f\|_{L^2}\le |f|_{L^2}$,
we obtain
\ba\label{oneB}
\|e^{L(\e) t}\partial_x f \|_{L^2},\,
\|e^{\tilde L(\e) t}\tilde \Pi^\e \partial_x f \|_{L^2}
&\le Ct^{-1/2} \|f\|_{L^2},\\
\ea
and, in particular, \eqref{BDuhamel1}. 
Likewise, we may obtain integrated bounds
\ba\label{twoB}
\Big\|\int e^{L(\e) t}\partial_x f \Big\|_{L^1},\,
\Big\|\int e^{\tilde L(\e) t}\tilde \Pi^\e \partial_x f \Big\|_{L^1}
&\le C \|f\|_{L^1}\\
\ea
using the divergence form of $L(\e) $,
by integrating the linearized equations with respect to $x$
to obtain linearized equations $U_t=\CalL(\e)  U$ for integrated variable
$$
U(\e,x,t):=\int_{-\infty}^x u(\e,z,t)dz,
\qquad
u(\cdot,t):=e^{L(\e) t}\partial_x f,
$$
with linearized operator $\CalL(\e) :=-A^\eps(x)\partial_x + \partial_x^2$
of the same parabolic form as $L(\e) $, then applying 
standard parabolic semigroup estimates
(alternatively, pointwise Green function bounds as in Proposition
\ref{greenbounds}) 
to bound 
$$
\|U(\cdot, t, \eps)\|_{L^1}=
\|e^{\CalL(\e) t}f\|_{L^1}\le C\|f\|_{L^1},
$$
the $e^{\tilde L t}\tilde \Pi^\e$ bound then following by
relation $e^{L(\e) t}=e^{L(\e) t}\Pi^\e + e^{\tilde L(\e) t}\tilde \Pi^\e$
together with
$\|\int \Pi^\e \partial_x f\|_{L^1}\le |f|_{L^1}$.
Combining \eqref{oneB} and \eqref{twoB}, we obtain \eqref{BDuhamel2}.
Bounds \eqref{XDuhamel1} and \eqref{XDuhamel2}
follow by identical arguments, substituting
for $L^2$ and $L^1$ the weighted spaces $\|(1+|x|) \cdot \|_{L^\infty}$
and $\|(1+|x|)^2 \cdot \|_{L^\infty}$.

Bound \eqref{epsDuhamel1} follows by the observation that
$v:=\partial_\eps e^{L(\eps)t}f$ solves variational equation
$
v_t - Lv= (\partial_\eps L) e^{L t}f$
with initial data $v_{t=0}=0$, hence, by Duhamel's principle,
$$
v(t)= \int_0^t e^{L(t-s)} (\partial_\eps L) e^{Ls}f \, ds.
$$
By \eqref{BDuhamel1}, therefore, noting that 
$\|(\partial_\eps L)g\|_{\CalB_1}\le
C(\|g\|_{\CalB_1} +\|\partial_x g\|_{\CalB_1})$,
\begin{equation} \label{11.01}
\|v(t)\|_{\CalB_1}\le
C\int_0^t (1+s^{-1/2})\|f\|_{\CalB_1} \, ds \le C_2 \|f\|_{\CalB_1},
\end{equation}
and \eqref{epsDuhamel1} follows from \eqref{11.01}, the $\e$-regularity of $\Pi^\e$ stated in Lemma \ref{efns}, \eqref{BDuhamel1} and \eqref{1-12.0}.

Finally, let
 $$
W(\e,x,t):=\int_{-\infty}^x w(\e,z,t)dz,
\qquad
w(\cdot,t):=(\d_\e e^{L(\e) t})\partial_x f.
$$
Using again the divergence form of $L,$ we find that $W$ satisfies a parabolic equation
 $$ \d_t W - {\cal L} W = \int_{-\infty}^x (\d_\e L) e^{t L} \d_x f \, dx = \d_{\e} ( \d_u F(\e, \bar u^\e)) e^{t L} \d_x f,$$
 where ${\cal L}$ is of the same form as $L;$ this implies
 $$ \| W (t) \|_{L^1} \leq \int_0^t \| \d_{\e} ( \d_u F(\e, \bar u^\e)) \|_{L^\infty} \| e^{s L} \d_x f \|_{L^1} \,ds,$$
 from which we obtain 
 \begin{equation} \label{1-12} \| W (t) \|_{L^1} \leq (C \int_0^t s^{-1/2} \,ds) \| f \|_{L^1},\end{equation} 
Here, we have used the standard bound
$ \|e^{t L} \d_x f \|_{L^p} \leq C t^{-1/2} \| f \|_{L^p},$
$ 1 \leq p \leq \infty,$
 following from
pointwise (high-frequency) resolvent estimates
obtained by asymptotic ODE theory \cite{Sat, PW, ZH}.
To establish \eqref{epsDuhamel2}, it remains to bound $\| w \|_{L^2}.$ This is done as in the proof of \eqref{epsDuhamel1}:
\ba \label{1-12.2}
\|w(t)\|_{L^p}&= \Big\|\int_0^t e^{L(t-s)} (\partial_\eps L) 
e^{Ls}\partial_x f \, ds\Big\|_{L^p} \\ 
& = \Big\| \int_0^t e^{L(t-s)} \d_x \Big( \d_{\e} ( \d_u F(\e, \bar u^\e)) e^{Ls}\partial_x f \Big) \, ds
\Big\|_{L^p}\\
&\le C \int_0^t (t-s)^{-1/2} \| e^{Ls} \partial_x f\|_{L^p} \, ds\\
&\le C \Big( \int_0^t (t-s)^{-1/2}  s^{-1/2} \, ds \Big) \|f\|_{L^p}\\
\ea
Applying to $p=2$ gives the result. Estimate \eqref{epsDuhamel2} follows from \eqref{1-12}, \eqref{1-12.2} and \eqref{1-12.0}.
\end{proof}

\subsection{Pointwise Green function bounds}\label{ptwise}

Supressing the parameter $\eps$,
denote by $a_j^\pm$, $l_j^\pm$, and $r_j^\pm$
the eigenvalues and left and right eigenvectors of 
$A^\eps_\pm= F_u(\e,u^\eps_\pm)$.  

In the following Proposition, $\eta$ and $M$ denote positive constants, \be \label{zjk}
z_{jk}^\pm(y,t):=a_j^{\pm}\left(t-\frac{|y|}{|a_k^{-}|}\right)
\ee
and
\be \label{barbeta}
\bar \beta^{\pm}_{jk}(x,t;y):= 
\frac{|x^\pm|}{|a_j^{\pm} t|} 
+
\frac{|y|}{|a_k^{-} t|} 
\left( \frac{a_j^\pm}{a_k^{-}}\right)^2 
\ee
are approximate scattered characteristic
paths and effective diffusion rates along them,
scattering coefficients $[c_{k,-}^{j,i}]$, $i=-,0,+$ are constant,
$x^\pm$ denotes the positive/negative
part of $x$,  
$\chi_{\{t\ge 1\}}$ is a smooth cutoff function in $t$,
identically one for $t\ge 1$ and identically zero for $t\le 1/2$,
and indicator function $\chi_{\{ |a_k^{-}t|\ge |y| \}}$ is 
one for $|a_k^{-}t|\ge |y|$ and zero otherwise.

\begin{prop} [\cite{ZH, Z4, MaZ3}] \label{greenbounds}
Under assumptions {\rm (H0)--(H2)}, {\rm (\D)},
\be\label{ourdecomp}
\tilde G =  \CalE+  \CalS + \CalR,
\ee
$\tilde G$ as in {\rm \eqref{transkernel}}, where, for $y\le 0$ and all $t\ge 0$:
\ba \label{E}
\CalE&(x,t;y):=
\sum_{a_k^- > 0}
[c^{0}_{k,-}]
\bar U'(x)
l_k^{-t}
\left(\hbox{\rm errfn} \left(\frac{y+a_k^{-}t}{\sqrt{4t}}\right)
-\hbox{\rm errfn} \left(\frac{y-a_k^{-}t}{\sqrt{4t}}\right)\right)
\ea
and
\ba \label{S}
&\CalS(x,t;y)\\
&:=
\chi_{\{t\ge 1\}} 
\sum_{a_k^{-}<0}r_k^{-}  {l_k^{-}}^t
(4\pi t)^{-1/2} e^{-(x-y-a_k^{-}t)^2 / 4t} 
\\
&+ 
\chi_{\{t\ge 1\}} 
\sum_{a_k^{-} > 0} r_k^{-}  {l_k^{-}}^t
(4\pi t)^{-1/2} e^{-(x-y-a_k^{-}t)^2 / 4t}
\left(\frac{e^{-x}}{e^x+e^{-x}}\right)\\
&+ 
\chi_{\{t\ge 1\}}
\sum_{a_k^{-} > 0, \,  a_j^{-} < 0} 
[c^{j,-}_{k,-}]r_j^{-}  {l_k^{-}}^t
(4\pi \bar\beta_{jk}^{-} t)^{-1/2} e^{-(x-z_{jk}^{-})^2 / 
4\bar\beta_{jk}^{-} t} 
\left(\frac{e^{ -x}}{ e^x+e^{-x}}\right)\\
&+ 
\chi_{\{t\ge 1\}}
\sum_{a_k^{-} > 0, \,  a_j^{+} > 0} 
[c^{j,+}_{k,-}]r_j^{+}  {l_k^{-}}^t
(4\pi \bar\beta_{jk}^{+} t)^{-1/2} e^{-(x-z_{jk}^{+})^2 / 
4\bar\beta_{jk}^{+} t} 
\left(\frac{e^{ x}}{ e^x+e^{-x}}\right)\\
\ea
denote 
excited and scattering terms and 
\ba \label{Rbounds}
\CalR(x,t;y)&= 
O(e^{-\eta (|x-y|+ t)})\\
&+\sum_{k=1}^n 
O \left( (t+1)^{-1/2} e^{-\eta x^+} 
+e^{-\eta|x|} \right) 
t^{-1/2}e^{-(x-y-a_k^{-} t)^2/Mt} \\
&+
\sum_{a_k^{-} > 0, \, a_j^{-} < 0} 
\chi_{\{ |a_k^{-} t|\ge |y| \}}
O ((t+1)^{-1/2} t^{-1/2})
e^{-(x-z_{jk}^{-})^2/Mt}
e^{-\eta x^+}, \\
&+
\sum_{a_k^{-} > 0, \, a_j^{+}> 0} 
\chi_{\{ |a_k^{-} t|\ge |y| \}}
O ((t+1)^{-1/2} t^{-1/2})
e^{-(x-z_{jk}^{+})^2/Mt}
e^{-\eta x^-}, \\
\ea
\ba \label{Rtbounds}
\CalR_t(x,t;y)
&= 
O(e^{-\eta (|x-y|+ t)})
+\sum_{k=1}^n 
O \left( (t+1)^{-1/2} e^{-\eta x^+} 
+e^{-\eta|x|} \right) \\
&\times
(t+1)^{1/2} t^{-3/2}
e^{-(x-y-a_k^{-} t)^2/Mt} \\
&+
\sum_{a_k^{-} > 0, \, a_j^{-} < 0} 
\chi_{\{ |a_k^{-} t|\ge |y| \}}
O ( t^{-3/2})
e^{-(x-z_{jk}^{-})^2/Mt}
e^{-\eta x^+}, \\
&+
\sum_{a_k^{-} > 0, \, a_j^{+}> 0} 
\chi_{\{ |a_k^{-} t|\ge |y| \}}
O ( t^{-3/2})
e^{-(x-z_{jk}^{+})^2/Mt}
e^{-\eta x^-}, \\
\ea
\ba\label{Rybounds}
\CalR_y(x,t;y)&= 
\partial_t  r(x,t;y)+
O(e^{-\eta (|x-y|+ t)})\\
&+\sum_{k=1}^n 
O \left( (t+1)^{-1/2} e^{-\eta x^+} 
+e^{-\eta|x|} \right) 
t^{-1}
e^{-(x-y-a_k^{-} t)^2/Mt} \\
&+
\sum_{a_k^{-} > 0, \, a_j^{-} < 0} 
\chi_{\{ |a_k^{-} t|\ge |y| \}}
O ((t+1)^{-1/2} t^{-1}) 
e^{-(x-z_{jk}^{-})^2/Mt}
e^{-\eta x^+} \\
&+
\sum_{a_k^{-} > 0, \, a_j^{+} > 0} 
\chi_{\{ |a_k^{-} t|\ge |y| \}}
O ((t+1)^{-1/2} t^{-1}) 
e^{-(x-z_{jk}^{+})^2/Mt}
e^{-\eta x^-}, \\
\ea
\ba\label{hatR}
r(x,t;y)&=
O \big( e^{-\eta|y|}
(t+1)^{-1/2} 
\big) \\
&\qquad \times
\Big( \sum_{a_k^-<0}
e^{-(x-y-a_k^{-} t)^2/Mt} 
+
\sum_{a_k^+>0}
e^{-(x-y-a_k^{+} t)^2/Mt} \Big), \\
\ea
\ba \label{Rytbounds}
\partial_t^k\CalR_{y}(x,t;y)
&= 
O(e^{-\eta (|x-y|+ t)})
+\sum_{k=1}^n 
O \left( (t+1)^{-1/2} e^{-\eta x^+} 
+e^{-\eta|x|} \right) \\
&\times
(t+1)^{1/2} t^{-(k+3)/2}
e^{-(x-y-a_k^{-} t)^2/Mt} \\
&+
\sum_{a_k^{-} > 0, \, a_j^{-} < 0} 
\chi_{\{ |a_k^{-} t|\ge |y| \}}
O (t^{-(k+3)/2})
e^{-(x-z_{jk}^{-})^2/Mt}
e^{-\eta x^+}, \\
&+
\sum_{a_k^{-} > 0, \, a_j^{+}> 0} 
\chi_{\{ |a_k^{-} t|\ge |y| \}}
O (t^{-(k+3)/2})
e^{-(x-z_{jk}^{+})^2/Mt}
e^{-\eta x^-} \\
&+
O \big( e^{-\eta|y|}
(t+1)^{-(k+2)/2} 
\big) \\
&\times
\Big( \sum_{a_k^-<0}
e^{-(x-y-a_k^{-} t)^2/Mt} 
+
\sum_{a_k^+>0}
e^{-(x-y-a_k^{+} t)^2/Mt} \Big), \\
\ea
$k\ge 1$,
a faster decaying residual.

Moreover, for $y\le 0$, $k\ge 1$, some $M>0$,
\ba \label{Gytbounds}
\partial_t^k \tilde G_{y}(x,t;y)
&=\sum_{k=1}^n 
O \left( t^{-1/2} e^{-\eta x^+} 
+e^{-\eta|x|} \right) \\
&\times
(t+1)^{k/2} t^{-(k+1/2)}
e^{-(x-y-a_k^{-} t)^2/Mt} \\
&+
\sum_{a_k^{-} > 0, \, a_j^{-} < 0} 
\chi_{\{ |a_k^{-} t|\ge |y| \}}
O\big((t+1)^{-(k+2)/2} \big)
e^{-(x-z_{jk}^{-})^2/Mt}
e^{-\eta x^+}, \\
&+
\sum_{a_k^{-} > 0, \, a_j^{+}> 0} 
\chi_{\{ |a_k^{-} t|\ge |y| \}}
O\big((t+1)^{-(k+2)/2} \big)
e^{-(x-z_{jk}^{+})^2/Mt}
e^{-\eta x^-} \\
&+
O \big( e^{-\eta|y|}
(t+1)^{-(k+2)/2} 
\big) \\
&\times
\Big( \sum_{a_k^-<0}
e^{-(x-y-a_k^{-} t)^2/Mt} 
+
\sum_{a_k^+>0}
e^{-(x-y-a_k^{+} t)^2/Mt} \Big). \\
\ea

Symmetric bounds hold for $y\ge0$.
\end{prop}


\begin{proof}
Evidently, it is equivalent to establish decomposition 
$$
G=\CalO + \CalE+\CalS+\CalR
$$
of the full Green function.
This problem has been treated in \cite{ZH, MaZ3}, 
starting with Inverse Laplace Transform representation
\be\label{ILT}
G(x,t;y)=e^{Lt}\delta_y(x)= \oint_\Gamma e^{\lambda t}(\lambda-L(\e))^{-1} 
\delta_y(x)d\lambda \, ,
\ee
where 
$$
\Gamma:= \partial \{ \lambda : \Re \lambda\le \eta_1 - \eta_2 |\Im \lambda|\}
$$
is an appropriate sectorial contour, $\eta_1$, $\eta_2>0$;
estimating the resolvent kernel 
$G^\eps_\lambda(x,y):=(\lambda-L(\e))^{-1}\delta_y(x)$
using Taylor expansion in $\lambda$,
asymptotic ODE techniques in $x$, $y$, and judicious decomposition
into various scattering, excited, and residual modes;
then, finally, estimating the contribution of various modes to \eqref{ILT}
by Riemann saddlepoint (Stationary Phase) method, moving contour
$\Gamma$ to a optimal, ``minimax'' positions for each
mode, depending on the values of $(x,y,t)$.

In the present case, we may first move $\Gamma$ to a contour
$\Gamma'$ enclosing (to the left) all spectra of $L(\e)$
except for the crossing pair $\lambda_\pm(\eps)$, to obtain
$$
G(x,t;y)= \oint_{\Gamma'} e^{\lambda t}(\lambda-L(\e))^{-1} d\lambda
+ \sum_{j=\pm} 
\Res_{\lambda_j(\eps)} \big( e^{\lambda t}(\lambda-L(\e))^{-1}
\delta_y(x) \big),
$$
where
$\Res_{\lambda_j(\eps)} \big( e^{\lambda t}(\lambda-L(\e))^{-1}
\delta_y(x) \big) = \CalO(x,t;y)$, then estimate the remaining term
$ \oint_{\Gamma'} e^{\lambda t}(\lambda-L(\e))^{-1} d\lambda$
on minimax contours as just described.
See the proof of Proposition 7.1, \cite{MaZ3}, for a detailed 
discussion of minimax estimates $\CalE+\CalS+\CalR$ and of Proposition 7.7, 
\cite{MaZ3}, for a complementary discussion of residues
incurred at eigenvalues in $\{\Re \lambda\ge 0\}\setminus\{0\}$.

We have repaired in \eqref{Rbounds}, \eqref{Rybounds} two minor omissions
in the bounds of \cite{MaZ3}, pointed out already in \cite{Raoofi, HR, HRZ}.
Specifically, (i) in the first term on the righthand side
of \eqref{Rbounds}, \eqref{Rybounds}, we have replaced the
term $O(e^{-\eta t}e^{-|x-y|^2/Mt)})$
appearing in \cite{MaZ3, Z1}, with a corrected version 
$O(e^{-\eta (|x-y|+ t)})$, and (ii) we have
added a missing term 
\ba\label{missing}
\partial_t r &= 
\partial_t O \big( e^{-\eta|y|}
(t+1)^{-1/2} 
\big) 
\Big( \sum_{a_k^-<0}
e^{-(x-y-a_k^{-} t)^2/Mt} 
+
\sum_{a_k^+>0}
e^{-(x-y-a_k^{+} t)^2/Mt} \Big), \\
\ea
which may alternatively be estimated as in \cite{Raoofi, HR, HRZ} as
\be\label{Raoofiest}
\partial_t r=
O \big( e^{-\eta|y|}
(t+1)^{-1} 
\big) 
\Big( \sum_{a_k^-<0}
e^{-(x-y-a_k^{-} t)^2/Mt} 
+
\sum_{a_k^+>0}
e^{-(x-y-a_k^{+} t)^2/Mt} \Big), 
\ee
to the righthand side of \eqref{Rybounds}.

As discussed in \cite{Raoofi}, the first correction concerns
only bookkeeping, and comes from the fact that approximate Gaussians
appearing in $\CalS$ decay only as $e^{-\eta |x-y|}$ in the far field
and not as does the total solution as $O(e^{-\eta t}e^{-|x-y|^2/Mt)})$.
This is only an artifact of our way of displaying the solution.
The second is more fundamental, and deserves further comment.
Namely, different from the constant-coefficient case, 
the (approximate) Gaussians appearing in the Green function for
the variable coefficient case involve convection, diffusion,
and also directional parameters that depend on $x$ and $y$
through $\bar u^\eps(x)$.
Thus, $x$- and $y$-derivatives do not always bring down additional
additional factors $t^{-1/2}$ of temporal decay as for constant-coefficient
Gaussians, but may instead, through derivatives falling on 
spatially-dependent parameters, 
bring down factors $e^{-\eta |x|}$ or $e^{-\eta |y|}$ 
of spatial decay; see \cite{ZH, Z4} for further discussion.
The neglected term \eqref{Raoofiest}, which is harmless for
the stability analyses of \cite{Raoofi, MaZ2, MaZ4, Z1}, 
may be recognized as a term of this form.

Unfortunately, term \eqref{Raoofiest}, though harmless in the
stability analysis, is unacceptable for our bifurcation analysis
(see Remark \ref{linuc}), and thus must be treated in a different
way.
Reviewing the analysis of \cite{MaZ3}, Proposition 7.7,
we see that term \eqref{Raoofiest} results from a term
in \eqref{ILT} of form
\be\label{format}
\partial_y \oint_{\Gamma_1}e^{\lambda t} \lambda \phi(x)\tilde \psi(y)
\, d\lambda,
\ee
where $\Gamma_1$ is a fixed finite arc of a larger contour $\Gamma$,
$\tilde \psi(y)$, $\partial_y \tilde \psi(y) \sim e^{-\eta|y|}$
is a fast-decaying mode and 
$\phi(x)\sim e^{c_1 \Re \lambda + c_2 \Re (\lambda^2)}$ is a slow-decaying
mode; the factor $\lambda$ reflects the fact that this is a Taylor
remainder term in the low-frequency expansion about $\lambda=0$.
The ``old'' estimate \eqref{Raoofiest} makes use of the fact that
each factor of $\lambda$ introduces (by a simple scaling argument)
an additional factor $(t+1)^{-1/2}$ of temporal decay in the
Riemann saddlepoint/Stationary Phase estimate from which the final
bounds are obtained.

The ``new'' estimate \eqref{missing} follows by a still simpler
argument, making use of the fact that
a factor of $\lambda$ corresponds identically to time-differentiation:
$$
\partial_y \oint_{\Gamma_1}e^{\lambda t} \lambda \phi(x)\tilde \psi(y)
\, d\lambda=
\partial_t \oint_{\Gamma_1}e^{\lambda t}\partial_y \phi(x)\tilde \psi(y)
\, d\lambda,
$$
and estimating
$\oint_{\Gamma_1}e^{\lambda t} \partial_y\phi(x)\tilde \psi(y) \, d\lambda$
as before.

Conversely, a time-derivative on the Green function
can always be traded for a factor of $\lambda$
in the integrand of \eqref{ILT}, yielding faster temporal 
decay by factor $(t+1)^{1/2}t^{-1}\sim t^{-1/2}$.  
We use this observation to generate estimates
\eqref{Rtbounds}, \eqref{Rytbounds}, and \eqref{Gytbounds} 
not stated in \cite{MaZ3, Z1}.
Note that conglomerate high-derivative bound \eqref{Gytbounds} does not
contain the ``bookkeeping'' error term $O(e^{-\eta (|x-y|+ t)})$
deriving from our way of decomposing the solution,
and is in fact much easier to get than the more detailed individual 
bounds, requiring only modulus-based stationary-phase estimates
like those used to bound residual $\CalR$.
\end{proof}

\br\label{xyt}
\textup{
The effect of $t$-derivatives for variable-coefficient equations
in bringing down additional factors of temporal decay
is in contrast to that of $x$- and $y$-derivatives,
which eventually bring down also $e^{-\eta|x|}$ and
$e^{-\eta|y|}$ factors not adding extra decay \cite{ZH, HZ}.
}
\er

\br\label{L2bd}
\textup{
Integrating the pointwise bounds of Proposition \ref{greenbounds}, we obtain
$$
\|(\tilde G-\CalE)_y\|_{L^2(x)}\le Ct^{-3/4},
\qquad
\|(\tilde G-\CalE)_{yt}\|_{L^2(x)}\le Ct^{-5/4}, 
$$
similarly as for the convected heat kernel $K$ in the
analysis of Section \ref{refined}.  
}
\er

\br\label{linuc}
\textup{
Note that $\tilde G$ bounds indeed agree with the model given in Section
\ref{xt}.  Specifically, $\CalE$ terms are exactly superpositions
of terms of form $J$.  Likewise, $\CalS$ terms near their centers are 
well-approximated by moving Gaussians $K$, so satisfy
similar bounds; see \cite{Raoofi}.  
Thus, the model analysis well-represents the primary terms
$\CalE$ and $\CalS$ in the transverse Green function expansion
for Lax-type viscous shocks.}

\textup{ Undercompressive shocks give a different model
$$
G\sim K(x,t; y)p(y)+J(x,t;y)q(y)
$$
with $|p_y|, |q_y|\sim e^{-\eta|y|}$, and thus
$$
(Kp)_y\sim K e^{-\eta |y|},
$$
insufficient for the argument.  For further discussion, see \cite{Z4, HZ}.
A related, more subtle point arising also in the Lax shock case is
that terms $\sim e^{-\eta|y|}|K_y|$
appearing in $R_y$ bound \eqref{Raoofiest} are also not acceptable, since
the $e^{-\eta|y|}$ term is no help in convolution against an $L^1$ source,
and $|K_y|$ does not have the sign cancellation of $K_y$.
Our replacement with bound \eqref{missing} amounts to showing that,
even though the $y$-derivative no longer introduces cancellation
for these terms, they already possess an additional $t$-derivative 
that leads to cancellation by another route.
}
\er

\subsubsection{Parameter-dependent bounds} \label{ptwise-eps}

 Bounds involving $\eps$-derivatives 
do not appear in \cite{ZH, MaZ3}, where dependence on parameters
was not considered. Such bounds are necessary here in the uniqueness analysis. Recall, Lipschitz continuity with respect to the bifurcation parameter is needed in Section \ref{framework} to obtain an uniqueness result.

\begin{prop}
Under assumptions {\rm (H0)--(H2)}, {\rm (\D)},
the 
residual ${\cal R},$ introduced in Proposition {\rm \ref{greenbounds}} 
satisfies 
\ba \label{Repsbounds}
\partial_\eps \CalR_{y}(x,t;y)&=
\partial_t \partial_\eps  r(x,t;y)+
\partial_t \hat r(x,t;y)+
O(e^{-\eta (|x-y|+ t)})\\
&+\sum_{k=1}^n 
O \left( (t+1)^{-1/2} e^{-\eta x^+} 
+e^{-\eta|x|} \right) 
t^{-1}
e^{-(x-y-a_k^{-} t)^2/Mt} \\
&+
\sum_{a_k^{-} > 0, \, a_j^{-} < 0} 
\chi_{\{ |a_k^{-} t|\ge |y| \}}
O ((t+1)^{-1/2} t^{-1}) 
e^{-(x-z_{jk}^{-})^2/Mt}
e^{-\eta x^+} \\
&+
\sum_{a_k^{-} > 0, \, a_j^{+} > 0} 
\chi_{\{ |a_k^{-} t|\ge |y| \}}
O ((t+1)^{-1/2} t^{-1}) 
e^{-(x-z_{jk}^{+})^2/Mt}
e^{-\eta x^-}, \\
\ea
\ba\label{epsr}
\partial_\eps r(x,t;y)&=
O\big(e^{-\eta|y|}\big) \Big( \sum_{a_k^-<0}
e^{-(x-y-a_k^{-} t)^2/Mt} 
+
\sum_{a_k^+>0}
e^{-(x-y-a_k^{+} t)^2/Mt} \Big), \\
\ea
\ba\label{epshatr}
\hat r(x,t;y)&=
\sum_{k=1}^n 
O \left( (t+1)^{1/2} e^{-\eta x^+} 
+(t+1)e^{-\eta|x|} \right) 
t^{-1}
e^{-(x-y-a_k^{-} t)^2/Mt} \\
&+
\sum_{a_k^{-} > 0, \, a_j^{-} < 0} 
\chi_{\{ |a_k^{-} t|\ge |y| \}}
O ((t+1)^{1/2} t^{-1}) 
e^{-(x-z_{jk}^{-})^2/Mt}
e^{-\eta x^+} \\
&+
\sum_{a_k^{-} > 0, \, a_j^{+} > 0} 
\chi_{\{ |a_k^{-} t|\ge |y| \}}
O ((t+1)^{1/2} t^{-1}) 
e^{-(x-z_{jk}^{+})^2/Mt}
e^{-\eta x^-}. \\
\ea

Moreover, for $y\le 0$, $k\ge 1$, some $M>0$,
\ba \label{Gytepsbounds}
\partial_t^k \partial_\eps \tilde G_{y}(x,t;y)
&=\sum_{k=1}^n 
O \left( t^{-1/2} e^{-\eta x^+} 
+e^{-\eta|x|} \right) \\
&\times
(t+1)^{(k+1)/2} t^{-(k+1)}
e^{-(x-y-a_k^{-} t)^2/Mt} \\
&+
\sum_{a_k^{-} > 0, \, a_j^{-} < 0} 
\chi_{\{ |a_k^{-} t|\ge |y| \}}
O\big((t+1)^{-(k+1)/2} \big)
e^{-(x-z_{jk}^{-})^2/Mt}
e^{-\eta x^+}, \\
&+
\sum_{a_k^{-} > 0, \, a_j^{+}> 0} 
\chi_{\{ |a_k^{-} t|\ge |y| \}}
O\big((t+1)^{-(k+1)/2} \big)
e^{-(x-z_{jk}^{+})^2/Mt}
e^{-\eta x^-} \\
&+
O \big( e^{-\eta|y|}
(t+1)^{-(k+2)/2} 
\big) \\
&\times
\Big( \sum_{a_k^-<0}
e^{-(x-y-a_k^{-} t)^2/Mt} 
+
\sum_{a_k^+>0}
e^{-(x-y-a_k^{+} t)^2/Mt} \Big). \\
\ea

Symmetric bounds hold for $y\ge0$.
%
\end{prop}

\begin{proof}   
The above bounds may be obtained in similar fashion
to those of Proposition \ref{greenbounds},
but require some
additional discussion.
Specifically, recall in the construction of \cite{MaZ3}
that, in the crucial low-frequency
regime, $\partial_y \CalR_\lambda(x,y)$ is expanded about $\lambda=0$ as
a sum of terms of the form
\be\label{genformat}
\partial_y \oint_{\Gamma_1}e^{\lambda t} d(\lambda) \phi(x,\lambda)
\tilde \psi(y,\lambda) \, d\lambda,
\ee
$d$ analytic, where $\phi$, $\psi$ are of form
\ba\label{conjugation}
\phi&= (I+\Theta_\pm(\lambda, x))e^{\mu(\lambda)x}v(\lambda),\\
\psi&= (I+\hat\Theta_\pm(\lambda, y))e^{\nu(\lambda)y}w(\lambda),\\
\ea
with $\mu$, $\nu$, $v$, $w$ analytic, and
\ba\label{Thetabounds}
|\partial_\lambda^j\partial_x^i \Theta_\pm(x,\lambda)|&\le Ce^{-\eta|x|}, \quad x\gtrless 0\\
|\partial_\lambda^j\partial_y^i \hat \Theta_\pm(y,\lambda)|&\le Ce^{-\eta|y|}, \quad y\gtrless 0\\
\ea
for $0\le i\le k$, $0\le j$.
These are then estimated by the Riemann Saddlepoint method as described
above.
Contributions in mid- and high-frequency regimes are of negligible
order $O(e^{-\eta(|x|+|y|+t)})$
and can be handled by simpler estimates \cite{ZH, MaZ3}.

To treat the parameter-dependent case, we have only to observe that
\eqref{genformat}--\eqref{conjugation} still hold, but with 
$d$, $\mu$, $\nu$, $v$, $w$ now depending on $\eps$ as well
as $\lambda$, and $\Theta_\pm$ and $\hat \Theta_\pm$
depending on $\eps$ as well as $(x,\lambda)$ and $(y,\lambda)$,
where the regularity inherited from (H0) is
$C^k$ in both $x$ and $\eps$, $k\ge 2$ as in Assumption (H0) from the introduction of Section \ref{notation}, and $\Theta_\pm$,
$\hat\Theta_\pm$ by the same fixed-point/contraction mapping
construction as in \cite{MaZ3} satisfy
\ba\label{epsThetabounds}
|\partial_\lambda^r\partial_x^i\partial_\lambda^j \Theta_\pm(x,\lambda,\eps)|&\le Ce^{-\eta|x|}, \quad x\gtrless 0\\
|\partial_\lambda^r\partial_y^i\partial_\lambda^j \hat \Theta_\pm(y,\lambda,\eps)|&\le Ce^{-\eta|y|}, \quad y\gtrless 0\\
\ea
for $0\le i,j\le k$, $0\le r$.
Further, just as in the parameter-independent case, growth rates $\mu$, $\nu$
are either {\it fast modes} $|\mu|\ge c_0>0$ 
(resp. $|\nu|\ge c_0>0$) for all $\eps$, 
or {\it slow modes} $\mu=\lambda \mu_0$ (resp. $\nu= \lambda \nu_0$),
where $\mu_0$ (resp. $\nu_0$) is analytic in $\lambda$ and $C^k$ in
$\eps$.
Moreover, either $e^{\mu x}$ and $e^{\nu y}$ are both decaying, $\Re \mu x,
\Re \nu y <0$ for $\Re \lambda >0$, or else $\mu\equiv -\nu$,
with $e^{\mu (x-y)}$ decaying, $\Re \mu (x-y)<0$ for $\Re \lambda >0$.

We treat $\eps$-derivative terms case-by-case, depending on
whether $\mu$ and $\nu$ are fast or slow, and which factor
the $\eps$-derivative falls on in the Leibnitz expansion.
If $\partial_\eps$ falls on $\Theta_\pm$ or $\hat \Theta_\pm$,
for example, then by \eqref{epsThetabounds} the result is no
worse than in the undifferentiated case, and so these terms
may be subsumed in the previously-obtained bounds \eqref{Rybounds}
for $\CalR_y$.  Likewise, derivatives falling on $v$ or $w$ are harmless,
yielding terms of form already estimated in \eqref{Rybounds}.
Derivatives falling on fast modes $e^{\mu x}$, $e^{\nu y}$,
or $e^{\mu (x-y)}$-- necessarily decaying by the case structure
described above-- yield factors $\mu_\eps x e^{\mu x}$, $\nu_\eps ye^{\nu y}$,
or $\mu_\eps (x-y)e^{\mu (x-y)}$ bounded by terms 
$C e^{\mu x/2}$, $Ce^{\nu y/2}$, or $Ce^{\mu (x-y)/2}$ of the same
order as the corresponding undifferentiated term.
Therefore, these contributions, too, may be bounded by 
estimate \eqref{Rybounds}.

It remains only to estimate terms for which the $\eps$-derivative
falls on a slow mode $e^{\mu x}$, $e^{\nu y}$, or $e^{\mu (x-y)}$--
equivalently,
$e^{\lambda\mu_0 x}$, $e^{\lambda\nu_0 y}$, or $e^{\lambda\mu_0 (x-y)}$--
yielding an additional factor of
$\lambda \partial_\eps \mu_0 x\sim \lambda x$, 
$\lambda\partial_\eps\nu_0 y\sim \lambda y$, or 
$\lambda\partial_\eps \mu_0 (x-y)\sim \lambda(x-y)$.
The $\lambda$ factor may either be accounted as a $t$-derivative,
or else as an extra factor of $(t+1)^{-1/2}$ decay arising through
the stationary phase estimate.
The remaining $x$, $y$, or $(x-y)$ factors, when multiplied against
the moving Gaussian bounds obtained by stationary phase, give
a contribution of order $(t+1)$.  For example, $(x-y)$
multiplied against $e^{-(x-y-at)^2/Mt}$ may be bounded
by $(t+1)$ times $e^{-(x-y-at)^2/2Mt}$, since $x-y\sim at + t^{1/2}$
where the Gaussian is of significant size.
Likewise, $x$ or $y$ times the typical term
$\chi_{\{ |a_k^{-} t|\ge |y| \}} e^{-(x-z_{jk}^{-})^2/Mt}$
coming from the product of slow-decaying modes
may be bounded by $(t+1)$ times 
$\chi_{\{ |a_k^{-} t|\ge |y| \}} e^{-(x-z_{jk}^{-})^2/2Mt}$,
since $x\sim z_{jk}+ t^{1/2}$ and $|z_{jk}|$, $|y/a_k^-|\le t$.
Finally, $y$ times the typical term
$e^{-(x-y-at)^2/Mt} e^{-\eta|x|}$ coming from the product of a slow-decaying
$y$-mode and a fast-decaying $x$-mode may be expanded as 
$$
(y-x) e^{-(x-y-at)^2/Mt} e^{-\eta|x|}
+x e^{-(x-y-at)^2/Mt} e^{-\eta|x|}
$$
and each term estimated as above.  Other cases go similarly.

Converting $\lambda$ to $(t+1)^{-1/2}$ decay for terms belonging to
$\partial_t r$ in \eqref{Rybounds}, we obtain an estimate
poorer by factor $(t+1)^{-1/2}(t+1)=(t+1)^{1/2}$ for $\partial_\eps r$
than for $r$ itself.
For all other terms of this type, we convert $\lambda$ into a time-derivative
and group them into $\partial_t \hat r$, where $\hat r$ is bounded
by $(t+1)$ times the remaining terms (besides $\partial_t r$) 
in \eqref{Rybounds}.
Combining these estimates with the previous ones, 
we obtain \eqref{Repsbounds} as claimed.
%

Proceeding similary for the expansion of $\tilde G$, but
converting all $\lambda$ factors to $(1+t)^{-1/2}$ decay, 
we obtain \eqref{Gytepsbounds}, similarly as in the proof of
\eqref{Gytbounds}.
\end{proof}

\br\label{epscont}
\textup{
The effect of $\eps$-derivatives, roughly, is
growth by factor $t^{1/2}$; that is, differentiation
with respect to $\eps$ degrades our decay estimates.
This makes derivative bounds much more delicate than sup-norm
bounds to obtain.
The $t^{1/2}$ growth rate may be understood by
the formal computation, neglecting commutators, 
$$
\partial_\eps e^{L(\eps)t} \sim
t(\partial_\eps L)  e^{L(\eps)t}
\sim t \partial_x e^{L(\eps)t} \sim t^{1/2} e^{L(\eps)t}.
$$
}
\er

\section{Return map construction}\label{returnsection}

 We show in this Section how the question of existence and uniqueness of time-periodic solutions to \eqref{sysepseqn} in a neighborhood of the stationary solution \eqref{profintro} can be formulated in the abstract framework of Section \ref{framework}.  We use the notations set out in Section \ref{notation}, in particular Assumptions (H0)-(H2) and (\D), and functional spaces $X_1, X_2$ and ${\cal B}_1,$ ${\cal B}_2,$ and the short-time bounds of Proposition \ref{masspres}. 
 
Given a family of stationary solutions $\bar u^\eps$ \eqref{prof} of the system \eqref{sysepseqn},
and a family of 
dynamic
solutions $\tilde u^\e$ of the system \eqref{sysepseqn},
define the perturbation variable 
\be\label{pert}
u(\e,x,t):=\tilde u^\eps (x,t) -\bar u^\eps(x),
\ee
satisfying nonlinear perturbation equations
\be\label{nonlin}
u_t- L(\e) u= Q(\e,u)_x,
\qquad
u(\e,x,0)=u_0(\e,x), 
\ee
where the nonlinear term 
\be\label{bounds}
Q(\e,u):= -F(\e,\tilde u^\eps) + F(\e,\bar u^\eps)
+ F_u(\e,\bar u^\eps) \bar u^\e,
\ee
satisfies the pointwise bound
\begin{equation} \label{Qbounds}
 | Q(\e,u)| + | \d_\e Q(\e, u)| \leq C |u|^2,
 \end{equation}
 $C$ being a nondecreasing function of $\| u \|_{L^\infty}.$
 
Decomposing
\be\label{pertdecomp}
u= w_+ \phi^\e_+ + w_- \phi^\e_- + v,
\ee
where (with the notations of Section \ref{projector}) $w_+ \phi^\e_+ + w_- \phi^\e_- = \Pi^\e u \in \Sigma^\e$,
$v:=\tilde \Pi^\e u\in \tilde \Sigma^\e$,
and coordinatizing as $(w, v),$ $w := (w_+, w_-),$ 
we obtain
\ba\label{PHeq2}
\dot w&= (\gamma(\eps) \mbox{Id}\, + \tau(\e) J) w + N (\e, w, v),\\
\dot v&= \tilde L(\e) v + \tilde \Pi^\e \d_x \tilde Q(\e, w, v),\\
\ea
where $J = \left(\begin{array}{cc} 0 & 1 \\ -1 & 0 \end{array}\right),$ and 
 $$ \begin{aligned} N(\e,w,v) & := \Pi^\e \d_x Q(\e, w_+ \phi^\e_+ + w_- \phi^\e_- + v), \\ \tilde Q(\e,w,v) & := Q(\e, w_+ \phi^\e_+ + w_- \phi^\e_- + v). \end{aligned}$$
From \eqref{bounds}, \eqref{Qbounds} and the $\Pi^\e$-bounds of Lemma \ref{efns},
 \begin{equation} \label{tildeQbounds}
  \sup_{x} \, (1 + |x|^2) \, | \tilde Q| \leq C (|w|+\|v\|_{X_1})^2, \qquad |N| \leq C(|w| + \| v \|_{L^q})^2,
  \end{equation}
 and, for $1 \leq q \leq \infty,$
   \begin{equation} \label{tildeQbounds2}
  \| \d_w \tilde Q\|_{{\cal L}(L^2,L^q)} \leq C (|w| + \| v \|_{L^\infty}), \qquad   \| \d_v  \tilde Q\|_{{\cal L}(L^q)} \leq C(|w| + \| v \|_{L^\infty}).  \end{equation}
From \eqref{Qbounds} and Lemma \ref{efns}, for $1 \leq q \leq \infty,$
 \begin{equation} \label{epsQbound}
  \| \d_\e \tilde Q \|_{L^q} \leq C (|w| + \| v \|_{L^{2q}} )^2.
  \end{equation}
 In \eqref{tildeQbounds}, \eqref{tildeQbounds2} and \eqref{epsQbound}, $C$ is a non-decreasing function of $|w| + \| v \|_{L^\infty}.$



Now, truncate $N,$ replacing it with
\be\label{truncN}
\hat N(\e,w,v) := N(\e, w, \hat v(w,v)),\\
\ee
where $\hat v(r,v)$ is the map
\be\label{hatv}
(x,t) \mapsto \left\{ \begin{aligned} \psi\big(C_0 \frac{|w(t)|}{|v(x,t)|_{\R^n}}\big) v(x,t), & \quad \mbox{if $v(x,t) \neq 0,$} \\ 0, & \quad \mbox{if $v(x,t) = 0,$} \end{aligned} \right. 
\ee
$C_0$ being a positive constant and $\psi \in \RR^1$ being a $C^\infty$ ``truncation'' function with
\be\label{psi}
\psi(z)=\begin{cases}
1 & z\ge 1,\\
z & z\le 1/2,\\ 
\end{cases}
\qquad
\psi'(z)=\begin{cases}
0 & z\ge 1,\\
1 & z\le 1/2,\\ 
\end{cases}
\ee
$|\psi| \leq 1,$ $|\psi'| \leq 1.$ With these notations, $|\hat v| \le C|w|$, and thus, by \eqref{tildeQbounds} with $q=\infty$,
\be\label{thetat}
|\hat N|\le C|w|^2,
\ee

The map $(w, v) \in \R^2 \times (L^q \cap L^\infty) \to \hat v \in (L^q \cap L^\infty)$ is $C^\infty$ at any $(w,v)$ such that $w \neq 0,$ and infinitely Fr\'echet differentiable in $v$ for any $(w,v),$ with first partial derivatives
$$ \d_v \hat v = \psi\big(C_0\frac{|w|}{|v|}\big)\Id_{{\cal L}(L^q)} -
C_0 |w| \psi'\big(C_0\frac{|w|}{|v|}\big) \frac{\langle v, \cdot \rangle_{\R^n}}{|v|^3} v,$$
(with the convention $(\d_v \hat v \cdot h)(x) := h(x)$ if $v(x) = 0$), and 
$$ \d_w \hat v = C_0 \frac{\langle w, \cdot \rangle_{\R^2}}{|w|} \frac{v}{|v|}\psi'\big(C_0\frac{|r|}{|v|}\big),$$
(with the convention $\d_w \hat v(x) := 0$ if $v(x) = 0$).
In particular, 
 \begin{equation} \label{hatvbds}  
  \| \d_v \hat v \|_{{\cal L}(L^q)} \leq C, \qquad \| \d_w \hat v \|_{{\cal L}(\R^2,L^q)} \leq C,\end{equation}
  for $1 \leq q \leq \infty,$ uniformly in $w,v,$ for $w \neq 0.$

\begin{lem}\label{hatN}
Under assumptions {\rm (H0)--(H2)}, {\rm (\D)},
as a function from $(w,v)\in \R^2\times 
L^q,$ for any $1\le q\le \infty$,
$\hat N$ is $C^k$, $k\ge 2$ as in Assumption {\rm (H0)} from Section {\rm \ref{notation}}, for $w$ away from 0.
At $w=0$, $\hat N$ is $C^1.$  Moreover, 
\be\label{gradbds}
| \d_w \hat N|_{{\cal L}(\R^2, \R^2)} + \|\d_v \hat N \|_{{\cal L}(L^q, \R^2)}  \leq C  |w|.
\ee
 Besides, $\hat N$ is Lipschitz in $\e,$ with the bound,
\begin{equation} \label{epsNbd} 
 | \d_\e \hat N|  \leq C  |w|^2.
  \end{equation} 
In {\rm \eqref{gradbds}} and {\rm \eqref{epsNbd}}, $C$ is a nondecreasing function of $|w| + \| v \|_{L^\infty}.$
\end{lem}

\begin{proof}
Estimate \eqref{gradbds} is immediate, using the chain rule, \eqref{tildeQbounds} and \eqref{hatvbds}. Estimate \eqref{epsNbd} follows from \eqref{epsQbound} and the $\e$-regularity stated in Lemma \ref{efns}.
\end{proof}

Consider the truncated system
\ba\label{truncPHeq2}
\dot w&= (\gamma(\eps) \mbox{Id } + \tau(\e) J) w + \hat N (\e, w, v),\\
\dot v&= \tilde L(\e) v + \tilde \Pi^\e \d_x \tilde Q(\e, w, v),\\
\ea

\begin{prop}\label{shortbounds}
Under assumptions {\rm (H0)--(H2)}, {\rm (\D)},
for $0\le t\le T$, any fixed $C_1,T>0$, some $C>0$, and $|a|$, $\|b\|_{X_1}$,
$|\eps|$
sufficiently small, system \eqref{truncPHeq2} with initial data
$(w_0,v_0)=(a,b)$ possesses a unique solution
$(w, v)(\e,a,b)\in L^\infty([0,T],\R^2 \times X_1)$ 
that for $a\ne 0$ is $C^{k+1}$ in $t$ and $C^k$ in $(\e,a,b)$, $k\ge 2$ as in Assumption {\rm (H0)} from Section {\rm \ref{notation}},
with respect to the weaker norm $\CalB_1$, and for $a=0$ is
$C^1$ in $t$ and Lipschitz in $(\e,a,b)$ 
with respect to $\CalB_1$, with
\ba\label{truncshort}
C^{-1}|a|&\le |w(t)|\le C|a|,\\
\|v(t)\|_{X_1}& \le C(\|b\|_{X_1}+|a|^2).\\
\ea
In particular, for $\|b\|_{X_1}\le C_1|a|$,
all $0\le t\le T$,
\be\label{truncjust}
\|v(t)\|_{X_1}\le C|w(t)|.
\ee
Besides, for all $0\le t\le T$,
\be\label{Bshorteq}
|\d_{a} w(t)| + \| \d_a v(t) \|_{{\cal B}_1} + \| \d_b w (t) \|_{{\cal L}({\cal B}_1, \R^2)} + \| \d_b v (t) \|_{{\cal L}({\cal B}_1)} \leq C,
\ee
and
\begin{equation} \label{Bshorteq-eps}
| \d_\e w(t)| + \| \d_\e v(t) \|_{{\cal B}_1} \leq C ( | a | + \| b \|_{{\cal B}_1}).
\end{equation} 
\end{prop}

\begin{proof}
Existence and uniqueness follow by a standard Contraction--mapping
argument, using Duhamel's Principle to express \eqref{truncPHeq2} as
\ba\label{fullDuhamel}
w(t)&= W(t) := e^{(\gamma(\eps) + \tau(\e) J)t}a
+\int_0^{t} e^{(\gamma(\eps) + \tau(\e) J)(t-s)} \hat N(\e,w, v)(s)ds,\\
v(t)&= V(t):= e^{\tilde L(\e) t}b
+\int_0^{t} e^{\tilde L(\e)(t-s)} \tilde \Pi^\e \d_x \tilde Q(\e, w, v)(s)ds.\\
\ea

Consider indeed the map ${\mathcal F}_{\e,a,b}: \R^2 \times X_1 \to \R^2 \times X_1,$ which 
maps $(w, v)$ to $(W, V),$ where $W, V$ are defined in \eqref{fullDuhamel}. Estimates \eqref{tildeQbounds}, \eqref{thetat} together with short-time bounds \eqref{XDuhamel1}, \eqref{XDuhamel2} show that given any $T >0,$ for $|\e| \leq \e_0,$ one can find $\rho(\e_0, T) > 0$ such that, for $|a| + \| b \|_{X_1} < \rho,$ ${\mathfrak F}_{\e,a,b}$ maps a closed ball in $L^\infty([0,T], \R^2 \times X_1)$ to itself; the norm in $L^\infty([0,T], \R^2 \times X_1)$ being defined as $\sup_{0 \leq t \leq T} (|w(t)| + \| v(t)\|_{X_1}).$ Estimates \eqref{tildeQbounds2}, \eqref{gradbds} together with short-time bound \eqref{BDuhamel2} show that, provided that $\rho$ is small enough, ${\mathfrak F}_{\e,a,b}$ is moreover a contraction in that same ball of $L^\infty([0,T], \R^2 \times X_1)$ with respect to the weaker norm $\sup_{0\leq t \leq T} (|w(t)| + \| v(t)\|_{{\cal B}_1}),$ whence we obtain existence and uniqueness by
the bounded in strong norm/contractive in weak norm principle
described in Section \ref{framework}.

 Bound \eqref{truncshort}(i) follows readily from \eqref{fullDuhamel}
by Gronwall's inequality, using the ``decoupled''
estimate $|\hat N|\le C|w|^2$ for $\hat N.$
Bound \eqref{truncshort}(ii) is a standard Gronwall bound that 
requires no comment.
Combining \eqref{truncshort}(i)--(ii),
we obtain evidently \eqref{truncjust} for $|a|$ sufficiently small.

Lipschitz regularity with respect to $a,b$ of the fixed point of ${\cal F}_{\e,a,b},$ still denoted $(w, v),$ is a consequence of the uniform bound $| \d_a W (t)| \leq C$ and of the uniform estimate $\| \d_b V(t) \|_{{\cal L}({\cal B}_1)} \leq C$ (itself a direct consequence of \eqref{BDuhamel1}). Estimate \eqref{Bshorteq} follows. 

Similarly, estimate \eqref{Bshorteq-eps} follows from the uniform estimates $| \d_\e W(t)| \leq C |a|$ (a consequence of \eqref{epsNbd}) and $\| \d_\e V(t) \|_{{\cal B}_1} \leq C(|a| + \|b \|_{{\cal B}_1})$ (a consequence of \eqref{epsDuhamel1}, \eqref{epsDuhamel2} and \eqref{epsQbound}). 

Finally, \eqref{truncshort}(i) shows for $a\ne0$ that $w(t)$
remains bounded from zero, whence the solution is $C^k$, $k\ge 2$
in $(\e,a,b)$ and $C^{k+1}$ in $t$, by the corresponding properties
of the righthand side of \eqref{fullDuhamel} away from $w=0$.
\end{proof}

\br\label{aposteriori}
\textup{
Since $\|v\|_{X_1}$ controls $\|v\|_{L^\infty}$,
the bound \eqref{truncjust} allows us to justify a posteriori
the truncation step provided that we can later find a solution
with $\|b\|_{X_1}\le C|a|^2$, and $|a|$, $\|b\|_{X_1}$ sufficiently
small.  That is, such a solution satisfies not only the truncated
system \eqref{truncPHeq2}
but also the original equations \eqref{PHeq2}.
}
\er


 
 Given $(a,b) \in \R^2 \times X_1,$ such that $a \neq 0,$ let $(w,v)$ be the solution of \eqref{truncPHeq2} over $[0,T]$ such that $(w,v)(t=0) = (a,b).$ Estimate \eqref{truncshort}(i) shows that $w$ never vanishes. Considering a double covering in polar coordinates of the equation in the plane $\Sigma^\e,$ 
we coordinatize $w$ as $(r, \theta) \in \R \times \TT,$ obtaining
\ba\label{truncPHeq}
\dot r&= \gamma(\eps)r + \hat N_r(\e, r, \theta, v),\\
\dot \theta&= \tau(\eps) + \hat N_\theta(\e, r, \theta, v),\\
\dot v&= \tilde L(\eps)v + \tilde \Pi^\e \d_x \tilde Q(\e, r e^{i \theta}, r e^{- i \theta}, v),\\
\ea
where 
$$ \hat N_r := \frac{1}{2} \langle e^{-  i \theta} \tilde \phi_+^\e + e^{i \theta} \tilde \phi_-^\e, \hat N\rangle, \qquad \hat N_\theta  := \frac{1}{2 i r} \langle e^{-  i \theta} \tilde \phi_+^\e - e^{i \theta} \tilde \phi_-^\e, \hat N \rangle.
$$ 

Observing, by \eqref{2-1-07}(ii), \eqref{thetat}, \eqref{truncshort}(i) 
and (P)(iii), that
\be\label{thetabd}
\dot \theta= \tau(\eps) + O(|a|)\ne 0,
\ee
we find that all values of $\theta \in [0,2\pi]$ are taken on
for $T$ chosen sufficiently large; say, $\ge 2\pi/\tau(\eps)$.
Thus, in searching for periodic solutions,
we may without loss of generality fix $\theta(0)=0$, or
$w(0)=(a,0)$, $a\in \R^1$, thereby factoring out invariance under
translation in $t$.

With this choice of initial condition, we
obtain after a brief calculation 
 \ba \label{2-1-07} 
  r(t)&= e^{\gamma(\eps)t}a
+\int_0^{t} e^{\gamma(\eps)(t-s)} \hat N_r(\e,r, \theta, v)(s)ds,\\
\theta(t)&= \tau(\eps)t
+\int_0^{t} \hat N_\theta(\e,r, \theta, v)(s)ds, \\ 
v(t)&= e^{\tilde L(\e) t}b
+\int_0^{t} e^{\tilde L(\e)(t-s)} \tilde \Pi^\e \d_x \tilde Q(\e,r e^{i \theta}, r e^{- i \theta}, v)(s)ds
\ea
for all $0 \leq t \leq T$.  Note that we have reduced the dimension
of $a$ from two to one; we will follow this convention for the
rest of the analysis.

\begin{lem}\label{Tfn}
Under assumptions {\rm (H0)--(H2)}, {\rm (\D)},
for $|\e|,$ $|a|$, $\|b\|_{X_1}$ sufficiently small, 
$a \neq 0\in \R^1,$ 
consider the solution $(w, v)$ of \eqref{truncPHeq2} 
on $0\le t\le 4 \pi /\tau(0)$ 
issuing 
from the initial condition $(w,v)(0) = ((a,0),b),$ and let $w = (r e^{i \theta}, r e^{-i \theta}).$ 
Then, there exists a unique smallest $T = T(\e,a,b)>0$ such that  
\be\label{Tperiod}
\theta(\e,a,b, T(\e,a,b))=2\pi.
\ee
Moreover,
the function $(\e,a,b) \mapsto T(\e,a,b)$ is Lipschitz in ${\cal B}_1$ norm,
with 
\be\label{Tprop}
T(\e,0,0)\equiv 2\pi/\tau(\eps).
\ee
\end{lem}

\begin{proof}
 The function $\theta$ is well defined and solves \eqref{2-1-07}(ii). Estimate \eqref{truncshort} shows for $(a,b)=(0,0)$ that
$(r,v)\equiv (0,0)$ for all $t$, hence $\theta(t)=\tau(\eps) t$,
giving \eqref{Tprop}. 
Local existence, uniqueness, and regularity of $T$ satisfying \eqref{Tprop} 
and \eqref{Tperiod} then follows by the Implicit Function Theorem 
applied around $(\e,a,b,T)=(0,0,0, 2\pi/\tau(0))$, using 
the fact that 
\eqref{thetabd} holds
for all $t \leq 4 \pi/\tau(0);$ the function $\theta$ inherits the regularity properties of $w$ described in Proposition \ref{shortbounds}, and $T$ inherits the regularity properties of $\theta.$ 
Finally, \eqref{thetabd} and continuity of $\tau$ show that the 
minimal solution of \eqref{Tperiod} must lie near $2\pi/\tau(0)$,
yielding global uniqueness as well.
\end{proof}

Substituting $t=T(\e,a,b)$ in \eqref{2-1-07},
we may express the Poincar\'e return map
$(\e,a,b)\to (\hat a, \hat b ):=(r,v)(\e,a,b, T(\e,a,b))$
as a discrete dynamical system 
\ba
\label{PHdyn}
\hat a&=R(\e,a,b)a+ N_1(\e,a,b),\\
\hat  b&=  S(\e,a,b)b+ N_2(\e,a,b),\\
\ea
of the form \eqref{dyn} studied in Section \ref{framework},
with $a,\, \eps, \, N_1 \in \RR^1$ and $b\in X_1$, $N_2\in X_2$,
where 
\ba\label{RS}
R(\e,a,b)&:=e^{\gamma(\eps)T(\e,a,b)},\\
S(\e,a,b)&:=e^{\tilde L(\e) T(\e,a,b)},\\
\ea
are primary and transverse linearized
solution operators for one time-step $T(\e,a,b)$ of 
continuous system \eqref{truncPHeq2}, and
\ba\label{PHDuhamel}
N_1(\e,a,b)&:= \int_0^{T(\e,a,b)} e^{\gamma(\eps)(T(\e,a,b)-s)}
\hat N_r(\e,r, \theta, v)(s)ds,\\
N_2(\e,a,b)&:= \int_0^{T(\e,a,b)} e^{\tilde L(\e) (T(\e,a,b)-s)}
\tilde \Pi^\e 
\d_x \tilde Q(\e,r e^{i \theta}, r e^{- i \theta}, v)(s)ds.\\
\ea

Evidently, small-amplitude periodic solutions of \eqref{truncPHeq2} 
with period $T$ close to $T(0,0,0)=2\pi/\tau(0)$ are equivalent
to fixed points of the Poincar\'e return map (equilibria of \eqref{PHdyn}).
Moreover, on the wedge $\{\|b\|_{X_1}\le C_1|a|\}$ 
(see \eqref{truncjust}), these
are equivalent to small-amplitude periodic solutions of the original
(untruncated) system \eqref{sys}, by Remark \ref{aposteriori}.

\begin{lem}\label{returnsetup}
Under assumptions {\rm (H0)--(H2)}
and {\rm (\D)} 
(which are stated in Section {\rm \ref{notation}} and, 
by Lemma {\rm \ref{HP}}, are equivalent to the assumptions 
of Theorem {\rm \ref{newPH}}),
$R$ satisfies \eqref{ass-R} and \eqref{lip-e-R}; $R(0,0,\cdot)$ is differentiable for all $\eps$ sufficiently
small and satisfies \eqref{ass-R-0}; $N_1$ satisfies \eqref{ass-N1} and \eqref{lip-e-N1}; $N_2$ satisfies \eqref{ass-N2} and \eqref{lip-e-N2}.
\end{lem}

\begin{proof} 
{}From Lemma \ref{Tfn}, $R$ satisfies \eqref{ass-R} and \eqref{lip-e-R};
{}from \eqref{Tprop} and (P)(iii) (see section \ref{formule}), we find that $R(\e,0,0)=e^{2\pi \gamma(\eps)/\tau(\eps)}$
is $C^k$, $k\ge 2$ for $\eps$ sufficiently small, satisfying \eqref{ass-R-0}.
 Likewise, the bounds on $N_1$ follow from bounds \eqref{truncshort}(i), \eqref{Bshorteq} and \eqref{gradbds}, 
and ODE bound 
$$|e^{\gamma(\eps)(T(\e,a,b)-s)}|\le C$$ 
on the primary linearized solution operator. Note that we find a control of $|N_1|$ and $|\d_{\e,a,b} N_1|$ by $|a|^2$ and $|a|$ respectively, stronger than \eqref{ass-N1}. 

 Bound \eqref{ass-N2}(i) follows from \eqref{XDuhamel2}, \eqref{tildeQbounds} and \eqref{truncshort}; \eqref{ass-N2}(ii) follows from \eqref{BDuhamel2}, \eqref{tildeQbounds2}, \eqref{Bshorteq} and \eqref{truncshort}, and finally \eqref{lip-e-N2} is a consequence of \eqref{epsDuhamel2}, \eqref{BDuhamel2}, \eqref{Bshorteq-eps}, \eqref{epsQbound} and \eqref{truncshort}. Note that the bounds for $\d_{\e,a,b} N_2$ involve $\| \d_x \tilde Q (T) \|_{{\cal B}_1},$ where $\tilde Q$ is evaluated at $\e,r,\theta,v,$ $(r e^{i \theta}, r e^{-i \theta},v)$ being the solution of \eqref{truncPHeq2} given by Proposition \ref{shortbounds}, and $T = T(\e,a,b);$ by parabolic smoothing, $\| \d_x \tilde Q(T) \|_{{\cal B}_1} \leq C(1 + T^{-1/2}) \| b \|_{{\cal B}_1}.$ 
\end{proof}
 
\br\label{parab}
\textup{
Note that 
we have made use
in Proposition \ref{shortbounds} and Lemma \ref{returnsetup}
 of parabolic smoothing,
reflected by integrability of singularity $t^{-1/2}$ 
in the bound on $e^{\tilde L(\e) t}\tilde \Pi^\e \partial_x$.
Compare the analogous computations in the proof of Theorem
1.3, \cite{TZ} on construction of a center manifold in the
context of model problem I.}
\er

At this point, we have reformulated model problem II in the
abstract framework of Sections \ref{LS} and \ref{bifu},
using only the conservative structure of the equations,
the elementary semigroup estimates of Lemma \ref{masspres},
and direct computation.
It remains to verify Assumptions \ref{B'''} and \ref{B'''bis}:
by Proposition \ref{B-to-A}, essentially linearized stability 
estimates on 
$\sum_{j=0}^\infty S^{j}$,
for which we shall require the detailed pointwise bounds
\eqref{ourdecomp}--\eqref{barbeta} of Section \ref{ptwise}.


\section{Proof of Theorem \ref{newPH}}\label{1dproof}

 
 We start from the perturbation equations \eqref{nonlin},
 where the linear operator $L(\e)$ is assumed to satisfy assumptions (H0)-(H2) and (\D) (which are stated in Section \ref{notation} and are shown in Lemma \ref{HP} to be equivalent to the assumptions of  Theorem \ref{newPH}). 
 
 Consider the Poincar\'e return map system \eqref{PHdyn} associated with the truncated system \eqref{truncPHeq2}, derived from \eqref{nonlin} by a proper choice of coordinatization and truncation. Periodic solutions of \eqref{truncPHeq2} are in one-to-one correspondence with fixed points of \eqref{PHdyn}. 
 Lemma \ref{returnsetup} implies that the basic assumptions of Section \ref{LS} are satisfied by the linear and nonlinear terms in \eqref{PHdyn}. 
 Thus we have reformulated the question of the existence and uniqueness of periodic solutions to the perturbation equations \eqref{nonlin} in the abstract framework of Section \ref{framework}. 
  So far we used only 
 the conservative structure of the equations, the elementary semigroup estimates of Lemma \ref{masspres}, 
and direct computation.

 We show in Section \ref{proof} that the assumptions of Proposition \ref{prop-B} are satisfied by systems \eqref{PHdyn}, thus yielding existence of periodic solutions to the perturbation equations \eqref{nonlin}.
 We show in Section \ref{uniq} that the assumptions of Corollary \ref{group-bif} are satisfied by system \eqref{PHdyn}, thus yielding uniqueness of the periodic solution. 
 Verification of these Assumptions require the detailed pointwise bounds of Sections \ref{ptwise} and \ref{ptwise-eps}. 
In a first step (Section \ref{revisit}), we verify these assumptions for a model linearized operator described in the introduction.

 We use in this proof the notations (in particular, the definitions of the functional spaces $X_1, X_2$ and ${\cal B}_1, {\cal B}_2$) set out in Section \ref{notation}.

\subsection{Key estimates}\label{revisit}

 We return in more detail to the model analysis of Section
\ref{refined}, taking 
\be\label{modeldef}
\tilde G(x,t;y)= K(x,t;y) + 
J(x,t;y),
\ee
where $K(x,t;y):= t^{-1/2}e^{-(x-y-at)^2/4t}$ models the
scattering term $\CalS$ and
$J(x,t;y):=\bar u'(x) \errfn ((-y-at)/2t^{1/2})$ (the function $\errfn$ is defined in \eqref{errfunction}) the excited term $\CalE$
in decomposition \eqref{ourdecomp}.

The aim of this Section is to check Assumptions \ref{B'''} for the model transverse operator $S$ defined by 
 $$ (S f)(x) := \int \tilde G(x,T;y) f(y) \, dy,$$
for some fixed $T > 0,$ where $\tilde G$ is given by \eqref{modeldef}, so that the iterated transverse solution kernel takes form
$$
S^j \delta_y(x)= K(x,jT ;y) +  J(x,jT ;y).
$$
By Proposition \ref{B-to-A}, to prove Assumption \ref{B'''}, it is sufficient to prove that the sequence $\sum_{j=0}^{N}S^j$ is conditionally convergent in ${\cal L}(X_2,{\cal B}_1)$ norm, and that the limit operator $\sum_{j=0}^{\infty}S^j$ is bounded in ${\cal L}({\cal B}_2, {\cal B}_1)$ and ${\cal L}({X}_2,{X}_1)$ norms.


We assume noncharacteristicity, $a<0$. Define the kernels 
 $$K_N:=\sum_{j=1}^{N} K(x,jT ;y), \quad J_N := \sum_{j=1}^{N} J(x,jT ;y),$$
 and associated operators
$$
\CalK_N f(x):=\int_{-\infty}^{+\infty}
K_N(x,j T;y)   f(y)\, dy, \quad \CalJ_N f(x):=\int_{-\infty}^{+\infty}
J_N(x,j T ;y)   f(y)\, dy.
$$
\subsubsection{Scattering estimate ($K$ term)}\label{scat}

Approximating the kernel $K_N$   
by integral
$
\hat K_N(x,y,T):=T^{-1}\int_{T}^{N T} 
K(x,t;y)\, dt,
$
define the ``continuization error'' kernel
$
\theta_N := K_N - \hat K_N,
$
and associated operators
$$
\hat {\cal K}_N f(x) := \int_{-\infty}^{+\infty}
\hat K_N(x,y,T) f(y)\, dy, \qquad 
\Theta_N f (x):=\int_{-\infty}^{+\infty}
\theta_N(x,y,T) f(y)\, dy.
$$

\begin{lem} \label{alg-B}
 There exists $C >0,$ such that, for all $T \geq T_0,$ $y \in \R,$ for $\a =0 ,1$ and $ 0 \leq \b \leq 3,$
  $$
    \| \d_t^\b \d_y^\a K (\cdot, T, y) \|_{L^2} \leq C T^{-(1 + 2\a + 2\b)/4}.
  $$
  In particular, for $T \geq T_0,$
   $$   \int_T^{+\infty} \|  K_{yy}(\cdot, t; y) \|_{L^2} \, dt + \int_T^\infty \|  K_{yt}(\cdot, t; y) \|_{L^2} \, dt \leq C T^{-1/4}. $$
 \end{lem}

\begin{proof} Direct computations, and use of $\| x^\a e^{x^2/ct} \|_{L^2_x} = O(t^{\a/2+1/4}).$ 
\end{proof}

\begin{lem} \label{alg-X}
 For all $x, y \in \R,$ all $T \geq T_0,$ there exists $C(T) > 0$ such that
  \begin{equation} \label{Kfacts-2}
    \quad | K(x,T;y|\le C(T) (1+|x-y|)^{-1}.
  \end{equation}
 There exists $C > 0$ such that, for all $x \in \R, T \geq T_0,$
    \begin{equation} \label{Kfacts-4}
     \int_T^{+\infty} \Big( |  K_{yy}(x, t; y) |,  |  K_{yt}(x, t; y) |  \Big) dt, \leq C (1 + | x - y|)^{-1}.
     \end{equation}
\end{lem}

\begin{proof} The bound \eqref{Kfacts-2} follows from
 $ K(x,T;y) \leq C(T) e^{-(x-y)^2/C T}.$
 Direct computations give the existence of $C >0,$  and $c_0, c_1 > 0,$ independent of $x,y, t,$ such that,  
\begin{equation} \label{0.0} |K_{yy}|\le Ct^{-3/2}e^{-(x-y-at)^2/c_0 t}, \quad |K_{yt}|\le Ct^{-3/2}e^{-(x-y-at)^2/c_1 t},
\end{equation}
 In the light of \eqref{0.0}, to establish \eqref{Kfacts-4}, it is sufficient to prove the uniform bound
$$
\int_T^{+\infty} |x-y|t^{-3/2} e^{-(x-y-at)^2/Ct}\, dt\le C.
$$
Observing that 
$$ \begin{aligned}
|x-y|t^{-3/2} e^{-(x-y-at)^2/Ct} &\le
(|x-y-at|/t^{1/2})t^{-1} e^{-(x-y-at)^2/Ct} \\
&\quad + 
|a|t^{-1/2} e^{-(x-y-at)^2/Ct}\\
&\le 
Ct^{-1/2} e^{-(x-y-at)^2/C_2 t},\\
 \end{aligned} $$
we find in turn that it is sufficient to show
\begin{equation} \label{cle}
\int_T^{+\infty} t^{-1/2} e^{-(x-y-at)^2/Ct}\, dt\le C,
\end{equation} 
which follows by
$$ \begin{aligned}
&\int_T^{+\infty} t^{-1/2} e^{-(x-y-at)^2/Ct}\, dt\\
\qquad &=
\Big(\int_{|x-y-at|\le t/C}+ \int_{|x-y-at|\ge t/C}\Big)
 t^{-1/2} e^{-(x-y-at)^2/Ct}\, dt\\
\qquad &\le 
C\int_{-\infty}^{+\infty} |x-y|^{-1/2} e^{-(x-y-at)^2/C|x-y|}\, dt
+
\int_{0}^{+\infty}t^{-1/2}e^{-t/C^3}dt.
\end{aligned}$$
\end{proof}

\begin{lem}\label{convolution}
For some $C>0$, all $x\in \RR^1$,
$$
\int_{-\infty}^{+\infty}
(1+|x-y|)^{-1}(1+|y|)^{-2} dy \le C(1+|x|)^{-1}.
$$
\end{lem}

\begin{proof}
Let $c(x,y) := (1+|x-y|)^{-1}(1+|y|)^{-2}.$ Dividing into commensurate and incommensurate parts, we have
$$ \begin{aligned}
&\int_{-\infty}^{+\infty} c(x,y) \, dy 
= \int_{|x-y|\le |x|/C} c(x,y) \, dy + \int_{|x-y|\ge |x|/C} c(x,y) \,dy \\
&\qquad  \qquad
\le 
C_2 \int_{x(1-1/C)}^{x(1+1/C)} (1+|x|)^{-2} dy 
+
C_2 \int_{-\infty}^{+\infty} (1+|x|)^{-1}(1+|y|)^{-2} dy \\
&\qquad  \qquad
\le C_3(1+|x|)^{-1}.
\end{aligned} $$ 
\end{proof}
 
\begin{lem} \label{discretization} 
 The sequence $\CalK_N$ is convergent in ${\cal L}({\cal B}_2, {\cal B}_1),$ uniformly in $T \in (T_0, \infty).$ Its limit ${\cal K}_\infty$ is Lipschitz in $T \in (T_0, \infty)$ as an operator in ${\cal L}({\cal B}_2, {\cal B}_1),$ and is bounded in ${\cal L}(X_2,X_1)$ norm, uniformly in $T \in (T_0, \infty).$
 \end{lem}

\begin{proof} That $\hat {\cal K}_N$ belongs to ${\cal L}({\cal B}_2, {\cal B}_1),$ for all $N,$ is a consequence of Lemma \ref{alg-B}. To prove convergence in ${\cal L}({\cal B}_2, {\cal B}_1),$ we use the identity,
 $$K_y= a^{-1}(K_t-K_{yy}),$$
 which implies
\begin{equation} \label{14.0}
a T\partial_y \hat K_N =
K(N T)-  K(T) 
- \int_T^{N T} K_{yy}(s)\, ds,
\end{equation} 
hence
$$\|(\hat \CalK_N-\hat \CalK_{N + p}) \partial_x f\|_{L^2} =
\Big\|\int_{-\infty}^{+\infty}
\partial_y (\hat K_N(x,y)-\hat K_{N + p})f(y)dy\Big\|_{L^2(x)},$$
is bounded by 
$$ \begin{aligned} C\sup_y \Big( &
\|K(x,NT;y)\|_{L^2(x)}  + \|K(x,(N + p)T;y)\|_{L^2(x)}
 \\ & + \int_{NT}^{+\infty} \|K_{yy}(x,t;y)\|_{L^2(x)}\, dt
\Big)
\|f\|_{L^1},\end{aligned}$$
which, by Lemma \ref{alg-B}, is in turn bounded by
$ C (NT)^{-1/4} \|f\|_{L^1}.$
This shows that $\hat \CalK_N (\d_x f)$ is Cauchy in ${\cal B}_1,$ uniformly in $T \geq T_0,$ for $\d_x f \in {\cal B}_2.$ By the uniform boundedness principle, the limit $\hat {\cal K}_\infty$ belongs to ${\cal L}({\cal B}_1, {\cal B}_2).$ By the Fundamental Theorem of Calculus, 
 \begin{equation} \label{theta-N}
  \theta_N = K( N T) + \frac{1}{T} \sum_{1}^N \int_{jT}^{(j+1)T} \int_{t}^{jT} K_{t}(s) \, ds \,dt,
  \end{equation}
hence
\be\label{1trunc}
|\partial_y \theta_N(x,y,T)|\le |K_y(x, N T;y)| + \int_T^{N T} |K_{yt}(x,y)|dt.
\ee
This bound, together with the above Lemmas, implies that $\Theta_N$ belongs to ${\cal L}({\cal B}_2, {\cal B}_1),$ and
$$ \begin{aligned}
\|(\Theta_N & -\Theta_{N + p}) \partial_x f\|_{L^2}
\le
\int_{-\infty}^{+\infty}
\int_{ N T}^{(N+ p )T}
\|K_{yt}(\cdot,y;t)\|_{L^2_x}
|f(y)|\, dy \, dt\\
& +  \int_{-\infty}^{+\infty} \big(\|K_y(\cdot, N T;y)\|_{L^2_x} + \|K_y(\cdot, (N+p)T;y\|_{L^2_x}\big) |f(y)| \, dy \\
& \leq 
C (NT)^{-1/4} \|f\|_{L^1},
\end{aligned}$$
using Lemma \ref{alg-B} again. Hence $\Theta_N (\d_x f)$ is Cauchy in ${\cal B}_1,$ uniformly in $T \geq T_0,$ for all $\d_x f \in {\cal B}_2.$ Combining, 
we have $\CalK_\infty=\hat \CalK_\infty + \Theta_\infty \in {\cal L}({\cal B}_1, {\cal B}_2).$

 The kernel of $\hat {\cal K}_\infty$ is given by the limit $N \to \infty$ in the right-hand side of \eqref{14.0}, where both terms are convergent in the sense of ${\cal L}({\cal B}_2, {\cal B}_1).$ In particular, given $\d_x f \in X_2,$
\be\label{hatKinfty}
   a T \hat {\cal K}_\infty (\d_x f) (x) = - \int K(x,T;y) f(y) dy - \int_T^\infty \int K_{yy}(x,t;y) f(y) \,dy \, dt,
\ee
  and, by Lemma \ref{alg-X},
  $$ \| \hat {\cal K}_\infty (\d_x f) \|_{X_1} \leq C \sup_{x} (1 + |x|) \int (1 + |x - y|)^{-1} |f(y)| \,dy,$$
  which, by Lemma \ref{convolution}, implies that $\hat {\cal K}_\infty \in {\cal L}(X_2, X_1),$ and that $\| \hat {\cal K}_\infty\|_{{\cal L}(X_2, X_1)}$ is uniformly bounded in $(T_0, \infty).$  
  
  Similarly, using \eqref{theta-N}, were both terms in the right-hand side were shown to converge as $N \to \infty,$
  $$ \| \Theta_\infty (\d_x f) \|_{X_1} \leq C \sup_{x} (1 + |x|) \int_T^{\infty} |K_{yt}(x,t;y)| |f(y)| \, dy \,dt,$$
  and Lemmas \ref{alg-X} and \ref{convolution} then imply that $\Theta_\infty \in {\cal L}(X_2, X_1)$ as well, and a uniform bound on $\| \Theta_\infty(T) \|_{{\cal L}(X_2, X_1)}.$ Combining, 
we have $\CalK_\infty$ uniformly bounded in ${\cal L}(X_2, X_1)$ norm. 

To establish Lipschitz regularity of $\hat {\cal K}_\infty,$ we compute 
$$ \| \d_T (T \hat {\cal K}_\infty) (\d_x f) \|_{L^2} \leq  \int (\|K_t(x,T;y)\|_{L^2_x} + \|K_{yy}(x,T;y)\|_{L^2_x}) |f(y)| \, dy,$$
so that, by Lemma \ref{alg-B}, $\hat {\cal K}_\infty$ is Lipschitz in $(T_0, \infty)$ as an operator in ${\cal L}({\cal B}_2, {\cal B}_1),$ with Lipschitz semi-norm bounded by  
$C T_0^{-1}( 1 + \| \hat {\cal K}_\infty \|_{{\cal L}({\cal B}_2, {\cal B}_1)}).$

Lipschitz regularity of $\Theta_\infty$ requires a bit more care,
and an additional observation of general use.
Namely, for a function like $K_y$ for which higher $t$-derivatives
decay successively faster, the total truncation error $\Theta_\infty$
can be partially evaluated as a time-independent function plus an arbitrarily
rapidly converging integral in time, by using successively higher
order numerical quadrature formulae.
 

 The second-order trapezoid rule,
  $$ \begin{aligned}
&(f(1)/2 + f(2)+ \dots+ f(N)+ f(N+1)/2)- \int_1^{N+1} f(t)\,dt\\
&\qquad
=\sum_{j=1}^{N}
\big(
\int_0^1\int_0^s (- s t)f''(j+t)\,dt \,ds
+
\int_0^1\int_1^s (s t+ s-t)f''(j+t)\,dt \,ds
\big)\end{aligned}
$$
gives
\ba\label{express}
\partial_y \theta_N &=
(1/2)(K_y(T) + K_y((N+1) T)\\
&\quad +
\sum_{j=1}^{N+1} \Big(
\int_0^1\int_0^s (-s \tau ) T^2 K_{ytt}(x,\big((j+\tau)T\big); y) \,d\tau \,ds\\
&\quad\quad
+
\int_0^1\int_1^s \big((s \tau + s-\tau)T^2\big) 
K_{ytt}(x,\big((j+\tau)T\big); y) \,d\tau \,ds \Big).
\ea
In \eqref{express}, taking the limit $N \to \infty$ in ${\cal L}({\cal B}_2, {\cal B}_1),$ then differentiating with respect to $T,$ we find that, given $\d_x f \in {\cal B}_2,$
 $$ \begin{aligned}
 \d_T & \Theta_\infty (\d_x f) (x)  = \frac{1}{2} \int K_{yt}(x,T;y) f(y) \,dy   \\
&+
\sum_{j=1}^{N+1} 
\int_0^1 \Big(\int_0^s (-s \tau ) (2 T K_{ytt} + (j + \tau) T^2 K_{yttt}) (x,\big((j+\tau)T\big); y) f(y) \, d\tau \\
&
+
\int_1^s (s \tau + s-\tau) ( (2 T K_{ytt} + (j + \tau) T^2 K_{yttt}) (x,\big((j+\tau)T\big); y) f(y) \, d\tau \Big)  ds.\end{aligned}$$
 Lemma \ref{alg-B} then implies that $\Theta_\infty$ is Lipschtiz in $(T_0, \infty),$ as an operator in ${\cal L}({\cal B}_2, {\cal B}_1),$ with Lipschitz semi-norm bounded by $C T_0^{-1/4}.$
 
\end{proof}



\br\label{negligible}
\textup{
Similarly, we could use Simpson's rule to express
$$
\partial_y \theta_\infty(x,y)=
(1/6)K_{y}(x,T;y)
+
(5/6)K_{y}(x,2T;y)
+
O\Big(\int_T^{+\infty} |K_{ytttt}(x,y)|dt\Big),
$$
obtaining convergence at rate 
$\int \|K_{ytttt}(x,y)\|_{L^2(x)} dt \sim \int t^{-11/4} dt \sim
t^{-7/4}$, and so on.
{\it Together with the favorable $t$-derivative bounds 
noted in \eqref{Gytbounds}, this 
allows us to neglect continuization error for any practical purpose}
of determining convergence or boundedness.
}
\er


\subsubsection{Excited estimate ($J$ term)}\label{excite}
 
%

Approximating the kernel $J_N$   
by integral
$
\hat J_N(x,y,T):=T^{-1}\int_{T}^{NT} 
J(x,t;y)\, dt,
$
define the continuization error kernel 
$\psi_N (x, y):=\hat J_N - J_N,$
and associated operators
$$\hat {\cal J}_N f (x):= \int_{-\infty}^{+\infty}
\hat J_N(x,y,T) f(y)\, dy
, \quad \Psi_N f (x):=\int_{-\infty}^{+\infty}
\psi_N(x,y,T)f(y)\, dy.
$$

\begin{lem} \label{Jdiscretization} 
 The sequence $\CalJ_N$ is convergent in ${\cal L}({X}_2, {\cal B}_1),$ uniformly in $T \in (T_0, \infty).$ Its limit ${\cal J}_\infty$ belongs to ${\cal L}({\cal B}_2, {\cal B}_1),$ is Lipschitz in $T \in (T_0, \infty)$ as an operator in ${\cal L}({\cal B}_2, {\cal B}_1),$ and is bounded in ${\cal L}(X_2,X_1)$ norm, uniformly in $T \in (T_0,\infty).$
\end{lem}
 
\begin{proof} By \eqref{profdecay}, Remark \ref{laxrmk}, 
 $$|J_y(x,t;y)|=|c_2\bar u'(x)K(0,t;y)|\le Ce^{-\eta|x|}|K(0,t;y)|.$$ 
 Given $\d_x f \in X_2,$ the crude bound $|K| \leq C$ then implies 
 $$ \Big\| \int_T^{NT} \int J_y(x;t,y) f(y) \,dy\,dt \Big\|_{L^2_x} \leq C N T |f|_{L^1},$$
 and $\hat {\cal J}_N \in {\cal L}(X_2, {\cal B}_1).$ To prove convergence in ${\cal L}(X_2, {\cal B}_1),$ we use $|f(y)|\le \|\partial_x f\|_{X_2} (1+|y|)^{-2},$ to bound
 $$ \Big\| \int_{NT}^{(N + p)T} \int J_y(x,t;y) f(y) \,dy\,dt \Big\|_{L^2_x},$$
 by
 \begin{equation} \label{001} C \Big( \int_{NT}^{(N+p)T} \int (1 + |y|)^{-2} K(0,t;y) \, dt \Big) \, \|\d_x f\|_{X_2},\end{equation}
 then split $(1 + |y|)^{-2}$ into $(1 + |y|)^{-3/4} (1 + |y|)^{-5/4},$ where the second factor is integrable, to bound \eqref{001} by
 \begin{equation} \label{002} C \Big( \int_{NT}^{(N+p)T} \| (1 + |y|)^{-3/4} K(0,t;y)\|_{L^\infty_y} \, dt \Big) \,  \| \d_x f\|_{X_2}.\end{equation}
 For $t$ large enough, 
\be\label{003}
\sup_y (1 + |y|)^{-3/4} t^{-1/2} e^{-(y - a t)^2/4t} 
\leq C t^{-1/2} (1 + |t|)^{-3/4},
\ee
so that for $N$ large enough, \eqref{002} is controlled by $C (NT)^{-1/4} \| \d_x f \|_{X_2},$ showing convergence of $\hat {\cal J}_N$ in ${\cal L}(X_2, {\cal B}_1),$ uniformly in $T \in (T_0, \infty).$ 
 
 The bound
   $$ \sup_x \, \Big( (1 + |x|) \int_T^\infty \int \bar u'(x) (1 + |y|)^{-2} K(0,t;y) \, dy \, dt \Big) \, < \infty,$$
 is proved similarly and implies that the limit $\hat {\cal J}_\infty$ also belongs to ${\cal L}(X_2, X_1),$ and that $\| {\cal J}_\infty(T) \|_{{\cal L}(X_2, X_1)}$ is uniformly bounded for $T \in (T_0, \infty).$ Note that it is a consequence of \eqref{cle} in Lemma \ref{alg-X} that 
 \begin{equation} \label{boundJ}
  \sup_y \int_{T}^{+\infty} \| J_y(x,t;y) \|_{L^2_x} \, dt  < \infty.
  \end{equation}
 This implies in particular that $\hat {\cal J}_\infty$ also belongs to ${\cal L}({\cal B}_2, {\cal B}_1),$ with norm controlled by \eqref{boundJ}.  
  Finally,
$$\d_T( T \hat {\cal J}_\infty) (\d_x f) = - \int J_y(x,T;y) f(y) \, dy,$$
which implies that in ${\cal L}({\cal B}_2, {\cal B}_1),$ $\hat {\cal J}_\infty$ is Lipschitz with respect in $T \in (T_0, \infty).$
 The treatment of the continuization error $\Psi_N$ and of its limit $\Psi_\infty$
is straightforward; see Remark \ref{negligible}.
The results for $\CalJ_N$ and $\CalJ_\infty$ then follow 
by $\CalJ_N=\hat \CalJ_N+\Psi_N$.
\end{proof}


\subsection{Existence}\label{proof}

 The nonlinear terms in \eqref{PHdyn} satisfy \eqref{ass-N1} and \eqref{ass-N2}. The linear operator $R$  in \eqref{PHdyn} satisfies \eqref{ass-R} and \eqref{ass-R-0}. In the $(r, \theta,v)$ coordinatization, the variable $r$ is scalar. Thus, with the following Lemma, we can apply Proposition \ref{prop-B} to obtain an existence result.

\begin{lem} \label{verif1} 
Under assumptions {\rm (H0)--(H2)}, {\rm (\D)},
the linear operator $S$ defined in {\rm \eqref{RS}}
 satisfies Assumption {\rm \ref{B'''}} and {\rm \ref{B'}}. 
\end{lem} 

\begin{proof} 
Following the model analysis of Section \ref{revisit},
approximate $\Sigma_N:= \sum_{j=1}^{N+1} S^j$ by its continuization
$\hat \Sigma_N$, with associated kernel
$$
\hat \sigma_N(x,y):=\hat \Sigma_N\delta_y(x)= T(\e,a,b)^{-1}
\int_{T(\e,a,b)}^{NT(\e,a,b)} \tilde G(x,t;y) dt,
$$
where $\tilde G$ is the Green function associated with $S$ and  described in Proposition \ref{greenbounds}.
By Remark \ref{negligible}, and the Green function bounds
\eqref{Gytbounds}, the continuization error 
$ \Theta_N:= \Sigma_N-\hat \Sigma_N $
is uniformly bounded as an operator from $X_2\to X_1$ and,
for $\|\cdot\|_{X_2}$ uniformly bounded, is
uniformly convergent from $\CalB_2\to \CalB_1$, with limit uniformly
Lipschitz with respect to $(a,b)$.
Thus, we may discard $\Theta_N$ in what follows, and work directly
with continuization $\hat \Sigma_N$.

Following \eqref{ourdecomp}, decompose
\be\label{Sigmadecomp}
\hat \Sigma_N= \hat \CalE_N+ \hat \CalS_N + \hat \CalR_N,
\ee
with kernels 
$\hat E_N:=\hat \CalE_N \delta_y(x)$,
$\hat S_N:=\hat \CalS_N \delta_y(x)$,
$\hat R_N:=\hat \CalR_N \delta_y(x)$ defined by
\ba\label{kernels}
\hat E_N&:=
T(\e,a,b)^{-1} \int_{T(\e,a,b)}^{NT(\e,a,b)} \CalE (\e,x,t;y) dt,\\
\hat S_N&:=
T(\e,a,b)^{-1} \int_{T(\e,a,b)}^{NT(\e,a,b)} \CalS (\e,x,t;y) dt,\\
\hat R_N&:=
T(\e,a,b)^{-1} \int_{T(\e,a,b)}^{NT(\e,a,b)} \CalR (\e,x,t;y) dt,\\
\ea
where $\CalE$, $\CalS$, $\CalR$ are as in \eqref{ourdecomp}.

The operators $\hat \CalE_N$ and $\hat \CalS_N$ may be estimated
as were $\hat \CalJ_N$ and $\hat \CalK_N$ in 
Sections \ref{scat} and \ref{excite}, respectively.
Indeed, $\hat\CalE_N$ is exactly a superposition
of terms of form $\hat {\cal J}_N$, while 
$\hat \CalS_N$ is a superposition of approximate Gaussian terms
obeying the same estimates used to bound $\hat {\cal K}_N$, so can
be handled in the same way; see \cite{Raoofi, HR, HRZ}
for similar calculations.

We now turn to the residual $\hat \CalR_N.$ By bounds \eqref{Rybounds} and \eqref{hatR},
$\CalR_y$ may be decomposed into terms
$$
\sum_{k=1}^n 
O( e^{-\eta|x|} )
t^{-1}
e^{-(x-y-a_k^{-} t)^2/Mt} 
$$
of order $|J_y|$, $J$ as defined in Section \ref{excite},
terms
$$ \begin{aligned} 
&\sum_{k=1}^n 
O ( (t+1)^{-1/2} e^{-\eta x^+} t^{-1})
e^{-(x-y-a_k^{-} t)^2/Mt} \\
&+
\sum_{a_k^{-} > 0, \, a_j^{-} < 0} 
\chi_{\{ |a_k^{-} t|\ge |y| \}}
O ((t+1)^{-1/2} t^{-1}) 
e^{-(x-z_{jk}^{-})^2/Mt}
e^{-\eta x^+} \\
&+
\sum_{a_k^{-} > 0, \, a_j^{+} > 0} 
\chi_{\{ |a_k^{-} t|\ge |y| \}}
O ((t+1)^{-1/2} t^{-1}) 
e^{-(x-z_{jk}^{+})^2/Mt}
e^{-\eta x^-}, \\
\end{aligned} $$
of order $|K_{yy}|$,
a neglible term of order $e^{-\eta (|x-y|+ t)}$, and
the time-derivative $\partial_t r(x,t;y)$ of terms
$$ \begin{aligned} 
r(x,t;y)&=
O \big( e^{-\eta|y|}
(t+1)^{-1/2} 
\big) \\
&\qquad \times
\Big( \sum_{a_k^-<0}
e^{-(x-y-a_k^{-} t)^2/Mt} 
+
\sum_{a_k^+>0}
e^{-(x-y-a_k^{+} t)^2/Mt} \Big) \\
\end{aligned} $$
of order $|K|$.

Terms of order $|J_y|$ may be estimated exactly as
was $J_y$ in the proof of Lemma \ref{Jdiscretization},
since those arguments depended only on modulus bounds.
Likwise, terms of order $|K_{yy}|$ may be handled by
the modulus-bound arguments used in the proofs of
Lemma \ref{alg-B} and Lemma \ref{discretization}
to bound terms of the same order.
Terms of order $e^{-\eta (|x-y|+ t)}$ may be handled by
a similar argument, using 
$
|\int_{NT}^{+\infty} e^{-\eta (|x-y|+ t)}dt|
\le Ce^{-\eta NT}e^{-\eta|x-y|}.
$
Finally, derivative terms may be estimated by
a cancellation argument like that used to bound $K_y$,
integrating in time to obtain
$ \int_{NT}^{(N+p)T}\partial_t r(x,t;y) dt = r(x,t;y)|_{NT}^{(N+p)T} $,
hence
$$
\Big\| T^{-1}\int_{NT}^{(N + p)T}\partial_t r(x,t;y) dt \Big\|_{L^2(x)}
\le C (NT)^{-1/4},
$$
giving convergence in ${\cal L}(\CalB_2,\CalB_1),$ uniformly in $T.$
Likewise, the limiting kernel
$$
 T^{-1} \int_{T}^{+\infty}\partial_t r(x,t;y) dt = - T^{-1}r(x,T;y)
=O\big(e^{-\eta |x-y|^2}\big)
$$
is clearly bounded by $C(1+|x-y|)^{-1}$, hence the limiting operator
belongs to ${\cal L}(X_2,X_1)$ by the argument used in Lemma
\ref{discretization} to treat $K_y$.

Thus, the sequence $\hat {\cal R}_N$ is convergent in ${\cal L}(X_2, {\cal B}_1),$ with limiting operator bounded in ${\cal L}({\cal B}_2, {\cal B}_1)$ and ${\cal L}(X_2, X_1).$ With Proposition \ref{B-to-A}, Lemmas \ref{discretization} and \ref{Jdiscretization}, and the above remark on the discretization error, we may conclude that $S$ satisfies Assumption \ref{B'''}. 

We now check that Assumption \ref{B'} is satisfied. The results of Lemmas \ref{discretization} and \ref{Jdiscretization} clearly extend to model scattered and excited terms 
\be \label{Keps}
K^\e(x,t;y):= t^{-1/2}e^{-(x-y-a(\eps)t)^2/4t}, \ee \be \label{Jeps}
J^\e(x,t;y):=(\bar u^\eps)'(x) \errfn ((-y-a(\eps)t)/2t^{1/2}),
\ee
now including $\eps$-dependence, as long as $a(\e)$ is bounded away from 0, with uniform convergence with respect to $\e,$ for $\e$ in a neigborhood of the origin. 

Returning to the actual Green function, 
notice that $\e$-dependence of ${\cal E},$ ${\cal S}$ and ${\cal R}$ is similar to \eqref{Keps} and \eqref{Jeps}, with additional multiplicative terms (defined in terms of the spectral elements $a_j^\pm,$ $l_j^\pm,$ $r^j_\pm$) that do not modify the estimates of Section \ref{revisit}. Besides, ${\cal E},$ ${\cal S}$ and ${\cal R}$ depend on $a,b$ only through $T,$ the Lipschitz function given by Lemma \ref{Tfn}. Thus, convergence in ${\cal L}(X_2, {\cal B}_1)$ of $\hat E_N + \hat S_N + \hat R_N$ is uniform with respect to $\e,a,b$ (with $b \in X_1,$ $|\e| + |a | + \| b \|_{X_1}$ small enough). Together with continuity of the Green function with respect to $\e$ (a consequence of bound \eqref{epsDuhamel1}) and of $S$ with respect to $T,$ this implies continuity of the right inverse $(\Id - S)^{-1}$ with respect to $\e,a,b.$ 

 Lipschitz continuity with respect to $b$ of $\hat R_\infty$ is straightforward, 
similarly as in the proof of Lemma \ref{discretization},
since $b$-derivatives of limiting operators
$T^{-1}\int_T^{+\infty}$ pass onto the factor $T^{-1}$ and the 
lower limit of integration only. Thus, Lipschitz continuity 
of $(\Id - S)^{-1}$ with respect to $b \in {\cal B}_1$ in norm ${\cal L}({\cal B}_2, {\cal B}_1)$ follows from Lemmas \ref{discretization} and \ref{Jdiscretization}, and from Lipschitz continuity of $T$ in $b \in {\cal B}_1.$  

 Finally, the uniform bounds \eqref{cont-0} follow from the corresponding 
uniform bounds proved in Lemmas \ref{discretization} and \ref{Jdiscretization}, and the above description of ${\cal R}.$ 
\end{proof}

 Consider now the Poincar\'e return map system \eqref{PHdyn} associated with the truncated equations \eqref{truncPHeq2}. We can apply Proposition \ref{prop-B}, for instance with $\bar \omega \equiv 0:$ there exists a function $a \to \e(a),$ defined for small $|a|,$ such that $\e(0) = 0,$ $\e(\cdot)$ continuous at 0, and periodic maps $r, \theta, v$ such that $(r,\theta,v)$ is a solution of \eqref{truncPHeq2}, with $r(0) = a,$ $v(0) = B(\e(a), a, 0),$ with period $$T^*(a) := T(\e(a), a, B(\e(a), a, 0)),$$ where $T$ is given by Lemma \ref{Tfn} and $B$ by Proposition \ref{prop-B}. From \eqref{Bybd} and \eqref{truncshort}(ii), we find that if $a$ is small enough, \eqref{truncjust} is satisfied. That is, $(r,\theta,v)$ satisfies not only the truncated system \eqref{truncPHeq2} but also the original system \eqref{PHeq2}, and 
 $$u^a := r e^{i \theta} \phi^\e_+ + r e^{-i\theta} \phi^\e_- + v,$$
 is a solution of the perturbation equations \eqref{nonlin} with $\e = \e(a).$ At time $t  = 0,$ 
 $$u^a(\cdot,0) = a (\phi_+^\e + \phi_-^\e) + B(\e(a), a, 0),$$ so that, with \eqref{Bybd},
 $$ C^{-1} a \leq \| u^a(\cdot,0) \|_{X_1} \leq C a.$$
 The map $u^a + \bar u^{\e(a)}$ is a time-periodic solution of the original, unperturbed equations \eqref{sysepseqn}, satisfying 
\eqref{xbd}.

\br\label{complementary}
\textup{
As discussed in Remark \ref{zeromass}, operator 
$(I-S)^{-1}=\sum_{j=0}^{\infty}S^j$ does not preserve zero mass, 
$\int u\, dx=0$, since $\sum_{j=0}^{N} S^j\partial_x$ 
includes terms of form $K(x,NT;y)r_j^\pm l_\pm^k$, $a_j^+>0$
(resp. $a_j^-<0$) with nonzero mass but zero $L^2$-limit as $N\to \infty$.
However, it does appear to preserve $\int \ell^\eps \cdot u\, dx=0$,
for $\ell^\eps$ orthogonal to $r_j^\pm$, $a_j^+>0$ (resp. $a_j^-<0$),
though we shall not pursue this here.
%
That is, it appears to preserve the invariant subspace $\kernel \Pi_0$
complementary to $\kernel L$, where $\Pi_0$ is the generalized
spectral projection defined in Remark \ref{simple}.
Compare the discussion of the finite-dimensional case in 
Remark \ref{restricted}.
%
}
\er

\subsection{Uniqueness} \label{uniq}

 We complete the proof of Theorem \ref{newPH} in this section, as we establish Lipschitz (and later $C^1$) regularity of 
$\eps(\cdot)$ and uniqueness up to translation of solutions,
using still further cancellation in the $\eps$-derivative estimates
to verify Assumption \ref{B'''bis}.

\subsubsection{Reduction to zero speed}\label{zero}

\begin{lem}\label{RHlem}
Under assumptions {\rm (H0)--(H2)}, {\rm (\D)},
let $\bar u^\eps$, 
$\lim_{z\to \pm \infty} \bar u^\eps=u_\pm^\eps$,
be a bounded standing-wave solution of \eqref{sysepseqn},
and $ {\bf u}^a(x-\sigma^a t,t)$
another bounded solution, with ${\bf u}^{a}(x,t)$ time-periodic and 
$\lim_{z\to \pm \infty} {\bf u}^a(z)=u_\pm^\eps$.
Then, $\sigma^a\equiv 0$ if $u_+^\eps\ne u_-^\eps$.
\end{lem}

\begin{proof}
First, observe, by standard parabolic smoothing estimates, 
that both $\partial_x\bar u^\eps$ and $\partial_x{\bf u}^a$
converge to zero as $x\to \pm \infty$.
Integrating $F(\bar u^\eps, \eps)_x-\bar u^\eps_{xx}=0$ over
$[-M,+M]$ for $M>0$ sufficiently large, we thus obtain the 
Rankine--Hugoniot conditions
\be\tag{RH}
F(u_+^\eps, \eps) - F(u_-^\eps, \eps)=0.
\ee
On the other hand, integrating
${\bf u}^a_t -\sigma^a {\bf u}^a_x+ F({\bf u}^a, \eps)_x-{\bf u}^a_{xx}=0$ 
by parts over the rectangle $(x,t)\in [-M,+M]\times [0,T]$, $M>0$
sufficiently large and $T$ the period of ${\bf u}^a$, we obtain 
\be\tag{aRH}
F(u_+^\eps, \eps) - F(u_-^\eps, \eps)=\sigma^a (u_+^\eps-u_-^\eps),
\ee
yielding $\sigma^a=0$ unless $u_+^\eps=u_-^\eps$.
\end{proof}

\subsubsection{Reduction to a wedge}\label{wedge}

Recall that the analysis of Section \ref{returnsection} applies to solutions in a wedge
$\|b\|_{X_1}\le C|a|$, $(a,b):=(w, v)|_{t=0}$, with
period $T=T(\e,a,b)$ determined by the Poincar\'e
return map construction.
We must therefore eliminate the possibility that there may exist
other periodic solutions outside the wedge or with other periods $T$.

Accordingly,
we look for periodic solutions of period $T$ 
of system \eqref{PHeq2},
where $T$ is now considered as an arbitrary parameter.
The Poincar\'e return map for this system is denoted by $(\tilde f,\tilde g);$ periodic solutions of \eqref{PHeq2} are in one-to-one correspondence with solutions of the system 
\ba\label{Piw}
0&= \tilde f(\e,a,b),\\
0&= \tilde g(\e,a,b)=
\big(S(\eps,T)-\Id\big)b + \tilde N_2(\e,a,b,T),\\
\ea 
$(a,b)=(w,v)|_{t=0} \in \R^2 \times X_1$, 
\begin{equation} \label{2.7.07} S(\eps,T):=e^{\tilde L(\eps) T}\tilde \Pi^\e.
\end{equation}

\begin{lem}\label{Atilde}
Under Assumptions {\rm (H0)--(H2)}, {\rm (\D)}, 
on subspace $X_2\subset \CalB_1$ containing $\Range(\tilde N_2|_{X_1})$,
$(\Id -S(\eps,T))$ has a right inverse $(\Id -S(\eps,T))^{-1},$
that is bounded in ${\cal L}(X_2,X_1)$ and ${\cal L}({\cal B}_2, {\cal B}_1),$ 
uniformly in $(\eps,T)$ for $T\in [T_0,T_1]$,
$0<T_0<T_1<\infty,$ and $|\e|$ small enough. 
\end{lem}

\begin{proof}
Identical to the proof of Lemma \ref{verif1} in Section \ref{proof},
which uses only the form of $e^{\tilde L(\eps)t}\tilde \Pi^\e$ 
and the fact that $T$ is uniformly bounded from zero and infinity.
\end{proof}

\begin{lem} \label{lem-1-27} 
Under Assumptions {\rm (H0)--(H2)}, {\rm (\D)}, 
the linear operator $S(\e,T)$ defined in {\rm \eqref{2.7.07}} satisfies Assumption {\rm \ref{group2}}.
\end{lem}

\begin{proof} Consider the family of maps
 $$ \Phi^\e: \R \times (\R^2 \times X_1) \longrightarrow \R^2 \times X_1,$$
 where $ \Phi^\e(c,a,b)$ is defined as
 $$ x \mapsto (a,b(x + c))  + \big( \bar u^\e(x+c) - \bar u^\e(x) \big).$$
 Decomposing $\bar u^\e$ into the sum of $\Pi^\e \bar u^\e \in \R^2$ and $\tilde \Pi^\e \bar u^\e \in X_1,$ we see that $\Phi^\e(c,a,b)$ does belong to  $\R^2 \times X_1.$
 
Let $(a,b)$ be an equilibrium (that is, a solution of \eqref{Piw}) for a given value of $\e;$ let $u^\e= w + v$ be the time-periodic solution of the perturbation equations \eqref{nonlin} issued from $a +  b,$ and let
 $$ (x,t) \mapsto \hat u^\e(x,t) = u^\e(x + c,t) + \big( \bar u^\e(x+c) - \bar u^\e(x) \big).$$
 The map $\hat u^\e + \bar u^\e$ is time-periodic, and a translate of $u^\e + \bar u^\e,$ itself a solution of the autonomous system \eqref{sysepseqn}. Hence $\hat u^\e + \bar u^\e$ is a solution of \eqref{sysepseqn}, and $\hat u^\e$ is a time-periodic solution of the perturbation equations \eqref{nonlin}. The initial datum for $\hat u^\e$ is $\Phi(c,a,b).$ Thus, for any $c,$ $\Phi^\e(c,\cdot,\cdot)$ preserves equilibria.

It is obvious that $\Phi^\e$ is a group action on the set of equilibria. The derivative of $\Phi^\e$ in its first variable at 0 is
 $$ \d_1 \Phi^\e(0,0,0) = (\bar u^\e)'(x),$$
 and it is a consequence of (\D) that $(\bar u^\e)'$ generates $\ker (\Id - S).$ 
 Finally, 
$\Phi^\e(c,a,b)$ is continous in $c$ by density in $X_1$ of 
$C^\infty\cap X_1$.
\end{proof}

\begin{cor} \label{translate} 
Under Assumptions {\rm (H0)--(H2)}, {\rm (\D)}, 
any time-periodic solution $\bar u^\e + u^\e$ of 
\eqref{sysepseqn} with $\|u^\e\|_{X_1}$ sufficiently small, has a translate $(\bar u^\e+\hat u^\e)(x)=(\bar u^\e+u^\e)(x+c)$ such that $(\hat a,\hat b) := \hat u^\e(t=0) \in \Sigma^\e \times \tilde \Sigma^\e$ satisfies $\|\hat b\|_{X_1}\le C|a|$, for any $C > 0.$ 
Thus, in seeking periodic solutions, it is sufficient to
look in the wedge $\|b\|_{X_1}\le C_1|a|$.
\end{cor}

\begin{proof} Just like the quadratic term ${N}_2$ in Lemma \ref{returnsetup}, the quadratic term $\tilde {N}_2$ in \eqref{Piw} satisfies \eqref{ass-N2}(i). Indeed, the change of coordinatization in the $(\phi_+^\e,\phi_-^\e)$ plane does not modify the Duhamel estimate. Besides, Lemma \ref{Atilde} states that the operator $S(\e,T)$ satisfies Assumption \ref{B'''} (up to a slight change in notation; indeed here $S$ depends on the real parameters $\e,T,$ when Assumption \ref{B'''} deals with operators depending on $\e,a,b \in \R \times \R^n \times X_1;$ it is readily seen that this makes no difference in the proof of Proposition \ref{group-inv}). 
Lemma \ref{lem-1-27} states that Assumption \ref{group2} is also satisfied, and we can conclude by Proposition \ref{group-inv}. 

\end{proof} 

\subsubsection{Lipschitz regularity: model computation}\label{lipreg}

Restricting to zero-speed solutions within the wedge, we
now investigate the Lipschitz regularity of $\eps(a)$. We first show how to verify 
Lipschitz regularity in $\e$ for the model problem
\be\label{epsmodel}
\tilde G^\eps(x,t;y)= K^\eps(x,t;y) + 
J^\eps(x,t;y),
\ee
now including $\eps$-dependence, where $K^\e$ is given by \eqref{Keps} and $J^\e$ by \eqref{Jeps}, where $a(\eps)< a_0<0$ and $a(\cdot)$ is $C^{1}$ in $\eps$.
By inspection, this is faithful to the approximation $\tilde G\sim
\CalE^\eps + \CalS^\eps$ obtained by neglecting residual $\CalR^\eps$
in \eqref{ourdecomp}.

By Remark \ref{negligible}, it is sufficient to treat the continuous
approximants $\hat {\cal K}_\infty$, $\hat {\cal J}_\infty$.
The $J$-term may be shown to be Lipschitz in $\eps$
by a modulus-bound computation,
 \ba \label{Jnaive} 
  \|\partial_\eps(T \hat {\cal J}^\eps_\infty) \partial_x f\|_{L^2} & \leq C \Big\| \int J^\e \d_x f (y) \, dy \Big\|_{L^2} \\ & + \Big\| \int_T^\infty \int (\d_\e \d_x \bar u^\e) K^\e(0,t;y) f (y) \, dy \, dt \Big\|_{L^2} \\ & + 
     \Big\| \int_T^\infty \int (\d_x \bar u^\e) \d_\e K^\e(0,t;y) f (y) \, dy \, dt \Big\|_{L^2}.
    \ea
 The first term in \eqref{Jnaive} is bounded by $C \| f \|_{L^\infty};$ the second term is bounded by $C \| \d_x f \|_{X_2},$ due to the uniform bound $\| \d_\e \d_x \bar u^\e \|_{L^2} \leq C,$ and using $\d_\e K^\e = (\d_\e a) t K^\e_y,$ we see that the third term is bounded by $C \| \d_x f \|_{X_2}$ as in \eqref{001}-\eqref{003}. 


The $K$-term by \eqref{hatKinfty} may be expressed as 
\ba \label{cont}
 \d_\e (a T \hat {\cal K}_\infty) \d_x f & = - \int \d_\e K^\e(x,t;y) f(y) \,dy + (\d_\e T) \int K^\e_{yy}(x,t;y) f(y) \, dy \\ & -  \int_T^\infty \int \d_\e K^\e_{yy}(x,t;y) f(y) \, dy \,dt.
 \ea
Lemma \ref{alg-B} implies that the $L^2$ norm of the first two terms in the right-hand side of \eqref{cont} is bounded by $C \| f \|_{L^\infty}.$ The third term may be bounded by a cancellation estimate like
that used to bound $\hat K_\infty$ in the first place, using $\d_\e K^\e_{yy} = (\d_\e a) t K^\e_{yyy},$ 
$K^\e_{yyy}=-a^{-1}K^\e_{yyt} + a^{-1}K^\e_{yyyy}$ and integration by parts
to obtain
$$ \begin{aligned} 
a\int_{T(\e,a,b)}^{+\infty} K^\eps_{yyy}(x,t;y) t \, dt &=
K^\eps_{yy}(x,T;y) T + \int_{T(\e,a,b)}^{+\infty}K^\e_{yy}(x,t;y)\, dt \\
&\qquad +
\int_{T(\e,a,b)}^{+\infty} K^\eps_{yyyy}(x,t;y) t \, dt,
\end{aligned} $$ 
where the last two terms by $|K^\e_{yyyy}|t\sim |K^\e_{yy}|$ are similar
order and uniformly absolutely convergent in $L^2(x)$.

\subsubsection{Lipschitz regularity: full computation}\label{fullcomp}

With these observations, it is straightforward to treat the full problem.

\begin{lem} \label{lipdep}
Under Assumptions {\rm (H0)--(H2)}, {\rm (\D)},  
the linear operator $S$ defined in 
{\rm \eqref{RS}} satisfies Assumption {\rm \ref{B'''bis}}.
 \end{lem}

\begin{proof}
That $S$ satisfies Assumptions \ref{B'''} and \ref{B'} was checked in Lemma \ref{verif1}. It remains thus only to prove that \eqref{lipdisplay}(ii) is satisfied. We check below that $\sum_{j=0}^\infty S^j$ is Lipschitz with respect to $\e$ in ${\cal L}(X_2, {\cal B}_1).$ The dependence of $S$ on $a$ is through $T = T(\e,a,b),$ and so Lipschitz regularity of $\sum_{j=0}^\infty S^j$ with respect to $a$ is proved in the same way. 

By Remark \ref{negligible} and Green function bound \eqref{Gytepsbounds}, 
we can neglect the error made in the approximation of $\Sigma_N:= \sum_{j=1}^{N+1} S^j$ by its continuization
$\hat \Sigma_N$,
with associated kernel
\be\label{ehatsigma}
\hat \sigma(x,y):=\hat \Sigma\delta_y(x)= T(\e,a,b)^{-1}
\int_{T(\e,a,b)}^{\infty} \tilde G(x,t;y) dt,
\ee
where $\tilde G$ is the Green function associated with $S$ and  described in Proposition \ref{greenbounds}.
Decomposing
$$
\hat \Sigma^\eps= \hat \CalE^\eps+ \hat \CalS^\eps + \hat \CalR^\eps,
$$
with kernels 
$\hat E^\eps:=\hat \CalE^\eps \delta_y(x)$,
$\hat S^\eps:=\hat \CalS^\eps \delta_y(x)$,
$\hat R^\eps:=\hat \CalR^\eps \delta_y(x)$ defined by
$$ \begin{aligned} 
\hat E^\eps&:=
T(\e,a,b)^{-1} \int_{T(\e,a,b)}^{\infty} \CalE^\eps (x,t;y) dt,\\
\hat S^\eps&:=
T(\e,a,b)^{-1} \int_{T(\e,a,b)}^{\infty} \CalS^\eps (x,t;y) dt,\\
\hat R^\eps&:=
T(\e,a,b)^{-1} \int_{T(\e,a,b)}^{\infty} \CalR^\eps (x,t;y) dt,\\
\end{aligned} $$ 
where $\CalE^\eps$, $\CalS^\eps$, $\CalR^\eps$ are as in \eqref{ourdecomp},
we may estimate operators $\hat \CalE^\eps$ and $\hat \CalS^\eps$ 
as we did $ K^\eps$ and $J^\eps$ in 
Section \ref{lipreg}.
Thus, we need only estimate the residual term $\hat \CalR^\eps$,
showing uniform boundedness of $\|\partial_\eps \CalR^\eps\|_{L^2(x)}$.

By the computations of Section \ref{proof},
we have
$$
\hat R^\eps_y=
T(\e,a,b)^{-1} \Big(
\int_{T(\e,a,b)}^{\infty} 
\big(\CalR^\eps_y- \partial_t r\big) (x,t;y) dt
-r(x,T(\e,a,b);y) \Big).  \\
$$
We have $\partial_\eps T$ bounded, by Lemma \ref{Tfn}.
Likewise, $\partial_\eps r(x,T;y)$ is bounded in $L^2(x)$
by \eqref{epsr} and boundedness of $|T|$.
Thus, we need only show uniform boundedness in $L^2(x)$ of
$$
\int_{T(\e,a,b)}^{\infty} 
\big(\partial_\eps\CalR^\eps_y- \partial_\eps\partial_t r\big) (x,t;y) dt,
$$
which, by \eqref{Repsbounds}, \eqref{epshatr} may be expressed as
$$ \begin{aligned} 
\int_{T(\e,a,b)}^{\infty} 
\big(\partial_\eps\CalR^\eps_y- &\partial_\eps\partial_t r 
-\partial_t \hat r\big) (x,t;y) dt
+ \lim_{M\to \infty}\hat r(x,t;y)|_{T}^M=\\
&\int_{T(\e,a,b)}^{\infty} 
\big(\partial_\eps\CalR^\eps_y- \partial_\eps\partial_t r 
-\partial_t \hat r\big) (x,t;y) dt
- \hat r(x,T;y)|,\\
\end{aligned} $$ 
where $\|r\|_{L^2(x)}$ is evidently bounded.
Observing that
$$
\big(\partial_\eps\CalR^\eps_y- \partial_\eps\partial_t r 
-\partial_t \hat r\big) (x,t;y)
$$
by \eqref{Repsbounds} may be decomposed into terms
$
\sum_{k=1}^n 
O \left( e^{-\eta|x|} \right) 
t^{-1}
e^{-(x-y-a_k^{-} t)^2/Mt} 
$
of order $|J_y|$ and terms
$$ \begin{aligned} 
&\sum_{k=1}^n 
O \left( (t+1)^{-1/2} e^{-\eta x^+} \right) 
t^{-1}
e^{-(x-y-a_k^{-} t)^2/Mt} \\
&+
\sum_{a_k^{-} > 0, \, a_j^{-} < 0} 
\chi_{\{ |a_k^{-} t|\ge |y| \}}
O ((t+1)^{-1/2} t^{-1}) 
e^{-(x-z_{jk}^{-})^2/Mt}
e^{-\eta x^+} \\
&+
\sum_{a_k^{-} > 0, \, a_j^{+} > 0} 
\chi_{\{ |a_k^{-} t|\ge |y| \}}
O ((t+1)^{-1/2} t^{-1}) 
e^{-(x-z_{jk}^{+})^2/Mt}
e^{-\eta x^-}, \\
\end{aligned} $$ 
of order $|K_{yy}|$, we find by the 
calculations of Sections \ref{excite} and \ref{scat}
that 
$$
\Big\|\int_{T(\e,a,b)}^{\infty} 
\big(\partial_\eps\CalR^\eps_y- \partial_\eps\partial_t r\big) (x,t;y) dt,
\Big\|_{L^2(x)}
$$
is uniformly bounded as well, completing the proof.
\end{proof}


\subsubsection{Conclusion: uniqueness up to group invariance}

 We now sum up the results of Sections \ref{wedge} to \ref{fullcomp}. 
  
 Consider a periodic solution $u^\e$ of the perturbation equations \eqref{nonlin}. By Corollary \ref{translate}, $u^\e$ has a translate $\hat u^\e$ 
whose initial condition $\hat u^\e(t=0) = 
((a\cos \theta_0, a\sin \theta_0), b)
\in \Sigma^\e \times \tilde \Sigma^\e$ satisfies 
 \begin{equation} \label{1-27} \| b\|_{X_1} \leq C |a|. \end{equation} 
 If $a = 0,$ then $\hat u^\e,$ hence $u^\e$ identically vanishes, by uniqueness of the solution of \eqref{nonlin}. If $a \neq 0,$ coordinatize $\hat u^\e$ by $(w,v)$ as Section \ref{returnsection}; $(w,v)$ solves \eqref{PHeq2}. Consider now the associated truncated system \eqref{truncPHeq2}, where $\psi$ is defined by \eqref{psi} with $C_0 : = C.$ ($C$ as in \eqref{1-27}.) It is a consequence of \eqref{1-27} that $(w,v)$ also satisfies system \eqref{truncPHeq2}.
We can thus apply Proposition \ref{shortbounds}; in particular, $w$ never vanishes and we can use polar coordinates $w = (r e^{i \theta}, r e^{-i\theta}),$ 
without loss of generality taking $\theta(0) 
= \theta_0
 = 0$, $w(0)=(a,0)$.
By Lemma \ref{Tfn}, $\theta(T(\e,a,b)) = 2 \pi.$
By Lemmas \ref{lipdep} and \ref{lem-1-27}, the assumptions of 
Corollary \eqref{group-bif} are satisfied, and we can conclude that, 
up to group invariance, $\e = \e_0(a)$ and $b = B(\e_0(a), a, 0),$ 
with the notations of Proposition \ref{lipenough}.


\subsection{$C^1$ dependence}\label{C1}

We have so far carried out the analysis entirely in the Lipschitz
framework of Section \ref{framework}, thereby obtaining Lipschitz
dependence of functions $\eps(\cdot)$, $T(\cdot)$, $u^{(\cdot)}$
on $a$.
However, it is straightforward to improve this regularity to $C^1$.
Focusing on the special solution $\omega\equiv 0$
of Proposition \ref{group-bif}, notice that
$$
\|b\|_{X_1}=
\|B(\eps, a,0)\|_{X_1} \le C|a|^2.
$$
Thus, we may replace linear truncation \eqref{hatv} by a
quadratic order truncation
$$
\hat v:= \psi(C|r|^2/|v|)v
$$
without affecting the validity of truncated equations \eqref{truncPHeq2}.
With this modification, it is readily calculated that $\hat N_r$
and $\hat N_\theta$ become, respectively, $C^2$ and $C^1$.  In particular,
we obtain from this a $C^1$ rather than Lipschitz 
period function $T(\eps, a,b)$ and thereby a $C^1$ return map system
\eqref{PHdyn}.
This improved regularity may be carried throughout the bifurcation analysis
to yield the result.


\section{The multi-d case: Proof of Theorem \ref{multinewPH}}\label{multidcase}

We now briefly describe the extension to the multidimensional case,
model problem III, Section \ref{model3}.
Consider a one-parameter family of standing planar viscous shock 
solutions $\bar u^\eps(x_1)$ of a smoothly-varying family of conservation laws 
\begin{equation}
\label{multisysepseqn2}
u_t =\CalF(\e, u):= \Delta_x u- \sum_{j=1}^d F^j(\e, u)_{x_j},
\qquad u\in \RR^n
\end{equation}
in an infinite cylinder 
$$
\CalC:= \{x:\, (x_1, \tilde x)\in \RR^1\times \Omega\},
\qquad
\tilde x=(x_2, \dots, x_d)
$$
$\Omega\in \RR^{d-1}$ bounded,
with Neumann boundary conditions
$$
\partial u/\partial_{\tilde x} \cdot \nu_\Omega=0
\quad \hbox{\rm for} \quad \tilde x\in \partial \Omega,
$$
(or, in the case that $\Omega$ is rectangular, periodic boundary
conditions), with associated linearized operators
\be\label{multiLdef}
L(\e) :=\partial \CalF/\partial u|_{u=\bar u^\e}
= -\sum_{j=1}^d \partial_{x_j} A^j(x_1, \eps) + \Delta_x,
\ee
$A^j(x,\eps):= F^j_u(\bar u^\eps(x), \eps)$, denoting
$A^j_\pm(\eps):=\lim_{z\to \pm \infty} A^j(z,\eps)=F_u(u_\pm, \eps)$.
Profiles $\bar u^\eps$ satisfy the standing-wave ODE
\be\label{multiode}
u'=F^1(u,\eps)- F^1(u_-,\eps).
\ee

Following \cite{Z4, Z1}, assume:
\medbreak

\quad (H0) \quad  $F^j\in C^{k}$, $k\ge 2$.
\medbreak
\quad (H1)  \quad 
$\sigma (A^1_\pm(\eps))$ real, distinct, and nonzero,
and $\sigma (\sum \xi_j A^j_\pm(\eps))$ real and semisimple for $\xi\in \RR^d$.
\medbreak
\quad (H2)  \quad  Considered as connecting orbits of \eqref{multiode}, 
$\bar u^\eps$ are transverse and unique up to translation,
with dimensions of the stable subpace $S(A^1_+)$ 
and the unstable subspace $U(A^1_-)$ 
summing for each $\eps$ to $n+1$.
\medbreak

As in Remark \ref{laxrmk},
(H2) implies that
$\bar u^\eps$ is of standard Lax type. 

\subsection{Rectangular geometry}\label{rectangular}

For clarity of exposition, we specialize now to the case of
a rectangular cross-section $\Omega$ with periodic boundary conditions,
without loss of generality 
$$
\Omega=\TT^d.
$$
This case closely resembles that of the whole space,
making the analysis particularly transparent.
In particular, we may take the discrete Fourier transform
in transverse directions $\tilde x=(x_2, \dots, x_d)$ to
obtain a family of linearized ordinary differential operators
in $x_1$,
\be\label{family}
L_{\tilde \xi}(\eps):= L_0(\eps)
 -\sum_{j=2}^d i\xi_j A^j(x_1, \eps) -|\tilde \xi|^2,
\ee
indexed by $\tilde \xi\in \ZZ^{d-1}$,
where $L_0(\eps):=  \partial_{x_1}^2 -\partial_{x_1} A^1(x_1, \eps)$ 
is the one-dimensional linearized operator of
\eqref{Ldef}, and $\tilde \xi=(\xi_2, \dots, \xi_d)$ are the
frequency variables associated with coordinates
$\tilde x=(x_2, \dots, x_d)$.

Associated to each $L_\xi(\eps)$, we may define an Evans function 
$
D^\eps(\tilde \xi, \lambda)
$
as in the one-dimensional case, whose zeroes correspond in location
and multiplicity with eigenvalues of $L_\txi$.
Moreover, the eigenvalues of $L$ consist of the union of eigenvalues
of $L_\txi$ for all $\txi\in \ZZ^{d-1}$, with associated eigenfunctions 
$
W(x)=e^{i\tilde \xi\cdot \tilde x}w(x_1),
$
where $w$ is the eigenfunction of $L_\txi$.

To (H0)--(H2) we adjoin the Evans function condition:

\medbreak

\quad (\D) \quad
On a neighborhood of $\ZZ^{d-1}\times \{\Re \lambda \ge 0\}\setminus \{0,0\}$,
the only zeroes of $D(\txi, \cdot)$ are (i) a zero of multiplicity one at
$(\txi,\lambda)=(0,0)$, and (ii) a crossing conjugate pair of zeroes
$\lambda_\pm(\eps)=\gamma(\eps)+i\tau(\eps)$ of some $L_{\txi_*}$,
with $\gamma(0)=0$, $\partial_\eps \gamma(0)>0$, and $\tau(0)\ne 0$.
\medbreak

\begin{lem}\label{multiHP}
Conditions {\rm (H0)--(H2)} and {\rm (\D)} 
 are equivalent to conditions 
{\rm (P)(i)--(iii)} of the introduction (Section {\rm \ref{formule}})
 together with 
$F\in C^k$, $k\ge 2$, 
simplicity and nonvanishing of $\sigma(A^1_\pm(\eps))$, 
semisimplicity of $\sigma (\sum \xi_j A^j_\pm(\eps))$ 
and the Lax condition
$$
\dim S(A^1_+(\eps))+ \dim U(A^1_-(\eps))=n+1,
$$
with $\bar u^\eps$ (linearly) stable
for $\eps <0$ and unstable for $\eps\ge 0$.
\end{lem}

\begin{proof}
Essentially the same as that of Lemma \ref{HP}, but in the final
assertion (here, just a comment) substituting
for the one-dimensional linearized stability analysis of 
\cite{ZH, MaZ3}
a ``one-and-one-half dimensional'' stability analysis 
like that used below to verify Assumption \ref{B'''}; see Remark \ref{setstab}.
\end{proof}

\br\label{crossnumber}
\textup{
Eigenvalues crossing at transverse wave-number $\txi_*=0$
correspond to the  one-dimensional case considered
in Sections \ref{linest}--\eqref{1dproof}, hence
the multidimensional subsumes the one-dimensional analysis.
Such crossings correspond to longitudinal
``galloping'' or ``pulsating'' instabilities.
Eigenvalues crossing at nonzero wave number correspond in the
rectangular geometry to transverse ``cellular'' instabilities
as discussed in \cite{KS}.
}
\er

Introduce Banach spaces $\CalB_1=L^2$,
$\CalB_2=\partial_x L^1\cap L^2$, 
$$
X_1=\{f:\, |f(x)|\le C(1+|x_1|)^{-1}\},
$$
and
$$
X_2=\partial_x \{f:\, |f(x)|\le C(1+|x_1|)^{-2}\} \cap X_1,
$$
equipped with norms
$\|f\|_{\CalB_1}=\|f\|_{L^2}$,
$\|\partial_x f\|_{\CalB_2}=\|f\|_{L^1}+\|\partial_x f\|_{L^2}$,
$$\|f\|_{X_1}=\|(1+|x_1|)f\|_{L^\infty},
\quad \hbox{\rm and } \quad
\|\partial_x f\|_{X_2}=\|(1+|x_1|)^2f\|_{L^\infty}+\|\partial_x f\|_{X_1},
$$
where $\partial_x$ is taken in the sense of distributions.
By inspection, we have that 
$\CalB_2\subset \CalB_1$,
$X_2\subset X_1$,
$X_1\subset \CalB_1$,
$X_2\subset \CalB_2$,
and the closed unit ball in $X_1$ is closed as
a subset of $\CalB_1$.

\subsection{Linearized estimates}\label{multilinest}

\begin{lem}\label{multiefns}
Under Assumptions {\rm (H0)--(H2)}, {\rm (\D)},
associated with eigenvalues $\lambda_\pm(\eps)$ of $L(\e)$
are right and left eigenfunctions $\phi^\eps_\pm=e^{i\txi_*\cdot \txi}
w_\pm(x_1)$ 
and $\tilde \phi^\eps_\pm=e^{i\txi_*\cdot \txi}
\tilde w_\pm(x_1) \in C^k(x,\eps)$, $k\ge 2$ as in Assumption {\rm (H0)} from Section \ref{multidcase},
exponentially decaying in up to $q$ derivatives as
$x_1\to \pm \infty$,  and $L(\e)$-invariant projection
$$
\Pi^\e f:= \sum_{j=\pm} \phi^\eps_j(x)\langle \tilde \phi^\eps_j, f\rangle,$$
onto the total (oscillatory) eigenspace $\Sigma^\eps:=\Span \{\phi^\eps_\pm\}$,
bounded from $L^q$ or $\CalB_2$ to $W^{2,p}\cap X_2$ 
for any $1\le q,p\le \infty$.  Moreover,
$
\phi^\eps_\pm= \partial_x \Phi^\eps_\pm,
$Är
with $\Phi^\eps\in C^{k+1}$ exponentially
decaying in up to $k+1$ derivatives as $x\to \pm \infty$.
\end{lem}

\begin{proof}
Essentially identical to that of Lemma \ref{efns}.
\end{proof}

Defining $\tilde \Pi^\eps,$ $\tilde \Sigma^\eps$ $\tilde L(\eps)$ as in Section \ref{projector}, denote by
$G$ as in \eqref{kernel} (where $x$ is now multi-dimensional)  
the Green kernel associated with the linearized
solution operator $e^{Lt}$ of the linearized evolution equations
$u_t=L(\eps) u$, and by $\tilde G$ as in \eqref{transkernel} 
the Green kernel associated with the transverse linearized
solution operator $e^{\tilde L(\eps) t}\tilde \Pi^\e$.
Evidently, $G = \CalO +  \tilde G$, where $\CalO$ is defined as in \eqref{O}. 


\begin{lem}[Short-time estimates]\label{multimasspres}
Under Assumptions {\rm (H0)--(H2)}, {\rm (\D)},
for $0\le t\le T$, any fixed $T>0$, and some $C=C(T)$,
\be\label{multiBDuhamel0}
 \| e^{\tilde L(\e)} \tilde \Pi f \|_{{\cal B}_1} \leq \| f \|_{{\cal B}_1}, \quad \| e^{\tilde L(\e)} \tilde \Pi f \|_{X_1} \leq \| f \|_{X_1},
 \ee
  \ba\label{multiBDuhamel}
\|e^{L(\eps) t}\partial_x f \|_{\CalB_2},\,
\|e^{\tilde L(\eps) t}\tilde \Pi^\e \partial_x f \|_{\CalB_2}
&\le Ct^{-1/2} \|f\|_{L^1\cap L^2}.\\
\ea
\ba\label{multiXDuhamel}
\|e^{L(\eps) t}\partial_x f \|_{X_2},\,
\|e^{\tilde L(\eps) t}\tilde \Pi^\e \partial_x f \|_{X_2}
&\le Ct^{-1/2} \|(1+|x_1|)^2f\|_{L^\infty},\\
\ea
and
  \begin{eqnarray}
  \| \d_\e (e^{\tilde L(\e) t} \tilde \Pi^\e) f \|_{{\cal B}_1} & \leq &  C \| f \|_{{\cal B}_1}, \label{multiepsDuhamel1} \\
  \| \d_\e (e^{\tilde L(\e) t} \tilde \Pi^\e) \d_x f \|_{{\cal B}_2} & \leq &  C ( \| f \|_{L^1} + \| f \|_{{\cal B}_1}). \label{multiepsDuhamel2}
\end{eqnarray}
\end{lem}

\begin{proof} Essentially identical to that of Lemma \ref{masspres}.
\end{proof}

Note that the Fourier transform completely decouples the linearized
problem in the sense that the Fourier-transformed operators
$$
\tilde G_\txi(x_1,t;y_1)= \CalF G \CalF^{-1}
$$ 
acting in frequency space,
$\CalF$ denoting Fourier transform, have simple form
\be\label{FTmultitranskernel}
\tilde G_\txi(x,t;y)= 
\begin{cases}
e^{L_\txi(\eps) t} \delta_{y_1}(x_1)& \txi\ne \txi_*,\\
e^{L_\txi(\eps) t}\tilde \Pi_{\txi_*} \delta_{y_1}(x_1)& \txi = \txi_*,\\
\end{cases}
\ee
$\tilde \Pi_\txi:= \Id - \Pi_\txi$, where
$$
\Pi_\txi f:= \sum_{j=\pm} w^\eps_j(x)\langle \tilde w^\eps_j, f\rangle
$$
is the one-dimensional projection onto the oscillatory subspace
of $L_{\txi_*}$.

\begin{prop}[Global bounds]\label{multilin}
Under Assumptions {\rm (H0)--(H2)}, {\rm (\D)}: (i) $\tilde G_0$ satisfies
the pointwise bounds stated in Proposition \ref{greenbounds}.
(ii) For $\xi\ne 0$, the Fourier transformed solution operators 
$e^{L_\txi(\eps) t} $  for $\txi\ne \txi_*$,
$e^{L_\txi(\eps) t}\tilde \Pi_{\txi_*}$ for $\txi = \txi_*$
satisfy for some $C,\, \eta>0$, the exponential bounds
\ba\label{expsemi}
\|(1+|x_1|)e^{L_\txi(\eps) t}f\|_{L^\infty(x_1)}
&\le Ce^{-\eta t} \|(1+|x_1|)f\|_{L^\infty(x_1)},\\
\|(1+|x_1|)e^{L_\txi(\eps) t}\tilde \Pi_{\txi_*}f\|_{L^\infty(x_1)}
&\le Ce^{-\eta t} \|(1+|x_1|)f\|_{L^\infty(x_1)}.
\ea
\end{prop}

\begin{proof}
Assertion (i) follows immediately from the observation that $L_0$ 
is exactly the one-dimensional linearized operator studied in 
Section \ref{linest}, together with the fact that the bounds 
on the restriction to $\tilde \Sigma$
in case $\txi_*=0$ are the same as the bounds on the full solution operator
in the case $\txi\ne 0$ of linearized one-dimensional stability.

Assertion (ii) follows for $|\txi|$ large by
standard semigroup/asymptotic ODE estimates.
For $|\txi|$ bounded but nonzero, it follows by standard
semigroup estimates, together with the assumption
that there are no eigenvalues of $L_\txi$ other than the crossing
pair at $\txi=\txi_*$ and the computation as in \eqref{disp} for
the one-dimensional case
that $\sigma_\ess (L_\txi)\subset \{\lambda: \, \Re \lambda
\le -2\eta |\txi|^2\}$ for some $\eta>0$.
\end{proof}

\br\label{notzero}
\textup{
The argument for (ii) of course fails at $\txi=0$, due to the
lack of spectral gap.
}
\er

\subsection{Proof of Theorem \ref{multinewPH}}\label{multimain}

We restrict for definiteness to the case of periodic 
boundary conditions. The proof in the Neumann case is
essentially identical.

Using the bounds of Lemmas \ref{multiefns} and \ref{multimasspres},
we may carry out the polar coordinate and
Poincar\'e return map construction of Section \ref{returnsection}
essentially unchanged to obtain a multidimensional version
of Lemma \ref{returnsetup} (statement unchanged), reducing
the problem to the abstract form studied in Section \ref{framework}.
Thus, it is sufficient to establish the linearized 
estimates of Assumptions \ref{B'''} and \ref{B'''bis}.
But, these follow easily using the one-dimensional bounds of
Proposition \ref{multilin}(i) and the calculations of the one-dimensional
case together with the $\xi\ne0$ bounds of 
Proposition \ref{multilin}(ii).

For example, uniform boundedness of $\int_{T}^{+\infty}e^{\tilde L(\eps)t}dt$
from $X_2\to X_1$ follows by
\ba\label{multikey}
\Big\|\int_{T}^{+\infty}e^{\tilde L(\eps)t}f\, dt\Big\|_{X_1}&:=
\Big\|(1+|x_1|)\int_{T}^{+\infty}e^{\tilde L(\eps)t}f\, dt\Big\|_{L^\infty(x_1,\tilde x)}\\
&\le
\sum_{\txi}
\Big\|(1+|x_1|)\int_{T}^{+\infty}
e^{\tilde L_\txi(\eps)t}\hat f\, dt\Big\|_{L^\infty(x_1) }\\
\ea
and the observation, by Proposition \ref{multilin}(ii), 
Hausdorff--Young's inequality, and the fact that $L^\infty$
controls $L^1$ on bounded domains, that
\ba\label{nez}
\sum_{\txi\ne 0}
\Big\|(1+|x_1|)\int_{T}^{+\infty}
e^{\tilde L_\txi(\eps)t}\hat f\, dt\Big\|_{L^\infty(x_1) }
&\le
\sum_{\txi\ne 0} \Big(\int_T^{+\infty}Ce^{-\eta |\txi|^2 t}
\, dt\Big)\\
&\qquad \times
\|(1+|x_1|)\hat f(\cdot, \txi)\|_{L^\infty(x_1)}\\
&\le
C_2 \sup_{\txi} \|(1+|x_1|)\hat f(\cdot, \txi)\|_{L^\infty(x_1)}\\
&\le
C_2 \| (1+|x_1|)f\|_{L^1(\tilde x; L^\infty(x_1))}\\
&\le
C_3 \| (1+|x_1|)f\|_{L^\infty(x)}\\
&=
C_3 \|f\|_{X_1},
\ea
together with the computation, using the one-dimensional estimates 
carried out in Section \ref{1dproof}, and
denoting $f=\partial_{x_1} F_1 + \partial_{\tilde x}\tilde F$,
so that $\hat f= \partial_{x_1} \hat F_1$, of
\ba\label{multione}
\Big\|(1+|x_1|)\int_{T}^{+\infty}
e^{\tilde L_0(\eps)t}\hat f\, dt\Big\|_{L^\infty(x_1) }
&=
\Big\|(1+|x_1|)\int_{T}^{+\infty}
e^{\tilde L_0(\eps)t}\partial_{x_1}\hat F_1\, dt\Big\|_{L^\infty(x_1) }\\
&\le
C\|(1+|x_1|)^2 \hat F(\cdot, 0)\|_{L^\infty(x_1)}\\
&\le 
C \sup_{\txi} \|(1+|x_1|)^2\hat F(\cdot, \txi)\|_{L^\infty(x_1)}\\
&\le
C \| (1+|x_1|)^2 F\|_{L^1(\tilde x; L^\infty(x_1))}\\
&\le
C_2 \| (1+|x_1|)^2 F\|_{L^\infty(x)}\\
&:=
C_2 \|f\|_{X_2}.
\ea
Other computations follow similarly.

\br\label{setstab}
\textup{
We point out that we have in passing set up a framework
suitable for the linearixe stability analysis of flow in a duct, 
a problem of interest in its own right and somewhat different
from either the one-dimensional case considered in \cite{ZH, MaZ3}
or the multi-dimensional case considered in \cite{ZS, Z1}.
It would be very interesting to try to carry out a full nonlinear stability
analysis by this technique.
}
\er

\subsection{General cross-sectional geometry} \label{gencross}

We may treat general cross-sections $\Omega$ by separation of
variables, decomposing
$$
L(\eps, x_1, \partial_x)=L_0(\eps, x_1, \partial_{x_1})
+ M(\eps, x_1, \partial_{\tilde x}),
$$
where 
$$
M:= -\sum_{j=2}^d A^j(x_1, \eps)\partial_{x_j} + \Delta_{\tilde x},
$$
and expanding the perturbation $u$ in eigenfunctions 
$w_j(\eps, x_1)$ of $M(\eps, x_1, \partial_{\tilde x})$
on domain $\Omega$ as 
$$
u(x,t)=\sum_{k=0}^{+\infty} \alpha_k(x_1) w_k(\eps, x_1)(\tilde x),
$$
$\alpha_k=:\hat u(k)$, to recover a decoupled system
$$ \begin{aligned} 
\partial_t \alpha_k&=L_k\alpha_k\\
&:= (L_0 + \nu_k) \alpha_k 
+ O(|\partial_{x_1}w_k|+
|\partial_{x_1}^2w_k|)\alpha_k + 
O(|\partial_{x_1}w_k|) \partial_{\tilde x}\alpha_k,
\end{aligned} $$ 
where $\nu_k=\nu_k(x_1)$ are the eigenvalues of $M$ associated with $w_k$,
with both $\nu_k$ and $w_k$ converging exponentially in $x_1$ to
limits as $x_1\to \pm \infty$, that is, for which $L_k$ is of the
same basic form as $L_{\txi}$ in Section \ref{multimain}.
Thus, we can carry out the entire Evans function construction
of Section \ref{rectangular}, with $k$ replacing $\txi$.

More, since $w_0\equiv 1$ for Neumann boundary conditions, we have
again that $L_k=L_0$ for $k=0$, justifying our notation.
Likewise, augmenting P(i) with the assumption that $\nu_k\le -\nu_*<0$
for all $k\ge 1$, we recover a spectral gap for all $L_k$, $k\ge 1$.
{}From these two observations, we obtain Proposition \ref{multilin} exactly
as before.
Lemma \ref{multimasspres} holds also, indeed was independent of
$\Omega$.
Thus, we may carry out the entire argument of Theorem \ref{multinewPH}
essentially unchanged, provided that we can establish the
generalized Hausdorff--Young inequalities
$$ \begin{aligned} 
\|\hat u\|_{\ell^\infty(k, x_1)}&\le C\| u\|_{L^\infty(x_1, L^1(\tilde x))}\\
\| u\|_{L^\infty(x)}&\le
C\|\hat u\|_{L^\infty(x_1, \ell^1(k, x_1))},
\end{aligned} $$ 
which follow, for example, if one can establish 
uniform supremum bounds on the right and left eigenfunctions of $M$
by asymptotic eigenvalue--eigenfunction estimates as $k\to \infty$.

Alternatively, we could perform a more general analysis allowing
viscosities with cross-derivatives, 
and also avoiding the need for asymptotic spectral analysis, 
by expanding $u$ instead in
eigenfunctions of $\Delta_{\tilde x}$ (or, more generally, the
second-order elliptic operator in derivatives of $\tilde x$ appearing
in operator $L$), and noting that for Neumann boundary conditions
the first eigenfunction $w_0$ is still $\equiv 1$.
Thus, though different wave numbers no longer completely decouple,
we still obtain a cascaded system in which the zero wave-number
equation is again just $\partial_t w_0=L_0 w_0$ as in the one-dimensional
case, and we can perform an analysis as before, but dividing only
into the two blocks $k=0$ and $k\ge 1$, treating coupling terms in
the $k\ge 1$ part as source terms driving the equation.

Estimating the one-dimensional part as before in \eqref{multione}
and for the $k\ge 1$ part substituting for \eqref{multikey}--\eqref{nez} the 
Sobolev/Parseval estimates
$$ \begin{aligned} 
\Big\|&\int_{T}^{+\infty}e^{\tilde L(\eps)t}\Pi_{k\ge 1}f\, dt\Big\|_{X_1}:=
\Big\|(1+|x_1|)\int_{T}^{+\infty}e^{\tilde L(\eps)t}f\, dt\Big\|_{L^\infty(x_1,\tilde x)}\\
&\le
\Big\|(1+|x_1|)\int_{T}^{+\infty}e^{\tilde L(\eps)t}\Pi_{k\ge 1}f\, dt
\Big\|_{L^\infty(x_1,H^{d-1}(\tilde x))}\\
&\le
\sqrt{
\sum_{k\ge 1}
\Big\|
\int_{T}^{+\infty}
(1+|x_1|)
\| (1+|\nu_k|)^{(d-1)/2}
e^{\tilde L_k(\eps)t}\hat f\|_{\ell^2(k)}
\, dt\Big\|_{L^\infty(x_1) }^{2}}\\
&\le
\|(1+|x_1|)\hat f(\cdot, \txi)\|_{L^\infty(x_1),\ell^2(k)}
\sqrt{
\sum_{k\ge 1} \Big(\int_T^{+\infty}(1+t^{-(d-1)/2})Ce^{-\eta |\nu_k| t}
\, dt\Big)^2}\\
&\le
C_2 
\|(1+|x_1|)\hat f(\cdot, \txi)\|_{L^\infty(x_1),\ell^2(k)}\\
&\le
C_2 \| (1+|x_1|)f\|_{L^2(\tilde x; L^\infty(x_1))}\\
&\le
C_3 \| (1+|x_1|)f\|_{L^\infty(x)}
=
C_3 \|f\|_{X_1},
\end{aligned} $$ 
we would then obtain the result without recourse to asymptotic eigenvalue
estimates.
We shall not pursue these issues further here.


\end{document}